\documentclass{amsart}

 \usepackage[dvips]{graphicx}
 \newtheorem{theo}{Theorem}[section] 
 \newtheorem{lemm}[theo]{Lemma}

\newtheorem{coro}[theo]{Corollary}
\newtheorem{prop}[theo]{Proposition}
\newtheorem{rema}[theo]{Remark}

\newcommand\NN{{\mathbb N}}
\newcommand\RR{{{\mathbb R}}}

\newcommand\R{{{\mathbb R}}}
\def\SS {\mathbb{S}}

\def\la{\langle}
\def\ra{\rangle}

\newcommand\Supp{\mathop{\rm supp}}

\let\a=\alpha
\let\b=\beta

\let\Lam=\Lambda

\newcommand\cD{\mathcal D}
\newcommand\cE{{\mathcal E}}
\newcommand\cF{{\mathcal F}}
\newcommand\cH{{\mathcal H}}

\newcommand\cS{{\mathcal S}}
\newcommand\cM{{\mathcal M}}
\newcommand\cT{{\mathcal T}}

\newcommand\pa{\partial}

\def\vk{\mbox{\boldmath $k$}}

\renewcommand{\theequation}{\thesection.\arabic{equation}}
\begin{document}

\title[The non cutoff Boltzmann equation]
{Regularizing effect and local existence\\
for non-cutoff Boltzmann equation}
\author{R. Alexandre}
\address{IRENAV Research Institute, French Naval Academy
Brest-Lanv\'eoc 29290, France}
\email{radjesvarane.alexandre@ecole-navale.fr}
\author{Y. Morimoto }
\address{Graduate School of Human and Environmental Studies,
Kyoto University, Kyoto, 606-8501, Japan}
\email{morimoto@math.h.kyoto-u.ac.jp}
\author{S. Ukai}
\address{17-26 Iwasaki-cho, Hodogaya-ku, Yokohama 240-0015, Japan}
\email{ukai@kurims.kyoto-u.ac.jp}
\author{C.-J. Xu}
\address{Universit\'e de Rouen, UMR 6085-CNRS, Math\'ematiques
Avenue de l'Universit\'e,\,\, BP.12, 76801 \newline \hskip0.458cm
Saint 
Etienne 
du Rouvray,
France,
\newline \hskip0.458cm 
and  School of Mathematics, Wuhan University 430072,
Wuhan, China} \email{Chao-Jiang.Xu\@@univ-rouen.fr}
\author{T. Yang}
\address{Department of mathematics, City University of Hong Kong,
Hong Kong, P.R. China}
\email{matyang@cityu.edu.hk}


\date{}


\begin{abstract}
The Boltzmann equation without Grad's angular cutoff assumption is
believed to have regularizing effect on the solution because of the
non-integrable angular singularity of  the cross-section. However,
even though so far this has been justified satisfactorily for the
spatially homogeneous Boltzmann equation, it is still basically
unsolved for the spatially inhomogeneous  Boltzmann equation. In
this paper, by sharpening the coercivity and upper bound estimates
for the collision operator, establishing the hypo-ellipticity of the
Boltzmann operator based on a generalized version of the uncertainty
principle, and  analyzing the commutators between the collision
operator and some weighted pseudo-differential operators,  we prove
the regularizing effect in all (time, space and velocity) variables
on solutions when some mild regularity is imposed on these
solutions. For completeness, we also show that when the initial data
has this mild regularity and Maxwellian type decay in velocity
variable, there exists a unique local solution with the same
regularity, so that  this solution acquires 
the $C^\infty$ regularity
for positive time.
\end{abstract}

\maketitle

\section{Introduction}\label{s1}

Consider the Boltzmann
equation,
\begin{equation}\label{1.1}
f_t+v\cdot\nabla_x f=Q(f, f),
\end{equation}
where $f= f(t,x,v)$ is the density distribution function
of particles  with position
$x\in \RR^3$ and velocity $v\in \RR^3$ at time $t$.
The right hand side of (\ref{1.1}) is given by the
Boltzmann bilinear collision operator
\[
Q(g, f)=\int_{\RR^3}\int_{\mathbb S^{2}}B\left({v-v_*},\sigma
\right)
 \left\{g(v'_*) f(v')-g(v_*)f(v)\right\}d\sigma dv_*\,,
\]
which is well-defined for
 suitable functions $f$ and $g$ specified later. Notice that the
collision operator $Q(\cdot\,,\,\cdot)$ acts only on the velocity
variable $v\in\RR^3$. In the following discussion, we will use the
$\sigma-$representation, that is, for $\sigma\in \mathbb S^{2}$,
$$
v'=\frac{v+v_*}{2}+\frac{|v-v_*|}{2}\sigma,\,\,\, v'_*
=\frac{v+v_*}{2}-\frac{|v-v_*|}{2}\sigma,\,
$$
which give the relations between the post and
pre collisional velocities.

It is well known that the Boltzmann equation is a fundamental
equation in statistical physics. For the mathematical theories on
this equation, we refer the readers to
\cite{Ce88,Cercignani-Illner-Pulvirenti, D-L,Grad, villani2}, and
the references therein also for the physics background.

In addition to the special bilinear structure of the collision
operator, the cross-section $B(v-v_*,\sigma)$ varies with different
physical assumptions on the particle interactions and it plays an
important role in the well-posedness theory for the Boltzmann equation.
In fact, except for the hard sphere model, for most of the other molecular
interaction potentials such as  the inverse power laws, the cross
section  $B(v-v_*,\sigma)$ has a non-integrable angular singularity.
For example, if the interaction potential obeys the inverse power law
 $r^{-(p-1)}$ for $2<
p<\infty$, where $r$ denotes the distance between
 two interacting  molecules,
the cross-section behaves like
\[
B(|v-v_*|, \cos \theta)\sim |v-v_*|^\gamma \theta^{-2-2s},\,\,\,\,\,
\cos \theta=\big<\frac{v-v_*}{|v-v_*|} \, ,\,\sigma\big>,\,\,\,
0\leq\theta\leq\frac{\pi}{2},
\]
with
\[
-3<\gamma=\frac{p-5}{p-1}<1,\,\,\, \,\,\,\,\,\,0<s=\frac{1}{p-1}<1.
\]
As usual, the hard and soft potentials correspond to $2<p<5$ and
$p>5$ respectively, and the Maxwellian potential corresponds to
$p=5$. The fact that the singularity $\theta^{-2-2s}$ is not
integrable on the unit sphere leads to the conjecture that the nonlinear
collision operator should  behave like a fractional Laplacian in the variable
$v$. That is,
$$
Q(f,\, f)\, \sim\,  -{(-\Delta_v )^{s}f} +\text{lower order terms}.
$$
Indeed, consider the Kolmogorov type equation
$$
f_t+v\,\cdot\,\nabla_x f=-{(-\Delta_v )^{s}f}.
$$
Straightforward calculation by Fourier transformation shows that the
solution is in Gevrey class when $0<s\le \frac 12$ and is
ultra-analytic if $\frac 12<s<1$ for initial data only in
$L^2(\mathbb{R}^3_x\times\mathbb{R}^3_v)$ if it admits a unique
solution (see \cite{mo-xu2} for a more general study). However,
for the Boltzmann equation, the gain of Gevrey regularity of
solution is a long lasting open problem which has only been proved
so far in the linear and spatially homogeneous setting, see \cite{MUXY-DCDS}.

The mathematical study on the inverse power law potentials can be
traced back to the work by Pao \cite{Pao} in the 1970s. And in the early
1980s, Arkeryd in \cite{Arkeryd} proved the existence of weak
solutions to the spatially homogeneous Boltzmann equation when
$0<s<\frac 12$, while  Ukai in \cite{ukai}  applied an abstract
Cauchy-Kovalevskaya theorem to obtain local solutions in the
functional space of functions, which are analytic in $x$ and Gevrey in $v$. However,
the smoothing effect of the collision operator was not studied at
that time. 

Since then, this problem  has attracted increasing
interests in the area of kinetic theory and a lot of progress has
been made on the existence and regularity theories. More precisely,
that the long-range interactions have smoothing effects on the
solutions to the Boltzmann equation was first proved by Desvillettes
for some simplified models, cf. \cite{D95,D97}. This is in contrast
with the hard sphere model and the potentials with Grad's angular
cutoff assumption. In fact, for the hard sphere model, the
cross-section has the form (in the $\sigma$ representation)
\begin{equation*}\label{hardball}
B(|v-v_*|, \cos\theta )=q_0 |v-v_*| ,
\end{equation*}
where $q_0$ is the surface area of a hard sphere.
For singular cross-sections,
Grad \cite{Grad}  introduced the idea to cut off the
singularity at $\theta=0$ so that
 $B(|v-v_*|, \cos\theta ) 
\in L^1( \SS^{2})$. This assumption has been widely accepted and is
now called Grad's angular cutoff assumption which influences a few
decades of mathematical studies on the Boltzmann equation. Under
this angular cutoff assumption, the solution has the same
regularity, at least in the Sobolev space, as the initial data. In
fact, it was shown in \cite{dly}, that the solution has the form
$$
f(t,x,v)=a(t,x,v) f(0, x-vt, v) +b(t,x,v),
$$
when the initial data $f(0,x,v)$ is in some weighted
$L^p_{x,v}$ space. Here, $a(t,x,v)$ and $b(t,x,v)$ are
in the Sobolev space $H^\delta_{t,x,v}$ for some $\delta>0$.
And the term $f(0,x-vt,v)$ just represents the free transport
so that it is clear that $f(t,x,v)$ and $f(0,x,v)$ have the same
regularity.

One of the main features of the Boltzmann equation is the celebrated
Boltzmann's H-theorem saying that the H-functional
$$%
H(t)=\int_{\R^3\times\R^3}f\log fdxdv,
$$
satisfies
$$
\frac{dH(t)}{dt} +D(t)=0,
$$
where
$$
D(t)= -\int_{\R^3 \times \R^3} Q(f,f)\log f dxdv \ge 0,
$$
which is called the entropy dissipation rate. Notice that $D(t)$ is
non-negative and vanishes only when $f$ is a Maxwellian. The
non-negativity of  $D$ indicates that the Boltzmann equation is a
dissipative equation. This fact is a basic ingredient in the $L^1$
theory of the Boltzmann equation, see for example \cite{D-L}.

By using  the entropy dissipation rate $D$ and the ``$Q^+$ smoothing
property'', the formal smoothing estimate was derived by P.-L. Lions (see the
complete references in \cite{al-1})
$$
\|\sqrt{f}(\sqrt f * \la v \ra^{-m} )\|^2_{{B}^{\delta,2}_\infty}
\le C\|f\|^{1-\theta}_{L^1}(
\|f\|_{L^1} +D(f)^{\frac 12})^{\theta}, \hskip0.5cm  \delta=\frac{s}{1+s},
\, \theta= \frac{1}{1+s}\,,
$$
for any constant $m>3$. Notice that the above regularity estimate is
on $\sqrt{f}$ but not $f$ itself. Later, some almost optimal estimates
together with some extremely useful results, such as the cancellation lemma,
were obtained in the work by Alexandre-Desvillettes-Villani-Wennberg
\cite{al-1}. By using these analytic tools, the mathematical theory
 regarding the regularizing effect
for the spatially homogeneous problems may now be considered as quite
satisfactory, see \cite{al-2-1,al-2-2,desv-F,Des-Fu-Ter,desv-wen1,HMUY,MUXY-DCDS,villani},
and the references therein.

However, for the spatially inhomogeneous equations, there are much
less results. The main difficulty comes from the coupling of the
transport operator with the collision operator, and the commutators
of the differential (pseudo-differential) operators with the
collision operator. There are two progresses which have been achieved so far. One is
about the local existence of solutions between two moving Maxwellians in
\cite{alex-two-maxwellian}, obtained by constructing upper and lower
solutions. The other one is about the global existence of renormalized solutions with defect measures constructed in \cite{al-3}, 
which becomes weak solutions if the defect measures vanish. Some results on similar but linear kinetic equations were also given
in \cite{amuxy-nonlinear-b} and \cite{bouchut}. In particular, a
generalized uncertainty principle \`a la Fefferman
\cite{Fefferman} (see also \cite{morimoto873, morimoto921,
morimoto-mori}) was introduced in \cite{amuxy-nonlinear-b} in order to prove
smoothing effects of the linearized and spatially inhomogeneous
Boltzmann equation  with non-cutoff cross sections, and get partial
smoothing effects for the nonlinear Boltzmann equation. In the following
analysis, this partial regularity result, together with its proof, will also be used.

This paper can be viewed as a continuation of our recent work
\cite{amuxy-nonlinear-b}. Under some mild regularity assumption on
the initial data, we will prove the existence of solutions and their $C^\infty$
regularity  with respect to all ( time, space and
velocity) variables. 

Even though it is still not known whether only
some natural bounds, such as total mass, energy and entropy on the
initial data,  can lead to the $C^\infty$ regularizing effect, as far as we know, 
the results shown in this paper are the first ones justifying the
$C^\infty$ regularizing effect for the nonlinear and spatially
inhomogeneous Boltzmann equation without Grad's angular cutoff
assumption.

In order to state our theorems, let us first introduce the notations
and assumptions used in this paper.

The non-negative cross-section $B(z, \sigma)$
(for a monatomic gas, which is the case considered herein) depends only on $|z|$ and the scalar product
$<\frac{z}{|z|}, \sigma>$. In most cases, the
 collision kernel  cannot be expressed explicitly, but
to capture the essential properties, it can be assumed to have the form
\[
B(|v-v_*|, \cos \theta)=\Phi (|v-v_*|) b(\cos \theta),\,\,\,\,\,
\cos \theta=\big<\frac{v-v_*}{|v-v_*|} \, ,\,\sigma\big>,\,\,\,
0\leq\theta\leq\frac{\pi}{2}.
\]

Furthermore, to keep the presentation as simple as possible, and in
particular to avoid the difficulty coming from the
vanishing of the cross-section at zero relative velocity,
 we assume that the kinetic factor  $\Phi$ in the cross-section is
modified as
\begin{equation}\label{E-soft-p}
\Phi (|v-v_*|)=\left(1+|v-v_*|^2\right)^{\frac \gamma 2}, \qquad \gamma\in\RR.
\end{equation}
This point can certainly be removed using our whole calculus, at the expense of technical and more complicated details.

Moreover, the angular factor is assumed to have
the following singular behavior
\begin{equation}\label{1.2}
\sin \theta \,\,b(\cos \theta)\,\, \approx \,\, K
\theta^{-1-2s},\,\,\,\mbox{when}\,\,\, \theta \rightarrow 0+ ,
\end{equation}
where  $0<s<1$ and $K$ is a positive constant. In fact $\gamma=0$
corresponds to the Maxwellian molecule,  $\gamma<0$ corresponds to
the modified soft potential, and $\gamma>0$ corresponds to the
modified hard potential.  The singularity will be called mild
for $ 0< s <{\frac 1 2}$ and strong for $\frac 12\le s<1$.
The case $s=\frac 12$ is critical in the sense that different computations
are required in many parts of our proofs for mild and strong singularities, as will be seen below. This is similar to the known fractional Laplacian studies.

It is now well known from the work \cite{al-1} that the singular behavior of the collision kernel
(\ref{1.2}) implies a sub-elliptic estimate in the
 velocity variable $v$.
In the following analysis,  we shall need a slightly
precised weighted sub-elliptic estimate 
in the velocity variable.
We shall show that for
$\gamma\in\RR$ and $0<s<1$, if $f\geq 0, \not \equiv 0\,, \,f\in L^1_2\bigcap L\log
L (\RR^3_v)$, there exists a constant $C>0$ such that for any
function $g\in H^{1}(\RR^3_v)$ we have
\begin{equation}\label{E-sub-estimate}
C^{-1}\|\Lambda^{s}_v W_{\gamma/2} g\|^2_{L^2(\RR^3_v)} \leq (-Q(f,
g),\, g)_{L^2(\RR^3_v)}+C\|f\|_{L^1_{\tilde{\gamma}}(\RR^3)}
\|g\|^2_{L^2_{\gamma^+/2}(\RR^3_v)},
\end{equation}
where $\tilde{\gamma}=\max(\gamma^+,{2-\gamma^+}), \gamma^+=\max
\{\gamma,\, 0\}$. Here $ W_l=W_l(v)=(1+|v|^2)^{l/2}= \left\langle
v\right\rangle^l, l\in\RR$, is the weight function in the variable
$v\in\RR^3$.

Similar sub-elliptic estimates, first proved in \cite{al-1} and then developed 
in many other works such as \cite{mouhot} in a linearized context, have been used crucially
at least for the following two aspects:

i) the proof of the regularizing effect on the
solutions to the spatially homogeneous Boltzmann equations, see
\cite{al-2-1,al-2-2, desv-wen1,HMUY, MUXY-DCDS};

ii) the proof of existence of solutions to the nonlinear and
spatially inhomogeneous Boltzmann equation
\cite{alex-two-maxwellian,al-3, villani2}.

In this paper, we will apply this tool in order to study the complete smoothing effect
for the
spatially inhomogeneous and nonlinear Boltzmann equation.

It is now well understood, see \cite{av-2} and references therein, that Landau equation
corresponds to the grazing limit of Boltzmann equation. However,
while Landau operator involves usual partial differential operators,
it should be kept in mind that fractional differential operators
appear in the Boltzmann case, see \cite{alex,alex-2}. Therefore,  the
analysis on the Boltzmann equation appears much more involved because it
requires the unavoidable use of
  Harmonic Analysis. In particular,
 we shall use a generalized uncertainty principle which
was introduced in \cite{amuxy-nonlinear-b}, and the estimation of
commutators used in the work \cite{mo-xu1} for the study of
hypo-elliptic properties.

Throughout this  paper, we shall use the following standard
weighted (with respect to the velocity variable $v\in\RR^3\,$)
Sobolev  spaces. For $m,\, l \in\RR$, set
$\RR^7=\RR_t\times\RR^3_x\times\RR^3_v$ and
$$
H^m_l  (\RR^7)=\Big\{f\in \cS^{\,'}(\RR^7);\,\, W_l(v)\, f\in
H^m(\RR^7)\,\Big\}\, ,
$$
which is a Hilbert space. Here $H^m$ is the usual Sobolev space. We
shall also use the functional spaces $H^k_l(\RR^6_{x,\, v})$ and
$H^k_l(\RR^3_{v})$, specifying the variables,  the
weight being always taken with respect to $v\in\RR^3$.

Since the regularity property to be proved here is local in space
and time, for convenience, we define the following
 local version of weighted Sobolev space. For
$-\infty\leq T_1<T_2\leq+\infty$, and any given open domain
$\Omega\subset\RR^3_x$, define
\begin{eqnarray*}
&&{\cH}^m_l  (]T_1, \, T_2[\times\Omega\times\RR^3_v)=
\Big\{f\in \cD^{\, '}(]T_1, \, T_2[\times\Omega\times\RR^3_v);\,\, \\
&&\hskip 3cm \varphi(t)\psi(x) f\in
H^m_l(\RR^7)\,,\,\,\forall\,\varphi\in C^\infty_0(]T_1, \,
T_2[),\,\,\psi\in C^\infty_0(\Omega)\, \Big\}.
\end{eqnarray*}

Our first main result is about the smoothing effect on the
solution and can be stated as follows

\begin{theo}\label{theo1}
{\bf (Regularizing effect on solutions)}

Assume that $0<s<1,\,\, \gamma\in\RR$, $-\infty\leq
T_1<T_2\leq+\infty$ and let $\Omega\subset\RR^3_x$ be an open domain.
 Let $f$ be a non-negative
 function belonging to $ {\mathcal
H}^5_l  (]T_1, \, T_2[\times\Omega\times{\mathbb R}^3_v)$ for all $l
\in{\mathbb N}$ and solving the Boltzmann equation (\ref{1.1}) in
the domain $]T_1, \, T_2[\times\Omega\times{\mathbb R}^3_v$ in the
classical sense. Furthermore, if $f$ satisfies  the non-vacuum
condition
\begin{equation}\label{1.4}
\|f(t, x, \cdot)\|_{L^1({\mathbb R}^3_v)}>0,
\end{equation}
for all $(t, x)\in ]T_1, \, T_2[\times\Omega$, then we have
$$
f\in \cH^{+\infty}_l(]T_1, \, T_2[\times\Omega\times\RR^3_v),
$$
for any $l\in\NN$, and hence
$$
f\in C^{\infty}(]T_1, \, T_2[\times\Omega;\, {\mathcal S}({\mathbb
R}^3_v)).
$$
\end{theo}

With this result in mind, a natural question is whether the Boltzmann
equation has solutions satisfying the assumptions imposed above. Let us recall that solutions constructed in
\cite{alex-two-maxwellian,al-3} do not work for
our purpose because of the lack of the weighted regularity
${\mathcal H}^5_l$. For Gevrey class solutions from \cite{ukai}, there is of course nothing to prove.

Thus, our second  main result is about the local existence and
uniqueness of solution for the Cauchy problem of the non-cutoff
Boltzmann equation. 

We consider the solution in the functional space
with Maxwellian type exponential decay in the velocity variable.
More precisely, for $m\in\RR$,  set
$$
\cE^m_0(\RR^6)=\Big\{g\in\cD'(\RR^6_{x, v});\, \exists \,\rho_0>0\,
\, s. t. \,\, e^{\rho_0 <v>^2} g\in H^m(\RR^6_{x, v}) \Big\},
$$
and for $T>0$
\begin{eqnarray*}
{\mathcal E}^m([0,T]\times{\mathbb R}^6_{x, v})&=&\Big\{f\in
C^0([0,T];{\mathcal D}'({\mathbb R}^6_{x, v}));\, \exists \,\rho>0
\\
&&\hskip 0.5cm s. t. \,\, e^{\rho \langle v \rangle^2} f\in C^0([0,
T];\,\, H^m({\mathbb R}^6_{x, v})) \Big\}.
\end{eqnarray*}

\begin{theo}\label{E-theo1}
Assume that $0<s<1/2$ and $\gamma+2s<1$. Let $f_0\geq 0$ and $ f_0 \in
\cE^{k_0}_0(\RR^6)$ for some  $4 \leq k_0 \in \NN$. Then, there exists
$T_* >0$ such that the  Cauchy
problem
\begin{equation}\label{1.1b}
\left\{\begin{array}{l}f_t+v\cdot\nabla_x f=Q(f, f),\\
f|_{t=0}=f_0,\end{array}\right.
\end{equation}
admits a non-negative and unique solution in the functional space
${\mathcal E}^{k_0}([0,T_*]\times{\mathbb R}^6)$.

Furthermore, if we assume that the initial data $ f_0$ is in ${\mathcal
E}^{5}_0({\mathbb R}^6)$ and does not vanish on a compact set
$K\subset {\mathbb R}^3_x$, that is,
$$
\|f_0(x,\,\cdot) \|_{L^1({\mathbb R}^3_v)}>0,\,\,\,\,\,\forall\,\, x\in K ,
$$
then we have the regularizing effect on the above solution,
that is,  there exist
$0<\tilde{T}_0\leq T_*$ and a neighborhood
$V_{0}$ of $K$ in ${\mathbb R}^3_x$ such that
$$
f \in C^\infty(]0, \tilde{T}_0[\times V_{0} ;\, {\mathcal
S}({\mathbb R}^3_{v})).
$$

Moreover, if $\gamma\leq 0$, the non-negative solution of the Cauchy
problem (\ref{1.1b}) is unique in the functional space $ C^0 ([0,
T_*];\,\, H^{m}_p({\mathbb R}^6))$ for $m>3/2+2s$, $p>3/2+4s$.
\end{theo}

\begin{rema}\label{rema1.1000}
For the inverse power law potential
 $r^{-(p-1)}$,
  the condition $0< s<1/2,\,
\gamma+2s<1$ corresponds to
$3< p<\infty$ which includes both soft and hard potentials.
\end{rema}

At the moment, it is not clear whether  we can relax the regularity
assumption initially made on the solutions. Note that for example,
the condition that $f\in L^1\cap L^\infty (\RR^7 )$ is enough to
give a meaningful sense to a weak formulation for the spatially
inhomogeneous Boltzmann equation. However, the analysis used here
can not be applied to this case, and so further study is needed. On
the other hand, the above two theorems give an answer to a long
lasting
 conjecture on the regularizing effect of the non-cutoff
cross-sections for the spatially inhomogeneous Boltzmann equation.

Finally in the introduction, let us review some related works on the
regularizing effect and  the existence of solutions for the Landau
equation. The  regularizing effect from the Landau collision
operator has been rather well studied. See
\cite{desv-vil,chen-li-xu1,av-2} for the spatially  homogeneous
case. For the spatially inhomogeneous problem, a regularizing result
was obtained in \cite{chen-desv-he}, where the $H^8$ regularity is
assumed on the solutions to start with. And similar result was also
recently proved for the Vlasov-Maxwell-Landau and  the
Vlasov-Poisson-Landau systems, cf.  \cite{chen} and the references
therein. As for the existence of solutions, see \cite{desv-vil}
where unique weak solutions  for spatially homogeneous case have
been constructed with rather general initial data, and see
\cite{guo} where the classical solutions for the spatially
inhomogeneous case have been constructed in a periodic box with
small initial data.

The rest of the paper will be organized as follows. First of all, in
the next section, we will use the pseudo-differential calculus to
study the upper bounds on the collision operator. We shall give a
precise coercivity estimate linked to the singularity in the
cross-section, and estimate the commutators between some
pseudo-differential operators and the nonlinear collision operators.
In Section 3, the regularizing effect will be proved under the
initial regularity assumption on the solution. The strategy of the
proof is as follows. We first choose  some suitable mollifiers such
that the mollified solutions can work as test functions for the
weak formulation of the problem.  We then establish a small gain of
the regularity in the velocity variable, by using the coercivity
estimate coming from the singularity of the cross section. On
account of the generalized uncertainty principle, a small gain of
the regularity in the space and time variables can be derived. The
$H^{+\infty}$ regularity will follow from an induction argument.
Finally, in Section 4,  local solutions to the non-cutoff Boltzmann
equation which meet the initialization condition of Theorem
\ref{theo1} are constructed, using a family of cutoff Boltzmann
equations with time local uniform bounds independent of cutoff
parameter in some weighted Sobolev space.  In particular, the uniform bounds are
established with the help of time dependent Maxwellian type weight
functions which were introduced in \cite{ukai,ukai-2}. The convergence of
the approximate solutions follows from compactness argument
, while the uniqueness of the solutions
can also be proved by using our sharp upper bounds on the
collision operator.

\renewcommand{\theequation}{\thesubsection.\arabic{equation}}
\vskip0.5cm
\section{Pseudo-differential calculus}\label{section2}
\smallskip

Under
 the non-cutoff cross section assumption, the Boltzmann collision operator is a (nonlinear) singular integral operator with respect to $v\in\RR^3_v$.
In the linearized case, it behaves like a pseudo-differential operator.

In this section, we study the pseudo-differential calculus on the
Boltzmann operator. It is one of the key 
analytic tools for proving the 
regularizing effect of the non-cutoff Boltzmann equation. 
Even though the regularity proved in this paper is local in space and time
variables, note that the
 collision operator is non-local in the 
space of $v$ variable. Moreover, since the kinetic 
factor in the cross-section is of the form  $\left\langle
v\right\rangle^\gamma$ which might be unbounded,
 we  need to consider the multiplication by the weight
function $W_l(v)$ of the pseudo-differential operators.
Hence, they are not the standard
pseudo-differential operators of order $0$ on the usual Sobolev
space. In other words, we shall consider
pseudo-differential operators with unbounded coefficients on the
weighted Sobolev space $H^m_l(\RR^3_v)$. The variables $(t, x)$ are
considered as parameters for the collision operators in this section.

\vskip0.5cm
\subsection{Upper bound estimates }\label{s2}
\setcounter{equation}{0}
\smallskip

We shall need some functional estimates on the Boltzmann collision
operator in the existence and regularization proofs below. The first one is 
 about the boundedness of the collision operator in
some weighted Sobolev spaces, see also \cite{alex,al-2,HMUY} .

\begin{theo}\label{theo2.1}
Let $0<s<1$ and $\gamma \in\R$. Then for any $m,
\,\alpha \in \RR$, there exists $C>0$ such that
\begin{equation}\label{2.2+001}
\|Q(f,\, g)\|_{H^{m}_{\alpha}(\RR^3_v)}\leq C\|f\|_{L^{1}_{{\alpha^+
+}( \gamma+2s)^+}(\RR^3_v)} \| g\|_{H^{m+2s}_{({\alpha +}
\gamma+2s)^+}(\RR^3_v)}\,
\end{equation}
for all $f\in L^{1}_{{\alpha^+ +}( \gamma+2s)^+}(\RR^3_v)$ and $
g\in H^{m+2s}_{({\alpha +} \gamma+2s)^+}(\RR^3_v)$ .
\end{theo}

\begin{rema}\label{rema2.2}.\\
{\rm (1) } The collision operator $Q(f,\, g)$ behaves differently with respect to $f$ and $g$:
 (\ref{2.2+001}) shows that, in some
sense, it is linear with respect to the second factor in the
 velocity variable $v$ because the action of differentiation of $Q(f,\,
g)$ with respect to $v$   goes only on 
$g$ when considered in the
Sobloev space. This is clear for the Landau
operator which is the grazing limit of the Boltzmann operator.
\\
{\rm (2) } The estimate (\ref{2.2+001}) is
in some sense optimal with respect to the
order of differentiation (exact order of $2s$) and also with respect to the
order of the
weight in $v$ coming from the cross-section. In \cite{HMUY}, 
the cases of both the modified hard potential
and Maxwellian molecule type cross-sections
 corresponding to $0\leq \gamma<1$ are
discussed. Let us also mention that a similar estimate 
was given in \cite{amuxy-nonlinear}, but it is not optimal in terms
of weight and differentiation.
However, its proof is more straightforward as it only uses
the  Fourier transformation of collision
operator (Bobylev's type formula \cite{bobylev} and see also the
Appendix of \cite{al-1}). For our purpose, the full precise
estimate  (\ref{2.2+001}) will be needed.
\end{rema}

\smallbreak \noindent {\bf Proof of Theorem \ref{theo2.1} : }

Firstly,  we consider the case when $\alpha =0$. To prove (\ref{2.2+001}) in
this case, it suffices to show that for any $m \in \R$
\begin{equation}\label{final2}
\left|\Big(Q(f,\,g),\,\,h\Big)_{L^2(\RR^3_v)}\right| \leq
C||f||_{L^1_{(\gamma+2s)^+}(\RR^3_v)}||g||_{H^{m+2s}_{(\gamma+2s)^+}
(\RR^3_v)}||h||_{H^{-m}(\RR^3_v)}.
\end{equation}
The proof needs some harmonic analysis tools based on the dyadic
decomposition. It is similar to the proof in \cite{HMUY}, where the
hard potential case $\gamma \geq 0$ was studied. Interested readers
may refer to the papers \cite{alex,al-1,HMUY} for more
details, though we will keep the paper self-contained.

Recall that
\begin{align*}
\Big(Q(f,\,g),\,\,h\Big)_{L^2(\RR^3_v)}& =
\displaystyle\int_{\R^{6}}\displaystyle\int_{\SS^{2}}b(\cos
\theta)f(v_*)\Phi(|v-v_*|)g(v)\{h(v^\prime)-h(v)\}d\sigma d v_* d v,
\end{align*}
where $\Phi(|v-v_*|)=\Phi(|v'-v_*'|)=\langle
v'-v_*'\rangle^{\gamma}$. Set
\[
F(v,v_*)=\Phi(|v-v_*|)g(v),
\]
and write
\begin{equation}\label{radja-1}
 \begin{split}
\Big(Q(f,\,g),\,\,h\Big)_{L^2(\RR^3_v)}&=\displaystyle\int_{\R^{6}}
\displaystyle\int_{\SS^{2}}b(\cos
\theta)f(v_*)F(v,v_*)\{h(v^\prime)-h(v)\}d\sigma d v_* d v\\
&=
 \int_{\R^3}f(v_*)(U_1-U_2)dv_*.
\end{split}
\end{equation}
Then we have (formally) by inverse Fourier formula,
\begin{equation*} 
U_1\equiv \int_{\R^{3}}\displaystyle\int_{\SS^{2}}b(\cos
\theta)F(v,v_*)h(v^\prime)d\sigma dv=\int_{\R^{3}}\int_{\R^{3}} H(\xi,\eta, v_*)
\hat{F}(\xi, v_*) \overline{\hat{h}(\eta)}d\xi d\eta,
\end{equation*}
where (also formally)
\begin{align*}
H(\xi,\eta, v_*)&=\int_{\R^{3}}\int_{\SS^{2}}b(k\cdot \sigma)e^{iv\cdot
\xi-iv'\cdot\eta } d\sigma dv
\\
&=\int_{\R^{3}} e^{iv\cdot
\xi-i\frac{v+v_*}{2}\cdot\eta}\Big[\int_{\SS^{2}}b(k\cdot \sigma)
e^{-i\frac{|v-v*|}{2} \sigma\cdot\eta } d\sigma \Big] dv
\\
&=\int_{\R^{3}} e^{iv\cdot
\xi-i\frac{v+v_*}{2}\cdot\eta}\Big[\int_{\SS^{2}}b(\tilde{\eta}\cdot
\sigma) e^{-i\frac{|v-v*|}{2}|\eta| \sigma\cdot k } d\sigma \Big]
dv\, , \enskip\enskip( \tilde{\eta} = \eta/|\eta| )
\\
&=\int_{\R^{3}} e^{iv\cdot
\xi-i\frac{v+v_*}{2}\cdot\eta}\Big[\int_{\SS^{2}}b(\tilde{\eta}\cdot
\sigma) e^{-i|\eta|\frac{v-v*}{2} \cdot\sigma } d\sigma \Big] dv
\\
 &=\int_{\SS^{2}}b(\tilde{\eta}\cdot \sigma)e^{-iv_*\cdot\eta^-}\Big[
 \int_{\R^{3}} e^{iv\cdot (\xi-\eta^+)}dv\Big]
 d\sigma
\\
 &=\int_{\SS^{2}}b(\tilde{\eta}\cdot \sigma)e^{-iv_*\cdot\eta^-}d\sigma
 \ \delta(\xi-\eta^+),
\end{align*}
with
$$
\eta^-=\frac 12(\eta-|\eta|\sigma), \quad \eta^+=\frac
12(\eta+|\eta|\sigma),
$$
so that
\begin{equation*}
 U_1=\int_{\R^{3}} \Big[\int_{\SS^{2}}b(\tilde{\eta}\cdot \sigma)
e^{-iv_*\cdot\eta^-}d\sigma\Big]
 \hat{F}(\eta^+, v_*)\overline{\hat{h}(\eta)} d\eta.
 \end{equation*}
 On the other hand,
\begin{align*}
 U_2\equiv &\int_{\R^{3}}\displaystyle\int_{\SS^{2}}b(\cos
\theta)F(v,v_*)h(v)d\sigma dv
\\
=& \Big[ \int_{\SS^{2}}b(\cos \theta)d\sigma \Big] \int_{\R^{3}}
\hat{F}(\eta, v_*) \overline{\hat{h}(\eta)} d\eta \qquad\quad
\\
=& \int_{\R^{3}} \Big[\int_{\SS^{2}}b(\tilde{\eta}\cdot\sigma)
d\sigma\Big]\hat{F}(\eta, v_*) \overline{\hat{h}(\eta)} d\eta,
\end{align*}
because (formally) we have
\begin{equation*}
 \int_{\SS^{2}}b(\cos \theta)d\sigma
=\int_{\SS^{2}}b(\tilde{\eta}\cdot\sigma)d\sigma={\rm const.}
\end{equation*}
Therefore, we have obtained the following generalized Bobylev formula
\begin{equation}\label{Qfourier}
 \begin{split}
&\Big(Q(f,\,g),\,\,h\Big)_{L^2(\RR^3_v)} \\
=&\int_{\R^3}
f(v_*)\Big[\int_{\R^3} \int_{\SS^{2}}b(\tilde{\eta}\cdot\sigma)
\Big\{e^{-iv_*\cdot\eta^-} \hat{F}(\eta^+, v_*)- \hat{F}(\eta,
v_*)\Big\} \overline{\hat{h}(\eta)} d\eta d\sigma \Big]dv_*
\\
 =&\int_{\R^3} f(v_*)\Big[\int_{\R^3}
\int_{\SS^{2}}b(\tilde{\eta}\cdot\sigma)
\Big\{e^{iv_*\cdot\eta^+} \hat{F}(\eta^+, v_*)-
e^{iv_*\cdot\eta} \hat{F}(\eta, v_*)\Big\}\\
 &\hskip6cm\times \overline{   e^{iv_*\cdot\eta}
\hat{h}(\eta)} d\eta d\sigma \Big]dv_* .
 \end{split}
\end{equation}
Notice that  the above derivation is only formal for
non-cutoff cross-section because we can not split the gain and
loss term in this case. However, the derivation
can be easily justified as a limit process of cutoff cross-sections
when combining the gain term and loss term together.

We now introduce a dyadic decomposition in $\R^3_v$ as follows:
\begin{equation*}\sum_{k=0}^{\infty} \phi_k(v) =1\,,
\enskip \enskip  \phi_k(v) = \phi(2^{-k} v) \enskip
\mbox{for} \enskip k \geq 1 \enskip \mbox{with} \enskip 0 \leq
\phi_0, \phi \in C_0^{\infty}(\R^3),
\end{equation*}
and
\begin{equation*}
\text{supp}~\phi_0 \subset \{|v|<2\}, \qquad\quad \text{supp}~\phi
\subset \{ 1 <|v|< 3\}.
\end{equation*}
Take also $\tilde \phi_0$ and $\tilde \phi \in C_0^\infty$ such that
\begin{eqnarray*}
&&\tilde \phi_0 = 1 \ \ \text{on} \ \  \{|v|\leq 2\}, \qquad \qquad
\quad \text{supp}~\tilde{\phi}_0 \subset \{|v|<3\},
\\
&&\tilde \phi = 1 \ \ \text{on} \ \  \{  1/2 \leq |v|\leq 3\},
\qquad \quad  \text{supp}~\tilde{\phi} \subset \{1/3 <|v|<4\}.
\end{eqnarray*}
Furthermore, we assume that all these functions are radial. From  $|v'-v_*| \leq |v-v_*| \leq \sqrt{2}
|v'-v_*|$, it follows that
\begin{equation*}
\tilde \phi_k(v'-v_*) \phi_k(v-v_*) = \phi_k(v-v_*)=\tilde
\phi_k(v-v_*) \phi_k(v-v_*),\enskip k \geq 0,
\end{equation*}
and thus we get
\begin{align*}
\Big(Q(f,\,g),\,\,h\Big)_{L^2(\RR^3_v)}& =\sum_{k=0}^\infty
\displaystyle\int_{\R^{6}}\displaystyle\int_{\SS^{2}}b(\cos
\theta)f(v_*)F_k(v,v_*)\{h_k(v^\prime,v_*)-h_k(v,v_*)\}d\sigma d v_*
d v,
\end{align*}
where
\begin{equation}\label{h-k}
F_k(v,v_*)=\phi_k(v-v_*) \Phi(|v-v_*|)g(v), \qquad \qquad
h_k(v,v_*)=\tilde \phi_k(v-v_*)h(v).
\end{equation}
Similarly to \eqref{Qfourier},  we also obtain
\begin{align*}
\Big(Q(f,\,g),\,\,h\Big)_{L^2(\RR^3_v)} =& \sum_{k=0}^\infty
\int_{\R^3} f(v_*)\Big[\int_{\R^3}
\int_{\SS^{2}}b(\tilde{\eta}\cdot\sigma) \Big\{e^{iv_*\cdot\eta^+}
\hat{F_k}(\eta^+, v_*)- e^{iv_*\cdot\eta}
\hat{F_k}(\eta, v_*)\Big\}\\
&\hskip3cm \times \overline{   e^{iv_*\cdot\eta} \hat{h_k}(\eta,v_*)}
d\eta d\sigma \Big]dv_* \\
=& \int_{\R^3} f(v_*)   \sum_{k=0}^\infty   K^k(v_*) dv_* .
\end{align*}
In the following, we will estimate $\sum_{k=0}^\infty   |K^k(v_*)|$, regarding $v_*$
as a parameter.

By setting
\begin{equation*}
\Omega_k = 
\left \{ \sigma \in {\mathbb S}^{2}\, ; \,
\tilde \eta \cdot \sigma \geq 1 - 2^{1-2k } \langle \eta
\rangle^{-2}\right \}
,
\end{equation*}
and
\begin{equation*}
\widetilde{\hat{F_k}}(\eta, v_*) =e^{iv_*\cdot\eta} \hat{F_k}(\eta,
v_*), \enskip \widetilde{\hat{h_k}}(\eta,v_*) =e^{iv_*\cdot\eta}
\hat{h_k}(\eta,v_*),
\end{equation*}
we split $K^k(v_*)$ into
\begin{align*}
K^k(v_*) =& \int_{\R^3} \int_{\SS^{2}\cap
\Omega_k}b(\tilde{\eta}\cdot\sigma)
\Big\{\widetilde{\hat{F_k}}(\eta^+, v_*)-
\widetilde{\hat{F_k}}(\eta, v_*)\Big\}
\overline{\widetilde{\hat{h_k}}(\eta,v_*)} d\eta d\sigma\\
+& \int_{\R^3} \int_{\SS^{2}\cap
\Omega_k^c}b(\tilde{\eta}\cdot\sigma)
\Big\{\widetilde{\hat{F_k}}(\eta^+, v_*)-
\widetilde{\hat{F_k}}(\eta, v_*)\Big\}
\overline{\widetilde{\hat{h_k}}(\eta,v_*)} d\eta d\sigma   \\
=& K^k_1(v_*) + K^k_2(v_*).
\end{align*}
Note that
\begin{align}\label{small-sing-1}
\int_{\SS^{2}\cap
\Omega_k}\theta^2\,\, b(\cos \theta) d\sigma &= 2\pi\int_{\{\theta
\in [0,\pi/2]; \,\sin(\theta/2) \leq 2^{-k} \langle
\eta \rangle^{-1}\}} \sin \theta\,\, b(\cos \theta) \theta^2 d\theta\\
&\notag \leq C \la
\eta \ra^{2s-2} 2^{k(2s-2)},\,\,\, \enskip \mbox{if $0<s<1$,}
\end{align}
\begin{align}\label{non-sing}
\int_{\SS^{2}\cap \Omega_k^c}\, b(\cos \theta) d\sigma 
&= 2\pi\int_{\{\theta
\in [0,\pi/2]; \,\sin(\theta/2) \geq 2^{-k} \langle
\eta \rangle^{-1}\}} \sin \theta\,\, b(\cos \theta)  d\theta
\\
&\leq  \notag C \la
\eta \ra^{2s} 2^{2ks},\,\,\, \enskip \mbox{for any $s>0$}.
\end{align}

It follows  from \eqref{non-sing} that 
\begin{align}\label{K2k}
|K^k_2(v_*)|\leq&\int_{\R^3} \int_{\SS^{2}\cap \Omega_k^c}b(\tilde{\eta}\cdot\sigma)
\Big|\widetilde{\hat{F_k}}(\eta^+, v_*)-\widetilde{\hat{F_k}}(\eta,
v_*)\Big|
\Big|\widetilde{\hat{h_k}}(\eta,v_*)\Big| d\eta d\sigma\\ \notag
\leq&\left( \int_{\R^3} \int_{\SS^{2}\cap
\Omega_k^c}b(\tilde{\eta}\cdot\sigma)\la \eta \ra^{2m+2s}
\Big(\Big|\widetilde{\hat{F_k}}(\eta^+, v_*)\Big|^2
+\Big|\widetilde{\hat{F_k}}(\eta, v_*)\Big|^2\Big)d\eta
d\sigma\right)^{1/2}\\ \notag
&\times\left( \int_{\R^3} \int_{\SS^{2}\cap
\Omega_k^c}b(\tilde{\eta}\cdot\sigma)\la \eta \ra^{-2m-2s}
\Big|\widetilde{\hat{h_k}}(\eta,v_*)\Big|^2 d\eta
d\sigma\right)^{1/2}\\ \notag
\leq& C 2^{2ks}\left\| \la D_v \ra^{m+2s}F_k(v,v_*)\}\right\|_{L^2}
||\la D_v \ra^{-m} h_k(v,v_*)||_{L^{2}}.
\end{align}
Here, we have used  the  change of variables $\eta\, \rightarrow\,\eta^+$, which is regular because the Jacobian can be computed, 
with $\tilde{\eta}=\eta/|\eta|$, as
\[
\Big|\frac{\partial (\eta^+)}{\partial (\eta)}\Big|
=\Big|\frac 12 I+\frac 12 \sigma\otimes \tilde{\eta}\Big|
=\frac 18(1+ \sigma\cdot\tilde{\eta})=\frac 14 \cos^2\frac{\theta}{2}.
\]
It should be noted that  after this change of variable,
$\theta$ plays no longer the role of the polar
angle because the ``pole" $\tilde{\eta}$ now moves with $\sigma$ and hence the measure
$d\sigma$ is no longer given by $\sin\theta d\theta d\phi$.
However, the situation is rather good because  if we   
take $\tilde{\eta}^+=\eta^+/|\eta^+|$
as a new pole which is independent of $\sigma$,  then the  new polar angle
$\psi$ defined by $\cos\psi=\tilde{\eta}^+\cdot\sigma$  satisfies
\[
\psi=\frac{\theta}{2}, \qquad d\sigma=\sin\psi d\psi d\phi, \quad \psi\in[0,\frac{\pi}{4}],
\]
and thus $\theta$ works almost as the polar angle. Therefore, since
$\la \eta \ra\leq 2\la \eta^+ \ra\leq 2\la \eta
\ra$ we have 
\[
 \int_{\R^3} \int_{\SS^{2}\cap
\Omega_k^c}b(\tilde{\eta}\cdot\sigma) \la\eta \ra^{2m+2s}
\Big|\widetilde{\hat{F_k}}(\eta^+, v_*)\Big|^2
d\eta
d\sigma
\le C
\int_{\RR^3_{\eta^+}}
D_0(\eta^+)
\Big|\widetilde{\hat{F_k}}(\eta^+, v_*)\Big|^2
d\eta^+
\]
with
\begin{align*}
D_0(\eta^+)& =
\int_{\SS^{2}\cap
\Omega_k^c}b(\tilde{\eta}\cdot\sigma) \la\,\eta(\eta^+, \sigma) \,\ra^{2m+2s}
d\sigma\\
&\leq C
\int_{\mathbb S^2 \cap \Omega^c_k} \la\,\eta(\eta^+, \sigma)\, \ra^{2m+2s}
\theta^{-2-2s}d\sigma
\\
&
\le C 
\la\eta^+ \ra^{2m+2s}\int_{2^{-k} \la \eta^+ \ra^{-1}}^{\pi/4}\psi^{-2-2s}\sin\psi d\psi \le 
2^{2ks} \la \eta^+ \ra^{2m+4s}\,,
\end{align*}
which implies \eqref{K2k}.
Notice that for $p =0, 1, 2$,
\begin{align}\label{bloc-estimate}\displaystyle
\Big|\frac{2^{k(2s-p)}|v-v_*|^p \phi_k(v-v_*) \Phi(|v-v_*|)}
{\la v_* \ra^{(\gamma+2s)^+}}\Big| &\leq\displaystyle C \frac{\la
v-v_* \ra^{\gamma+2s}}{\la v_* \ra^{(\gamma+2s)^+}}\phi_k
(v-v_*)\\
&\leq C\la v\ra^{(\gamma+2s)^+}\phi_k(v-v_*).\nonumber
\end{align}
Then, recalling \eqref{h-k} and using \eqref{bloc-estimate} with $p=0$ we have
\begin{align*}
|K_2^{k}(v_*)| \leq& C \la v_* \ra^{(\gamma+2s)^+} \left\|\frac{ \la
D_v \ra^{m+2s}}{\la v_* \ra^{(\gamma+2s)^+}}\{2^{2ks}
F_k(v,v_*)\}\right\|_{L^2}
||\la D_v \ra^{-m} h_k(v,v_*)||_{L^{2}}\\
\leq & C \la v_* \ra^{(\gamma+2s)^+} \Big(\|\tilde\phi_k(v-v_*) \la
D_v \ra^{m+2s}g\|_{L^2_{(\gamma+2s)^+}}^2+ 2^{-k} \|\la D_v
\ra^{m+2s}g\|_{L^2_{(\gamma+2s)^+}}^2
\Big)^{1/2}\\
 &\times\Big(\|\tilde \phi_k(v-v_*) \la D_v
\ra^{-m}h\|_{L^2}^2 +
2^{-k}\|\la D_v \ra^{-m}h\|_{L^2}^2\Big)^{1/2}\\
 := & C\, \Gamma_k(v_*),
\end{align*}
where $\Gamma_k(v_*)$ stands for the quantity  defined by this right hand
side up to a constant multiple.

On the other hand, in order to estimate $K^k_1(v_*)$, write
\begin{align*}
\Big\{\widetilde{\hat{F}_k}(\eta^+, v_*)-
\widetilde{\hat{F}_k}(\eta, v_*)\Big\}
\overline{\widetilde{\hat{h}_k}(\eta,v_*)} &=
\Big\{\widetilde{\hat{F}_k}(\eta^+, v_*)-
\widetilde{\hat{F}_k}(\eta, v_*)\Big\}
\Big\{ \overline{\widetilde{\hat{h}_k}(\eta,v_*)} -
\overline{\widetilde{\hat{h}_k}(\eta^+,v_*)} \Big\}\\
- & \eta^- \cdot \Big(\nabla \widetilde{\hat{F_k}}\Big)(\eta^+, v_*)
\overline{\widetilde{\hat{h_k}}(\eta^+,v_*)} \\
-\int_0^1 \Big\{\Big(\nabla \widetilde{\hat{F}_k}\Big)(\eta^+ + \tau
& (\eta - \eta^+), v_*) -   \Big(\nabla
\widetilde{\hat{F}_k}\Big)(\eta^+ , v_*) )d \tau\Big\} \cdot (\eta -
\eta^+) \overline{\widetilde{\hat{h}_k}(\eta^+,v_*)}.
\end{align*}
Correspondingly, we decompose $K^k_1(v_*)$ into
\begin{equation*}
K^k_1(v_*) = K^{k,1}_1(v_*) + K^{k,2}_1(v_*)+K^{k,3}_1(v_*).
\end{equation*}
For the variable transformation $\eta \longrightarrow  \eta^+ =
\frac{1}{2}(\eta + |\eta|\sigma)$, we denote its inverse
transformation $\eta^+ \longrightarrow   \eta$ by
$\psi_{\sigma}(\eta^+)$. Then
\begin{align*}
K^{k,2}_1(v_*) &=-\int_{\R^3} \int_{\SS^{2}\cap \Omega_k}b \Big(
\frac{\psi_{\sigma}(\eta^+)}{|\psi_{\sigma}(\eta^+)|}\cdot\sigma
\Big) \left|\frac{\partial (\psi_{\sigma}(\eta^+))}
{\partial(\eta^+)}\right|\\
&\hskip3cm \times \eta^-(\sigma)\cdot \Big(\nabla \widetilde{
\hat{F_k}}\Big)(\eta^+, v_*)
\overline{\widetilde{\hat{h_k}}(\eta^+,v_*)} d\eta^+ d\sigma  \\
& = 0 \,, \enskip \mbox{with} \enskip \eta^-(\sigma) =
\psi_{\sigma}(\eta^+)-\eta^+,
\end{align*}
because $ \sigma_1, \sigma_2 \in \SS^2 \cap \Omega_k$ are symmetric with
respect 
to each other in the sense that, cf  Figure \ref{fig2},
\begin{equation*}
\eta^{-}(\sigma_1) = \psi_{\sigma_1}(\eta^+) -\eta^+ =
-(\psi_{\sigma_2}(\eta^+) -\eta^+) = -\eta^{-}(\sigma_2).
\end{equation*}

Write $K^{k,1}_1(v_*)$ into
\begin{align*}
K^{k,1}_1(v_*) =& -\int_0^1\int_0^1\Big( \int_{\R^3}
\int_{\SS^{2}\cap \Omega_k}b(\tilde{\eta}\cdot\sigma)
\Big\{\Big(\nabla \widetilde{\hat{F_k}}\Big)(\eta^++
\tau(\eta-\eta^+),v_*) \cdot (\eta-\eta^+) \Big\}
\\
&\hskip3cm \times \Big\{ \overline{\Big(\nabla
\widetilde{\hat{h_k}}\Big)(\eta^++ s(\eta-\eta^+),v_*)\cdot (\eta-
\eta^+)\Big\}} d\eta d\sigma   \Big) d \tau ds  .
\end{align*}
Since $|\eta- \eta^+|^2=|\eta^-|^2=|\eta|^2\sin^2 (\theta/2)$ and the change of variable $\eta^++ \tau(\eta-\eta^+)\,
\rightarrow\,\eta$ is also regular (see Page 2044 of \cite{amuxy-nonlinear-b}),  \eqref{small-sing-1}
implies
\begin{align*}
|K^{k,1}_1(v_*)|\leq & C\int_0^1\int_0^1\Big( \int_{\R^3}
\int_{\SS^{2}\cap \Omega_k}\theta^2\,\,b(\tilde{\eta}\cdot\sigma)
\la\eta\ra^2 \Big|\Big(\nabla \widetilde{\hat{F_k}}\Big)(\eta^++
\tau(\eta-\eta^+),v_*)\Big|
\\
&\hskip3cm \times \Big|{\Big(\nabla
\widetilde{\hat{h_k}}\Big)(\eta^++ s(\eta-\eta^+),v_*)\Big|} d\eta
d\sigma   \Big) d \tau ds\\
\leq & C\int_0^1\int_0^1\Big( \int_{\R^3} \int_{\SS^{2}\cap
\Omega_k}\theta^2\,\,b(\tilde{\eta}\cdot\sigma) \la\eta\ra^{2+2s+2m}
\Big|\Big(\nabla \widetilde{\hat{F_k}}\Big)(\eta^++
\tau(\eta-\eta^+),v_*)\Big|d\eta d\sigma   \Big)^{1/2}
\\
&\times\Big( \int_{\R^3} \int_{\SS^{2}\cap
\Omega_k}\theta^2\,\,b(\tilde{\eta}\cdot\sigma)
\la\eta\ra^{2-2s-2m}\Big|{\Big(\nabla
\widetilde{\hat{h_k}}\Big)(\eta^++ s(\eta-\eta^+),v_*)\Big|} d\eta
d\sigma   \Big)^{1/2} d \tau ds \\
\leq& C2^{k(2s-2)}\|\la\eta\ra^{2s+m}\Big(\nabla
\widetilde{\hat{F_k}}\Big)\|_{L^2(\RR^3)}
\|\la\eta\ra^{-m}\Big(\nabla
\widetilde{\hat{h_k}}\Big)\|_{L^2(\RR^3)}  .
\end{align*}
Hence, we have obtained, by using \eqref{bloc-estimate} with $p=1$
\begin{align*}
|K_1^{k,1}(v_*)| \leq& C \la v_* \ra^{(\gamma+2s)^+}
\left\|\frac{\la D_v \ra^{m+2s}}{\la v_*
\ra^{(\gamma+2s)^+}}\{2^{k(2s-1)}(v-{v_*})
F_k(v,v_*)\}\right\|_{L^2}
\\
&\hskip3cm \times \|2^{-k}(v -{v_*})h_k(v,v_*)\|_{H^{-m}}\\
\leq &C \la v_* \ra^{(\gamma+2s)^+} \Big(\|\tilde\phi_k(v-v_*) \la
D_v \ra^{m+2s}g\|_{L^2_{(\gamma+2s)^+}}^2+ 2^{-k}\|\la
D_v\ra^{m+2s}g\|_{L^2_{(\gamma+2s)^+}}^2
\Big)^{1/2}\\
&\times\Big(\|\tilde \phi_k(v-v_*) \la D_v \ra^{-m}h\|_{L^2}^2 +
2^{-k}\|\la D_v\ra^{-m}h\|_{L^2}^2\Big)^{1/2},
\end{align*}
which has the same bound $\Gamma_k(v_*)$ as in the previous case, up
to a constant factor. Finally, we consider
\begin{align*}
K^{k,3}_1(v_*)=& -\int_0^1 \int_0^1\Big( \int_{\R^3}
\int_{\SS^{2}\cap \Omega_k}b(\tilde{\eta}\cdot\sigma)
\Big\{\Big(\nabla^2 \widetilde{\hat{F_k}}\Big)(\eta^++ \tau
s(\eta-\eta^+),v_*) \tau(\eta-\eta^+)^2 \Big\}
\\
&\hskip5cm \times \Big\{ \overline{
\widetilde{\hat{h_k}}(\eta^+,v_*)\Big\}} d\eta d\sigma  \Big)d \tau
ds.
\end{align*}
Then, by using \eqref{bloc-estimate} with $p=2$, we have 
\begin{align*}
|K_1^{k,3}(v_*)| \leq& C \la v_* \ra^{(\gamma+2s)^+}
\left\|\frac{\la D_v \ra^{m+2s}}{\la v_*
\ra^{(\gamma+2s)^+}}\{2^{k(2s-2)}(v-{v_*})^2
F_k(v,v_*)\}\right\|_{L^2}
\|h_k(v,v_*)\|_{H^{-m}}\\
\leq &C \Gamma_k(v_*).
\end{align*}
Therefore, it follows from Schwarz's inequality that
\begin{align*}
&\left|\Big(Q(f,\,g),\,\,h\Big)_{L^2(\RR^3_v)}\right| \leq
C ||f||_{L^1_{(\gamma+2s)^+}}\times \\
&\times \Big(\sum_{k=0}^\infty   \{ \|\tilde\phi_k(v-v_*) \la D_v
\ra^{m+2s}g\|_{L^2_{(\gamma+2s)^+}}^2+ 2^{-k} \|\la D_v
\ra^{m+2s}g\|_{L^2_{(\gamma+2s)^+}}^2   \}
 \Big)^{1/2}\\
&\times\Big(\sum_{k=0}^\infty \{ \|\tilde \phi_k(v-v_*) \la D_v
\ra^{-m}h\|_{L^2}^2 + 2^{-k}\|\la D_v \ra^{-m}h\|_{L^2}^2\}
\Big)^{1/2}\\
&\leq C||f||_{L^1_{(\gamma+2s)^+}}||g||_{H^{m+2s}_{(\gamma+2s)^+}
}||h||_{H^{-m}},
\end{align*}
which yields \eqref{final2}. Now the proof of Theorem \ref{theo2.1}
is complete for the case $\alpha =0$.

To prove (\ref{2.2+001}) for the case $\alpha \ne 0$, it suffices
to show that 
\begin{equation}\label{final2+general}
\left|\Big(Q(f,\,g),\,\,\langle v \rangle^\alpha h\Big)_{L^2(\RR^3_v)}\right| \leq
C||f||_{L^1_{\alpha^++ (\gamma+2s)^+}(\RR^3_v)}||g||_{H^{m+2s}_{(\alpha+ \gamma+2s)^+}
(\RR^3_v)}||h||_{H^{-m}(\RR^3_v)}.
\end{equation}
The argument is similar to the one for
$\alpha=0$, up to the estimation on  $h_k(v,v_*)$ in \eqref{h-k}
which
must be replaced by
$$\tilde \phi_k(v-v_*) \la v \ra^\alpha h(v) = \la v \ra^\alpha h_k(v,v_*). $$
We can write
\begin{align}
\la v \ra^\alpha h_k(v,v_*) =& (\la v_* \ra^\alpha 
+ 2^{k \alpha } )  \psi_k(v,v_*) h_k(v,v_*),& \mbox{if $\alpha >0$}, \label{alpha+}\\
\la v \ra^\alpha h_k(v,v_*) =&
\left(\frac{\la v_* \ra}{2^k}\right)^{\min\{(\gamma+2s)^+, -\alpha\}}\psi_k(v,v_*)
 h_k(v,v_*),&   \mbox{if $\alpha < 0$}. \label{alpha-}
\end{align}
with a suitable $\psi_k(v,v_*)$ belonging to $C_b^\infty (\RR_v^3)$, uniformly with respect to $k, v_*$.
For $p =0, 1, 2$, we have
\begin{align}\label{+bloc-estimate}\displaystyle
&\Big|\frac{2^{k(\alpha + 2s-p)}|v-v_*|^p \phi_k(v-v_*) \Phi(|v-v_*|)}
{\la v_* \ra^{(\alpha + \gamma+2s)^+}}\Big| \\
&\leq C \displaystyle \frac{\la
v-v_* \ra^{\alpha+ \gamma+2s}}{\la v_* \ra^{(\alpha+ \gamma+2s)^+}}\phi_k
(v-v_*)
\leq C \la v\ra^{(\alpha+ \gamma+2s)^+}\phi_k(v-v_*) \notag ,
\end{align}
which is similar to \eqref{bloc-estimate}.
We first consider the case  $\alpha >0$. 
It follows  from \eqref{non-sing} that 
\begin{align*}\label{K2k}
|K^k_2(v_*)|\leq& (\la v_* \ra^\alpha 
+ 2^{k \alpha } )\int_{\R^3} \int_{\SS^{2}\cap \Omega_k^c}b(\tilde{\eta}\cdot\sigma)
\Big|\widetilde{\hat{F_k}}(\eta^+, v_*)-\widetilde{\hat{F_k}}(\eta,
v_*)\Big|
\Big|\widetilde{\widehat{\psi_k h_k}}(\eta,v_*)\Big| d\eta d\sigma\\ \notag
\leq& (\la v_* \ra^\alpha 
+ 2^{k \alpha } )\left( \int_{\R^3} \int_{\SS^{2}\cap
\Omega_k^c}b(\tilde{\eta}\cdot\sigma)\la \eta \ra^{2m+2s}
\Big(\Big|\widetilde{\hat{F_k}}(\eta^+, v_*)\Big|^2
+\Big|\widetilde{\hat{F_k}}(\eta, v_*)\Big|^2\Big)d\eta
d\sigma\right)^{1/2}\\ \notag
&\times\left( \int_{\R^3} \int_{\SS^{2}\cap
\Omega_k^c}b(\tilde{\eta}\cdot\sigma)\la \eta \ra^{-2m-2s}
\Big|\widetilde{\widehat{\psi_k h_k}}(\eta,v_*)\Big|^2 d\eta
d\sigma\right)^{1/2}\\ \notag
\leq& C 2^{2ks}(\la v_* \ra^\alpha 
+ 2^{k \alpha } )\left\| \la D_v \ra^{m+2s}F_k(v,v_*)\}\right\|_{L^2}
||\la D_v \ra^{-m} h_k(v,v_*)||_{L^{2}}.
\end{align*}
Then, recalling \eqref{h-k}, and using \eqref{bloc-estimate} and  \eqref{+bloc-estimate} with $p=0$, we have
\begin{align*}
|K_2^{k}(v_*)| \leq& C \left \{\la v_* \ra^{\a+ (\gamma+2s)^+} \left\|\frac{ \la
D_v \ra^{m+2s}}{\la v_* \ra^{(\gamma+2s)^+}}\{2^{2ks}
F_k(v,v_*)\}\right\|_{L^2} \right. \\
& + \left .\la v_* \ra^{(\a+ \gamma+2s)^+} \left\|\frac{ \la
D_v \ra^{m+2s}}{\la v_* \ra^{(\a +\gamma+2s)^+}}\{2^{k(\a+2s) }
F_k(v,v_*)\}\right\|_{L^2}\right\}
||\la D_v \ra^{-m} h_k(v,v_*)||_{L^{2}}\\
\leq & C \la v_* \ra^{\a+ (\gamma+2s)^+} \Big(\|\tilde\phi_k(v-v_*) \la
D_v \ra^{m+2s}g\|_{L^2_{(\a+\gamma+2s)^+}}^2+ 2^{-k} \|\la D_v
\ra^{m+2s}g\|_{L^2_{(\a+\gamma+2s)^+}}^2
\Big)^{1/2}\\
 &\times\Big(\|\tilde \phi_k(v-v_*) \la D_v
\ra^{-m}h\|_{L^2}^2 +
2^{-k}\|\la D_v \ra^{-m}h\|_{L^2}^2\Big)^{1/2}\\
 := & C\, \Gamma_k^\a(v_*) ,
\end{align*}
where $\Gamma_k^\a(v_*)$ stands for the quantity  defined by this right hand
side up to a constant multiple. 

Performing the same computation as above for $K_2^{k}(v_*)$,
it follows from \eqref{small-sing-1} that
$$|K_1^{k,1}(v_*)| + |K_1^{k,3}(v_*)|\leq C \, \Gamma_k^\alpha (v_*) ,$$
so that  \eqref{final2+general} holds in this case.

The estimation on the
case $\alpha <0$ is also similar by using \eqref{alpha-}
if one considers the cases $\gamma+ 2s \leq 0$, $0 < \gamma +2s \leq -\alpha$ and
$\gamma+2s \geq -\alpha$ separately. Details are omitted.
And this completes the proof of Theorem 2.1.\\

In the following, we need also estimates on the commutator 
between the collision operator $Q$ and
 the weight $W_l$. For this purpose, estimations on 
$|W_l-W'_l|$ are needed.

\begin{lemm}\label{lemm2.2.11}
Let $l\in\NN$. There exists $C>0$ depending only on $l$ such that
\begin{equation}\label{E-3.105+3}
 \big|W_l-W'_l\big|\leq C\sin\left(\frac{\theta}{2}\right)\Big(W'_l+ W'_{l, *}\Big)
 \leq C\sin\left(\frac{\theta}{2}\right) W'_lW'_{l, *},
\end{equation}
and
\begin{equation}\label{E-3.105+3+1}
 \big|W_l-W'_l\big|\leq C\sin\left(\frac{\theta}{2}\right) \left(
 W'_l+ W'_{l-1} W'_{1, *}+ \sin^{l-1}\left(\frac{\theta}{2}\right) W'_{l,\, *}\right).
\end{equation}
\end{lemm}

\noindent{\bf Proof :} It follows from $|v-v_*|=|v'-v'_*|$ and
$|v|^2+|v_*|^2=|v'|^2+|v'_*|^2$ that, for any $\lambda>0$
$$
|v|^2\leq |v'|^2 + |v_*'|^2,\,\,\,\,\,\,W_\lambda \leq 2^\lambda
(W_{\lambda}' + W_{\lambda,\, *}').
$$
On the other hand
$$
|v-v'|^2=\sin^2\left(\frac{\theta}{2}\right)|v-v_*|^2,
$$
where $0\leq\theta\leq \pi/2$. Taylor formula yields
\begin{eqnarray*}
\big|W_l-W'_l\big|&\leq& C |v-v'|\Big(W_{l-1}+W'_{l-1}\Big)\\
&\leq&
C\sin\left(\frac{\theta}{2}\right)\, |v-v_*|\Big(W'_{l-1}+W'_{l-1, *}+W'_{l-1}\Big)\\
 &\leq& C\sin\left(\frac{\theta}{2}\right)\, |v'-v'_*|\Big(W'_{l-1}+W'_{l-1, *}\Big)\\
& \leq& C \sin\left(\frac{\theta}{2}\right)
\Big(W'_{1}+W'_{1,\, *}\Big)\Big(W'_{l-1}+W'_{l-1,\, *}\Big)\\
& \leq &C \sin\left(\frac{\theta}{2}\right) \big(W'_l+W'_{l, \,
*}\big) \leq C \sin\left(\frac{\theta}{2}\right) W'_lW'_{l,\, *} \, ,
\end{eqnarray*}
which gives \eqref{E-3.105+3}. For \eqref{E-3.105+3+1}, we have
\begin{eqnarray*}
\big|W_l-W'_l\big|&\leq& C |v-v'|\Big(W_{l-1}+W'_{l-1}\Big)\\
&\leq& C\sin\left(\frac{\theta}{2}\right)\, |v-v_*|
\left(W'_{l-1}+\Big(1+|v-v'+v'|^2\Big)^{\frac{(l-1)}{2}}\right)\\
 &\leq& C\sin\left(\frac{\theta}{2}\right)\,
 |v'-v'_*|\Big(W'_{l-1}+|v-v'|^{l-1}\Big)\\
& \leq& C \sin\left(\frac{\theta}{2}\right)
\left(\Big(W'_{1}+W'_{1,\,
*}\Big)W'_{l-1}+\sin^{l-1}
\left(\frac{\theta}{2}\right) |v'-v'_*|^{l}\right)\\
& \leq & C\sin\left(\frac{\theta}{2}\right) \left(
 W'_l+ W'_{l-1} W'_{1, *}+ \sin^{l-1}\left(\frac{\theta}{2}\right) W'_{l, *}\right)\, .
\end{eqnarray*}
And this completes the proof of the lemma.\\


\begin{lemm}\label{lemm3.4} 
Let $l \in \NN,\, m\in\RR$. 

$(1)$  \, If \, $0<s <1/2$, there exists $C >0$ such
that
\begin{equation}\label{3.11}
\left|\Big(\big(W_l  \,\, Q(f,\,\,g)-Q(f,\,\,W_l  \,\,
g)\big),\,\,\, h\Big)_{L^2(\RR^3_v)}\right| \leq C
\|f\|_{L^1_{l+\gamma^+}(\RR^3_v)} \|g\|_{L^{2}_{l+\gamma^+}
(\RR^3_v)} \|h\|_{L^{2}(\RR^3_v)}.
\end{equation}

Moreover, if $l \geq 3 $ (actually, we need only $l>\frac 32+2s$),
then
\begin{equation}\label{3.11-+}
\left|\Big(\big(W_l  \,\, Q(f,\,\,g)-Q(f,\,\,W_l  \,\,
g)\big),\,\,\, h\Big)_{L^2(\RR^3)}\right| \leq C
\|f\|_{L^2_{l+\gamma^+}(\RR^3_v))}\|g\|_{L^2_{l+\gamma^+}  (\RR^3)}
\|h\|_{L^2(\RR^3)}.
\end{equation}

$(2)$ \,  If $1/2 < s < 1$, then for any $\varepsilon>0$, there
is a constant $C_\varepsilon>0$ such that
\begin{eqnarray}\label{3.11+1}
\left|\Big(\big(W_l  \,\, Q(f,\,\,g)-Q(f,\,\,W_l  \,\,
g)\big),\,\,\, h\Big)_{L^2(\RR^3_v)}\right|  \hskip2.5cm \\
\hskip2cm \leq C_\varepsilon \|f\|_{
L^1_{l+2s-1+\gamma^+}(\RR^3_v)}\|g\|_{
H^{2s-1+\varepsilon}_{l+2s-1+\gamma^+}(\RR_v^3)}
\|h\|_{L^2(\RR^3_v)}\, , \nonumber
\end{eqnarray}
and
\begin{eqnarray}\label{3.11+1*}
\left|\Big(\big(W_l  \,\, Q(f,\,\,g)-Q(f,\,\,W_l  \,\,
g)\big),\,\,\, h\Big)_{L^2(\RR^3_v)}\right|  \hskip2.5cm \\
\hskip2cm \leq C_\varepsilon \|f\|_{
L^1_{l+2s-1+\gamma^+}(\RR^3_v)}\|g\|_{L^2_{l+2s-1+\gamma^+}(\RR_v^3)}
\|h\|_{H^{2s-1+\varepsilon}_{l}(\RR^3_v)}. \nonumber
\end{eqnarray}

$(3)$ \, 
When $s=1/2$, we have the same estimates  as $(2)$  with $2s-1$ replaced by any
small $\kappa >0$.
\end{lemm}

With Lemma 2.2,  we immediately have
the following improved upper bound estimate with respect to the
weight.

\begin{coro}\label{lemm104} $ $

$(1)$ When $0<s<1/2$, we have
$$
\|Q(f, g)\|_{H^{m}_l(\RR^3_v)}\leq C \|f\|_{L^{1}_{\max\{l+
\gamma^+, \,(\gamma+2s)^+ \}}(\RR^3_v)} \| g\|_{H^{m+2s}_{l+(2s+
\gamma)^+}(\RR^3_v)},
$$
provided that $m\leq 0$ and $0\leq m+2s$.

$(2)$ When  $1/2 <s <1$, we have
\begin{equation}\label{estimate-with-moment}
\|Q(f, g)\|_{H^{m}_l(\RR^3_v)}\leq C
\|f\|_{L^{1}_{\max\{l+2s-1+\gamma^+,\, (2s+ \gamma)^+\}}(\RR^3_v)}
\| g\| _{H^{m+2s}_{l+ \max\{2s-1+\gamma^+,\,(2s+
\gamma)^+\}}(\RR^3_v)},
\end{equation}
provided that $-1<m\leq 0$.

$(3)$ When $s=1/2$, we have the same form of estimate as
\eqref{estimate-with-moment} with $2s-1$ replaced by any small
$\kappa>0$.

\end{coro}

In fact, this corollary is a direct consequence of Theorem \ref{theo2.1} and
Lemma \ref{lemm3.4}.

\bigbreak \noindent
{\bf Proof of Lemma \ref{lemm3.4} :}
\smallbreak

{\bf Proof of $(1)$: the case $0<s<1/2$.}  By using $\Phi(|v'-v'_*|) \leq
\langle v' \rangle^{\gamma^+}\langle v'_* \rangle^{\gamma^+}$,
we have
\begin{eqnarray*}
&&\Big|\Big(\big(W_l\,\,Q(f,\,\,g)-Q(f,\,\,W_l\,\, g)\big),\,\,\,
h\Big)_{L^2(\RR^3_v)}\Big| \\
&=&\Big|\iiint
b\, \Phi \,f'_* g'(W'_l-W)\, h\, dv dv_* d\sigma\Big|\\
& \leq& C \iiint  b\, |\theta|\,\,
|(W_{l+ \gamma^+} f)'_*|\,\,|(W_{l+ \gamma^+}g)'|\,\,|h|\,\,dvdv_*d\sigma\\
& =&C\iiint  b\, |\theta| |(W_{l+ \gamma^+}f)_*|\,\,
|(W_{l+ \gamma^+}g)|\,\,|h'|\, dvdv_*d\sigma\\
& \leq& C \Big(\iiint  b\, |\theta|\,\, |(W_{l+ \gamma^+}\,
 f)_*|\,\,|(W_{l+ \gamma^+}\, g)|^2 dv dv_*
d\sigma \Big)^{1/2}\\
&&\times \Big(\iiint  b\, |\theta|\,\, |(W_{l+ \gamma^+}\, f)_*|\,
|h'|^2\, dvdv_*d\sigma\Big)^{1/2}\\
 &=&C J_1 \times J_2.
\end{eqnarray*}
Clearly, one has
\[
J^2_1\leq C \|f\|_{L^1_{l+ \gamma^+}}\|g\|^2_{L^2_{l+
\gamma^+}}\int_{\SS^2}b(\cos\theta)\,|\theta|\, d\sigma \leq C
\|f\|_{L^1_{l+ \gamma^+}}\|g\|^2_{L^2_{l+ \gamma^+}}.
\]
Next, by the regular change of variables $v\to v'$, cf. 
\cite{al-1, al-3}, we have
\[
J^2_2=\iint   D_0(v_*,v')|(W_{l+ \gamma^+}f)_*||h'|^2dv_*dv',
\]
where
$$
D_0(v,v')=2\int_{\SS^2} \frac{\theta(v_*, v',
\sigma)}{\cos^2(\theta(v_*, v', \sigma)/2)} b(\cos\theta(v_*,
v',\sigma))d\sigma \leq C\int_0^{\pi/4}\psi^{-1-2s}\sin\psi\,\,
d\psi,
$$
and
$$
\cos\psi=\frac{v'-v_*}{|v'-v_*|}\cdot \sigma, \quad \psi=\theta/2,
\qquad d\sigma=\sin\psi d\psi d\phi.
$$
Thus,
\[
J^2_2\leq C \|f\|_{L^1_{l+ \gamma^+}}\|h\|^2_{L^2},
\]
and this, together with the estimate on
$J_1$, gives (\ref{3.11}).

\smallbreak
We now prove \eqref{3.11-+} by using \eqref{E-3.105+3+1}
instead of \eqref{E-3.105+3}. We have
\begin{eqnarray*}
&&\Big|\Big(\big(W_l\,\,Q(f,\,\,g)-Q(f,\,\,W_l\,\, g)\big),\,\,\,
h\Big)_{L^2(\RR^3_v)}\Big| \\
& \leq& C \left\{ \iiint  b\, |\theta|^l\,\,
|(W_{l+ \gamma_+} f)'_*|\,\,|(W_{\gamma_+}g)'|\,\,|h|\,\,dvdv_*d\sigma\right.\\
&&+\iiint  b\, |\theta|\,\,
|(W_{1+\gamma^+} f)'_*|\,\,|(W_{l-1+ \gamma_+}g)'|\,\,|h|\,\,dvdv_*d\sigma\\
&& +\left. \iiint  b\, |\theta|\,\,
|(W_{\gamma_+} f)'_*|\,\,|(W_{l+ \gamma_+}g)'|\,\,|h|\,\,dvdv_*d\sigma \right\}\\
& =&\cM_1 + \cM_2 + \cM_3.
\end{eqnarray*}
$\cM_2, \cM_3$ can be estimated similarly to (\ref{3.11}), and we have
$$\cM_2 \leq C ||f||_{L^1_{1+\gamma^+}}\|g\|_{L^2_{l-1+\gamma^+}}\|h\|_{L^2},$$
$$\cM_3 \leq C ||f||_{L^1_{\gamma^+}}\|g\|_{L^2_{l+\gamma^+}}\|h\|_{L^2}.$$
 $\cM_1$ can be estimated as follows. Firstly, we have
\begin{align*}
\cM_1^2 =& C^2  \left(\iiint  b\, |\theta|^l
|(W_{l+ \gamma_+} f)_*||(W_{\gamma_+}g)|\,\,|h'|\,\,dvdv_*d\sigma\right)^2  \\
 \le &C^2  \iiint b \, |\theta|^{\, l-\frac{3}{2}}
 |(W_{\gamma_+}g)||(W_{l+ \gamma_+} f)_*|^2 dv dv_* d\sigma\\
& \times \iiint b \, |\theta|^{\,l+\frac{3}{2}}|(W_{\gamma_+}g)||h'|^2 dv dv_* d \sigma\\
= & \cM_{1,1} \times \cM_{1,2}.
\end{align*}
Then, if $l-\frac 32-2s-1>-1$, that is,  $l>2s+\frac 32$, we have
$$\cM_{1,1} \leq C  \|g\|_{L^1_{\gamma^+}} \|f\|_{L^2_{l+\gamma^+}}^2.$$
On the other hand, for $\cM_{1,2}$ we need to apply
the singular change of
variables $v_* \rightarrow v'$. The Jacobian of this transform
is, with $\vk=(v-v_*)/|v-v_*|$, 
\begin{equation}\label{jacobian2}
\Big|\frac{\pa v_*}{\pa v'}\Big|=\frac{8}{ \Big|I- \vk\otimes
\sigma\Big|}=\frac{8}{|1-\vk\cdot\sigma|} =\frac{4}{\sin^2
(\theta/2)}\le 16\theta^{-2}, \ \ \theta\in[0,\pi/2].
\end{equation}
Notice that this gives rise to an additional singularity
in the angle $\theta$ around $0$. Actually, the
situation is even worse in the following sense. 
Recall that $\theta$ is no longer
legitimate polar angle. In this case, the best choice of the pole is
$\vk''=(v'-v)/|v'-v|$ for which polar angle $\psi$ defined by
$\cos\psi=\vk''\cdot\sigma $  satisfies (cf. \cite[Fig. 1]{al-1})
\[
\psi=\frac{\pi-\theta}{2}, \qquad d\sigma=\sin\psi d\psi d\phi,
\qquad \psi\in[\frac{\pi}{4},\frac{\pi}{2}].
\]
This measure does not cancel any of the singularity of $b(\cos\theta)$,
unlike the case in the usual polar coordinates.
Nevertheless, this singular change of variables yields
\begin{align*}
\cM_{1,2} &= C \iiint  b\, |\theta|^{l+ \frac{3}{2}}
|(W_{\gamma_+}g)|\,\,|h'|^2\,\,dvdv_*d\sigma \\
&\leq C\iint D_1(v, v') |(W_{\gamma_+}g)|\,\,|h'|^2dv dv',
\end{align*}
when $l>\frac 32+2s$ because 
\[
D_1(v,v')=\int_{\SS^2}\theta^{l+\frac 32-2}b(\cos\theta)d\sigma \le
C \int_{\pi/4}^{\pi/2}(\frac{\pi}{2}-\psi)^{-2-2s+l+\frac 32-2}d\psi
\le C.
\]
Therefore,
\begin{equation*}
\cM_{1,2} \leq C \|g\|_{L^1_{\gamma^+}} \|h\|_{L^2}^2.
\end{equation*}
Now the proof of \eqref{3.11-+} is completed by the embedding estimate
for $l>\frac 32$,
$$
\|g\|_{L^1_{\gamma^+}}\leq C\|g\|_{L^2_{l+\gamma^+}}.
$$

\smallbreak 
{\bf Proof of $(2)$: the case 
$1/2 \leq s < 1$.}  Since we
look for an upper bound estimate and $\varepsilon >0$, 
it is sufficient to 
assume $s>1/2$ for our purpose. Write
\begin{eqnarray*}
&&\Big(\big(W_l\,\,Q(f,\,\,g)-Q(f,\,\,W_l\,\, g)\big),\,\,\,
h\Big)_{L^2(\RR^3_v)} =\iiint
B\, f'_* g'(W'_l-W_l)\, h\, dv dv_* d\sigma\\
&& = \iiint B\,\, f_* g(W_l-W'_l)\, h'\, dv dv_* d\sigma =\iiint
B \,\, f_* g'(W_l-W'_l)\, h'\, dv dv_* d\sigma \\
&&+\iiint B\,\, f_* (g-g')(W_l-W'_l)\, h'\, dv dv_* d\sigma
=I_1+I_2.
\end{eqnarray*}
 Taylor expansion gives 
$$
W_l-W'_l= \nabla W_l(v')(v-v')- \int_0^1 (1-\tau)\nabla^2
W_l(v'+\tau(v-v'))d\tau (v-v')^2\, ,
$$
so that
$$
I_1 = - \int_0^1 (1-\tau)\iiint
B \,\, f_* \{\nabla^2 W_l(v'+\tau(v-v'))\}(v-v')^2 g'\, h'\, dv dv_*
d\sigma d\tau\,.
$$
By using  the symmetry property shown in Figure
\ref{fig2} ( see also  Figure \ref{fig1} below, and \S3 in
\cite{HMUY}),  the first order term in the Taylor expansion 
vanishes, that is,
\begin{eqnarray*}
&& \iiint B \,\, f_* g'\nabla W_l(v')(v-v')\, h'\, dv dv_* d\sigma \\
&=& \iint f_*  \left \{ \int_{{\mathbb S}^{2}}
b\left(\frac{\psi_{\sigma}(v')-v_*}{ |\psi_{\sigma}(v')-v_*|}\cdot
\sigma\right) \Phi(|\psi_{\sigma}(v')-v_*|) \Big|\frac{\partial
(\psi_\sigma(v'))} {\partial (v')}\Big|
(\psi_{\sigma}(v')-v') d \sigma \right\}\\
&&\hspace{3cm}\cdot\,\,\,\nabla W_l(v')\,  g' h' dv'dv_* =0.
\end{eqnarray*}
Here, we have used the notation that
 for a transformation $v \rightarrow v'$, its
inverse transformation is denoted by $v' \rightarrow \psi_{\sigma}(v')=v$.
And $\sigma_1, \sigma_2$ are symmetric with respect to each other,
in the sense that
$\psi_{\sigma_1}(v')-v'=-(\psi_{\sigma_2}(v')-v')$.


Furthermore, since
$$
\begin{array}{lcl}
\displaystyle \left|\{\nabla^2 W_l(v'+\tau(v-v'))\}(v-v')^2 \right|
&\leq& \displaystyle C \theta^2 |v_* -v'|^2\{W_{l-2}(v_*)+
W_{l-2}(v'+\tau(v-v')-v_*)\}\\
&\leq &C \displaystyle \theta^2 \{W_{l}(v_*)+ W_{l}(v')\} \leq C \displaystyle \theta^2 W_{l}(v_*)W_{l}(v')
\end{array}
$$
and $\Phi(|v-v_*|) \leq (\sqrt 2 \langle v'-v_* \rangle)^{\gamma^+}
\leq \sqrt 2^{\gamma^+} \langle v_* \rangle^{\gamma^+} \langle
v'\rangle^{\gamma^+}$, we get by the regular change of variables $v\to v'$ and the Schwartz inequality
\begin{equation}\label{I-1}
|I_1| \leq C
||f||_{L^1_{l+\gamma^+}(\RR_v^3)}||g||_{L^2_{l+\gamma^+}(\RR_v^3)}||h||_{L^2(\RR^3_v)}.
\end{equation}

In order to estimate $I_2$, we shall apply the Littlewood-Paley
decomposition $\{\triangle_j\}_{j=0}^\infty$, which is a dyadic
decomposition in the Fourier variable (see also
\cite{bony,xu,alex}),
$$
\triangle_j g
(v)=\cF^{-1}\Big(\phi_j(\eta)\hat{g}(\eta)\Big),\,\,\,\,
g=\sum^\infty_0 \triangle_j g,
$$
and for $m\in\RR$,
$$
\|\triangle_j g\|_{H^m}\approx 2^{j m}\|\triangle_j
g\|_{L^2},\,\,\,\,\,\,\,\, \|g\|^2_{H^m}\approx \sum
2^{2j\,m}\|\triangle_j g\|^2_{L^2} .
$$
Then we have the following decomposition
$$
\begin{array}{lcl}
I_2 &=& \displaystyle \sum_{j=0}^\infty \int_0^1
\Big(\int_{\RR^6}\Big\{ \int_{\Omega_j} B\, f_*\nabla_v (\triangle_j
g) (v'+\tau(v-v'))(v-v')
(W_l-W'_l)\, h'\,d\sigma \Big\}dv dv_*   \Big) d\tau \\
&&+ \displaystyle \sum_{j=0}^\infty  \int_{\RR^6}\Big\{
 \int_{\Omega_j^c} B\, f_* \big\{(\triangle_j g)(v) - (\triangle_j g) (v')\big\}
(W_l-W'_l)\, h'\,d\sigma  \Big\} dv dv_*\\
&=& \displaystyle \sum_{j=0}^\infty \big(I_{2,j}^1+ I_{2,j}^2\big)\enskip,
\end{array}
$$
where
$$
\Omega_j= \Omega_j (v,v_*)= \left \{ \sigma \in {\mathbb S}^{2}\, ; \,
\frac{v-v_*}{|v-v_*|}\cdot \sigma \geq 1 - 2^{1-2j } \langle v-v_*
\rangle^{-2}\right \}.
$$
Note that if $1/2<s < 1$, then
\begin{eqnarray}\label{first}
\int_{\Omega_j} b(\cos \theta)\, \theta^2  d\sigma&=& 2\pi\int_{\{\theta
\in [0,\pi/2]; \,\sin(\theta/2) \leq 2^{-j} \langle
v-v_* \rangle^{-1}\}} \sin \theta\,\, b(\cos \theta) \theta^2 d\theta\\
&\leq& C 2^{j(2s-2)} \langle v-v_* \rangle^{2s-2},\nonumber
\end{eqnarray}
and
\begin{eqnarray}\label{second}
\int_{\Omega_j^c} b(\cos \theta)\, \theta\,  d\sigma&=& 2\pi
\int_{\{\theta \in [0,\pi/2]; \,\sin(\theta/2) \geq 2^{-j} \langle
v-v_* \rangle^{-1}\}} \sin \theta \,\,b(\cos \theta)\, \theta \,d\theta\\
&\leq& C 2^{j(2s-1)} \langle v-v_* \rangle^{2s-1}.\nonumber
\end{eqnarray}

To estimate $I_{2,j}^1$, we need the change of variables
\begin{equation}\label{change-z-v}
v\to
z=v'+\tau(v-v')=\frac{1+\tau}{2}v+\frac{1-\tau}{2}(|v-v_*|\sigma+v_*).
\end{equation}
The Jacobian of this transform
 is  bounded from below uniformly in $v_*, \ \sigma,
\ \tau$, because
\begin{align*}
\Big|\frac{\partial (z)}{\partial (v)}\Big|&=\Big|\text{det}\Big(
\frac{1+\tau}{2}I+\frac{1-\tau}{2}\sigma\otimes \vk\Big)\Big|
\quad\qquad
 (\vk=\frac{v-v_*}{|v-v_*|})
\\
& =\frac{(1+\tau)^3}{2^3} \Big|1+\frac{1-\tau}{1+\tau} \vk\cdot\sigma
\Big| = \frac{(1+\tau)^3}{2^3}
\Big|\frac{2\tau}{1+\tau}+2\frac{1-\tau}{1+\tau}\cos^2\frac{\theta}{2}
\Big|
\\
&\ge \frac{(1+\tau)^3}{2^3}\Big|\frac{2\tau}{1+\tau}+
\frac{1-\tau}{1+\tau} \Big|=\frac{(1+\tau)^3}{2^3}\ge \frac{1}{2^3}.
\end{align*}
Recall, cf. \cite{al-1}, that the cross-section $B(v-v_*,\theta)$ is assumed
to be supported in $0\le \theta\le \pi/4$.
Furthermore, we have
\begin{align}\label{z-v-equiv}
|z-v_*|&=\Big|\frac{1+\tau}{2}(v-v_*)+\frac{1-\tau}{2}|v-v_*|\sigma
\Big|
\\
& =|v-v_*|\Big|\Big(\frac{1+\tau}{2}\Big)^2+
\Big(\frac{1-\tau}{2}\Big)^2+\frac{1-\tau^2}{2}\vk\cdot\sigma
\Big|^{1/2} \notag
\\
& =|v-v_*|\Big|\tau^2+(1-\tau^2)\cos^2\frac{\theta}{2}\Big|^{1/2}
\ge \frac{1}{\sqrt{2}}|v-v_*|, \notag
\end{align}
which implies  $\langle v-v_* \rangle^{2s}
\Phi(|v-v_*|) \leq C \langle z \rangle^{2s+\gamma_+} \langle
v_* \rangle^{2s+ \gamma_+}$. Since 
$$
|(v-v') (W_l-W'_l)| \leq C \theta^2 |v-v_*|^2 (W_{l-1}(z) +
{W_*}_{l-1})\leq C \theta^2 |v-v_*|^2 W_{l-1}(z){W_*}_{l-1},
$$
we have from (\ref{first}) that for any $\varepsilon >0$
\begin{align*}
|I_{2,j}^1|&\leq C \int_0^1 \iiint_{\Omega_j} b \theta^2
|(W_{l-1+2s+\gamma_+}f)_*|\,\, |\big(W_{l-1+2s+\gamma_+}(\nabla_v
\triangle_j g)\big)(z)| \langle v-v_* \rangle^{2-2s}|h'|d
\sigma dv dv_* d \tau\\
&\leq C \left [ \int_0^1 \Big(\iint\Big( \int_{\Omega_j} b \theta^2
|(W_{l-1+2s+\gamma_+}f)_*|\,\, |\big(W_{l-1+2s+\gamma_+}(\nabla_v
\triangle_j g)\big)(z)|^2  \right.\\
& \hskip7cm \left.  \times \langle v-v_* \rangle^{2-2s}|d
\sigma \Big )dv dv_* \Big)^{1/2} d \tau \right ] \\
& \hskip3cm \times  \left[ \iint\Big( \int_{\Omega_j} b \theta^2
|(W_{l-1+2s+\gamma_+}f)_*|\,\,  \langle v-v_* \rangle^{2-2s}|h'|^2d
\sigma \Big) dv dv_* \right]^{1/2}\\
&\leq  C2^{-\varepsilon j}
\|f\|_{L^1_{l+2s-1+\gamma_+}(\RR^3_v)}\|g\|_{H^{2s-1 +
\varepsilon}_{l+2s-1+\gamma_+}(\RR_v^3)} \|h\|_{L^2(\RR_v^3)} ,
\end{align*}
where we used the regular change of variables $v \rightarrow  z$ defined by \eqref{change-z-v} and the regular change of variables
$v \rightarrow v'$.
The estimate (\ref{second}) yields the
same bound for $I_{2,j}^2$. Therefore, we obtain
\begin{equation}\label{I-2}
|I_2| \leq C ||f||_{L^1_{l+2s-1}(\RR_v^3)} ||
g||_{H^{2s-1+\varepsilon}_{l+2s-1}(\RR_v^3) } ||h||_{L^2(\RR^3_v)}.
\end{equation}
Estimates (\ref{I-1}) and (\ref{I-2}) together give the desired
estimate (\ref{3.11+1}).

For the convenience of the readers, we postpone the
 proof of (\ref{3.11+1*}) to the end of section \ref{s4+0}. And this
completes the proof of Lemma 2.2 because (3) comes from (2) for the case $s=1/2+ \kappa$.

\subsection{Coercivity estimates}\label{s4+00}
\setcounter{equation}{0}
\smallskip

We establish coercivity estimates of the Boltzmann collision
operator. We will show
that  the angular singularity in the cross-section yields
the  sub-elliptic estimates which are lower bounds of the collision operator, see \cite{al-1}.
Notice that we need 
precise weighted sub-elliptic estimates as given in the following
theorem. For more detailed explanations and notations,
interested readers can refer to
 \cite{alex,HMUY}. 

\begin{theo}\label{E-theo-0.1}
Assume that  $\gamma \in\R,\, 0<s<1$. Let $g\geq 0, \not \equiv 0, \, \,g\in
L^1_{\max\{\gamma^+,\, 2-\gamma^+\}}\bigcap L\log L (\RR^3_v)$. Then
there exists a constant $C_g >0$ depending only on $B(v-v_*,\theta)$,
$\|g\|_{L^1_{\max\{\gamma^+,\, 2-\gamma^+\}}}$ and $\|g\|_{L\, \log\, L}$,
and $C>0$ depending on $B(v-v_*,\theta)$ such that for any smooth function $f\in
H^{1}_{\gamma/2}(\RR^3_v)\cap L^2_{\gamma^+/2}(\RR^3_v)$, we have
\begin{equation}\label{E-sub-estimate+}
-\Big(Q(g,\, f),\,\, f\Big)_{L^2(\RR^3_v)}\geq C_g\|
W_{\gamma/2}f\|^2_{H^s(\RR^3_v)}- C ||g||_{L^1_{\max\{\gamma^+,\,
2-\gamma^+\}} (\RR^3_v)}\|f\|^2_{L^2_{\gamma^+/2}(\RR^3_v)}\,\, .
\end{equation}
\end{theo}

\begin{rema}\label{coercivity}
>From the proof of the theorem, the constant $C_g$ is seen to be an increasing function of 
$\| \tilde g \|_{L^1}$, $\| \tilde g \|_{L^1_1}^{-1}$ and $\| \tilde g \|_{L \log L}^{-1}$ where
$\tilde g =\la v \ra^{-|\gamma|} g$.
 If the function $g$ depends continuously on a parameter $\tau\in \Xi$,
then the constant $C_g$ depends on $\inf_{\tau\in\Xi}
\| \la v \ra^{-|\gamma|} g_\tau\|_{L^1}$, $\sup_{\tau\in\Xi}
\|g_\tau\|_{L\, \log\, L}$  and $\sup_{\tau\in\Xi}  \|g\|_{L^1_{\max\{\gamma^+,\, 2-\gamma^+\}}}$.
In the later application, we take $\tau =(t,x)$.
\end{rema}

\smallbreak \noindent
{\bf Proof.} Firstly,  we have
\begin{align*}
 (Q(g,f),f)&=\displaystyle
\int_{\R^{6}}\displaystyle\int_{\SS^{2}}\Phi(|v-v_*|)b(\cos
\theta)g(v_*) f(v)\{ f(v^\prime)
 - f(v)\}d\sigma d v_*dv  \\
&= {\displaystyle\frac{1}{2}}\displaystyle\int_{\R^{6}}\displaystyle
\int_{\SS^{2}}\Phi(|v-v_*|)b(\cos \theta)g(v_*)\{f(v^\prime)^2
 -f(v)^2\}d\sigma d v_*dv  \\
& -{\displaystyle\frac{1}{2}}\displaystyle\int_{\R^{6}}
\displaystyle\int_{\SS^{2}}\Phi(|v-v_*|)b(\cos
\theta)g(v_*)\{f(v^\prime)
 -f(v)\}^2d\sigma d v_*dv  \\
&= {\mathcal R}_1-{\mathcal R}_2.
\end{align*}
For ${\mathcal R}_1$, according to the cancellation lemma, Corollary 2 of
\cite{al-1}, we have
\begin{align*}
{\mathcal R}_1&={\displaystyle\frac{1}{2}}\displaystyle\int_{\R^{6}}
\displaystyle\int_{\SS^{2}}\Phi(|v-v_*|)b(\cos
\theta)g(v_*)\{f(v^\prime)^2
 -f(v)^2\}d\sigma d v_*dv\\
 &={\displaystyle\frac{1}{2}}\displaystyle\int_{\R^{6}}
 \displaystyle\int_{\SS^{2}}\left\{\Phi\left(\frac{|v-v_*|}{\cos
\frac{\theta}{2}}\right)\frac{1}{\cos^3
\frac{\theta}{2}}-\Phi(|v-v_*|)\right\}
b(\cos \theta)g(v_*)f(v)^2d v d\sigma d v_*\\
&={\displaystyle\frac{1}{2}}\displaystyle\int_{\R^{6}}
\displaystyle\int_{\SS^{2}}\Phi\left(\frac{|v-v_*|}{\cos
\frac{\theta}{2}}\right)\left\{\frac{1}{\cos^3
\frac{\theta}{2}}-1\right\} b(\cos
\theta)g(v_*)f(v)^2d v d\sigma d v_*\\
&+{\displaystyle\frac{1}{2}}\displaystyle\int_{\R^{6}}
\displaystyle\int_{\SS^{2}}\left\{\Phi\left(\frac{|v-v_*|}{\cos
\frac{\theta}{2}}\right)-\Phi(|v-v_*|)\right\} b(\cos
\theta)g(v_*)f(v)^2d v d\sigma d v_*\\
&={\mathcal R}_{11}+{\mathcal R}_{12}.
\end{align*}
For the first term ${\mathcal R}_{11}$, from $1- \cos^3
\frac{\theta}{2}\leq 3(1-\cos \frac{\theta}{2})=6 \sin ^2
\frac{\theta}{4} $, it follows that
\begin{equation*}
{\mathcal R}_{11}\leq
C\|g\|_{L^1_{\gamma^+}}\|f\|_{L^2_{\gamma^+/2}}^2,
\end{equation*}
because $\Phi \leq 1$ when $\gamma <0$. For the second term
${\mathcal R}_{12}$, we first note that
the mean value theorem gives
\begin{align*}
&\Phi\left(\frac{|v-v_*|}{\cos
\frac{\theta}{2}}\right)-\Phi(|v-v_*|)\\
&=-(\frac{1}{\cos
\frac{\theta}{2}}-1)|v-v_*|^2(1+(\frac{|v-v_*|}{a})^2)^{\frac{\gamma}{2}-1}\frac{2}{a^3}\\
&\leq  C (\frac{1}{\cos \frac{\theta}{2}}-1)\Phi(|v-v_*|),
\end{align*}
where $\frac{\sqrt{2}}{2}\leq \cos \frac{\theta}{2}<a<1$. Similar
to ${\mathcal R}_{11}$, we can obtain
\begin{equation*}
{\mathcal R}_{12}\leq \displaystyle
C\|g\|_{L^1_{\gamma^+}}\|f\|_{L^2_{\gamma^+/2}}^2.
\end{equation*}
For the term ${\mathcal R}_2$, we first note that
\begin{equation*}
\Phi(|v-v_*|)=(1+|v-v_*|^2)^{\frac{\gamma}{2}}\geq\displaystyle
 \frac{\langle v \rangle^\gamma}
{ \langle v_* \rangle^{|\gamma|}}.
\end{equation*}
Then, by using  the fact that $(a-b)^2 \geq a^2/2 - b^2$,  we have
\begin{align*}
{\mathcal R}_2&={\displaystyle\frac{1}{2}}\displaystyle\int_{\R^{6}}
\displaystyle\int_{\SS^{2}}\Phi(|v-v_*|)b(\cos
\theta)g(v_*)\{f(v^\prime)
 -f(v)\}^2d\sigma d v_*dv\\
 &\geq  C \displaystyle\int_{\R^{6}}\displaystyle\int_{\SS^{2}}b(\cos
\theta)\frac{g(v_*)}{\langle v_* \rangle^{|\gamma|}} \langle v
\rangle^\gamma \{f(v^\prime)
 -f(v)\}^2d\sigma d v_*dv\\
 &=  C
\displaystyle\int_{\R^{6}}\displaystyle\int_{\SS^{2}}b(\cos
\theta)\frac{g(v_*)}{\langle v_* \rangle^{|\gamma|}} \{\langle v
\rangle^{\frac{\gamma}{2}}f(v^\prime)
 -\langle v
\rangle^{\frac{\gamma}{2}} f(v)\}^2d\sigma d v_*dv\\
&= C \displaystyle\int_{\R^{6}}\displaystyle\int_{\SS^{2}}b(\cos
\theta)\frac{g(v_*)}{\langle v_* \rangle^{|\gamma|}} \{\langle
v^\prime \rangle^{\frac{\gamma}{2}}f(v^\prime)
 -\langle v
\rangle^{\frac{\gamma}{2}}f(v)\\
&\hspace{4.5cm}+\langle v
\rangle^{\frac{\gamma}{2}}f(v^\prime)-\langle v^\prime
\rangle^{\frac{\gamma}{2}}f(v^\prime)\}^2d\sigma d v_*dv\\
&\geq C_1
\displaystyle\int_{\R^{6}}\displaystyle\int_{\SS^{2}}b(\cos
\theta)\frac{g(v_*)}{\langle v_* \rangle^{|\gamma|}} \{\langle
v^\prime \rangle^{\frac{\gamma}{2}}f(v^\prime)
 -\langle v
\rangle^{\frac{\gamma}{2}}f(v)\}^2d\sigma d v_*dv\\
&-C_2\displaystyle\int_{\R^{6}}\displaystyle\int_{\SS^{2}}b(\cos
\theta)\frac{g(v_*)}{\langle v_* \rangle^{|\gamma|}}\{\langle v
\rangle^{\frac{\gamma}{2}}f(v^\prime)-\langle v^\prime
\rangle^{\frac{\gamma}{2}}f(v^\prime)\}^2d\sigma d v_*dv\\
&={\mathcal R}_{21}-{\mathcal R}_{22}.
\end{align*}
For the first term $ {\mathcal R}_{21}$, by using Corollary 3 and
Proposition 2 of \cite{al-1},  we have
\begin{equation}\label{5.22}
 \begin{split}
{\mathcal R}_{21}&=C_1
\displaystyle\int_{\R^{6}}\displaystyle\int_{\SS^{2}}b(\cos
\theta)\frac{g(v_*)}{\langle v_* \rangle^{|\gamma|}} \{\langle
v^\prime \rangle^{\frac{\gamma}{2}}f(v^\prime)
 -\langle v
\rangle^{\frac{\gamma}{2}}f(v)\}^2d\sigma d v_*dv
\\&
\ge C\int_{\RR^3}|\mathcal{F}(W_{\gamma/2}f)(\xi)|^2
\Big\{\int_{\SS^2}b(\tilde{\xi}\cdot \sigma)\Big(\cF(\tilde g)(0)-\Big|\cF(\tilde g)(\xi^-)\Big|
\Big)
\Big\}d\xi\\
 &\geq  \tilde{C}_{g}\|W_{\gamma/2}f\|_{H^s}^2-
C\|\tilde{g}\|_{L^1}\|f\|_{L^2_{\gamma^+/2}}^2,
 \end{split}
\end{equation}
 where
$\tilde g =\la v \ra^{-|\gamma|} g$.
Here $\tilde{C}_{g}$ 
is  an increasing function of 
$\| \tilde g \|_{L^1}$, $\| \tilde g \|_{L^1_1}^{-1}$ and $\| \tilde g \|_{L \log L}^{-1}$, 
according to
the proof in the last part of \cite{al-1} (see also Lemma 2.1 of
\cite{MUXY-DCDS}). 

For the second term $ {\mathcal R}_{22}$,
note that for some $\tau \in (0,1)$, we have
$$
\begin{array}{lll}\displaystyle
\frac{\langle v \rangle^{\frac{\gamma}{2}}-\langle v^\prime
\rangle^{\frac{\gamma}{2}}}{\langle v_*
\rangle^{\frac{|\gamma|}{2}}} &\leq& C \displaystyle  \frac{\langle
v'+ \tau(v-v') \rangle^{\frac{\gamma-2}{2}}}{\langle v_*
\rangle^{\frac{|\gamma|}{2}} }|v-v'|\\\\
&\leq& C \displaystyle \frac{\langle v_* \rangle^{\frac{|2-
\gamma|}{2} - \frac{|\gamma|}{2}}} {\langle v'+ \tau(v-v')-v_*
\rangle^{\frac{2-\gamma}{2}}}|v'-v_*| \tan(\theta/2)\\\\
&\leq& \displaystyle C \langle v_* \rangle^{\frac{|2-
\gamma|}{2} }  \langle
v'-v_*
\rangle^{\frac{\gamma}{2}}\tan(\theta/2)\\\\
&\leq & \displaystyle C   \left\{ \begin{array}{ll} \langle v_*
\rangle^{\frac{|2-\gamma|}{2}}
\langle v' \rangle^{\frac{\gamma}{2}} \tan(\theta/2),
& \mbox{if $  \gamma \geq 0$},\\
\langle v_* \rangle \tan(\theta/2),& \mbox{otherwise}. \end{array}
\right.
\end{array}$$
Hence, we get
\begin{align*}
{\mathcal R}_{22}&=C_2\displaystyle \int_{\R^{6}}\displaystyle
\int_{\SS^{2}}b(\cos \theta) g(v_*) \left \{ \frac{ \langle v
\rangle^{\frac{\gamma}{2}}-\langle v^\prime
\rangle^{\frac{\gamma}{2}}} {\langle v_*
\rangle^{\frac{|\gamma|}{2}}} \right \}^2  f(v^\prime)^2d\sigma d
v_*dv
\\
&\leq C_2\displaystyle\int_{\R^{6}}\displaystyle\int_{\SS^{2}}b(\cos
\theta)\tan^2(\theta/2) \langle v_* \rangle^{|2-\gamma^+|}
g(v_*)\{\langle v^\prime
\rangle^{\frac{\gamma^+}{2}}f(v^\prime)\}^2d
\sigma d v_*dv\\
&\leq  C_2 ||g||_{L^1_{|2-\gamma^+|} }\|f\|_{L^2_{\gamma^+/2}}^2.
\end{align*}
This completes the proof of Theorem \ref{E-theo-0.1}.

\bigbreak
In the following analysis, we
 shall also need the following interpolation inequality
between weighted Sobolev spaces in 
$v$, see for instance \cite{desv-wen1, HMUY}.

\begin{lemm}\label{lemm2.3}
For any $k\in\RR, p\in\RR_+, \delta>0$,
\begin{equation}\label{2.4}
\|f\|^2_{H^{k}_p(\RR^3_v)}\leq C_\delta \|f\|_{H^{k-\delta}_{2
p}(\RR^3_v)} \| f\|_{H^{k+\delta}_0(\RR^3_v)}.
\end{equation}
\end{lemm}

\vskip0.5cm
\subsection{Commutator estimates}\label{s4+0}
\setcounter{equation}{0}
\smallskip

We are now going to study the commutators of a family
of pseudo-differential operators with the 
Boltzmann collision operator. This is  a key step in the
regularity analysis of weak solutions because it requires the 
mollifiers defined by pseudo-differential operators. 
Below, we denote $(\cdot\,, \,\cdot)_{L^2(\RR_v^3)}$ by
$(\cdot\, ,\, \cdot)$ for simplicity of notations, without
any confusion.

\begin{prop}\label{coro-comm-v}
Let $\lambda \in \RR$ and  $M(\xi)$ be a positive symbol of
pseudo-differential operator in $S^{\lambda}_ {1,0}$ of the form of 
$M(\xi) = \tilde M(|\xi|^2)$.  Assume that for any $c>0$  there exists
a constant $C>0$ such that
 for any $s, \tau>0$
\begin{equation}\label{equivalence}
c^{-1}\leq \frac{s}{\tau}\leq c \,\,\,\,\,\,\mbox{implies}
\,\,\,\,\,\,\,\,C^{-1}\leq \frac{\tilde M(s)}{ \tilde M(\tau)}\leq
C.
\end{equation}
Furthermore assume that $M(\xi)$ satisfies
\begin{equation}\label{keep-form}
|M^{(\alpha)}(\xi)| = |\partial_\xi^\alpha M(\xi)| \leq C_{\alpha}
M(\xi) \la \xi \ra^{-|\alpha|},
\end{equation}
for any $\alpha\in\NN^3$.  Then the followings hold.

$(1)$ \,  If $0<s<1/2$, for any
$N>0$ there exists a $C_N  >0$ such that 
\begin{eqnarray}\label{10.8-2}
\left|( M(D_v) Q(f,\, g)-  Q(f,\, M(D_v) g) ,\,\, h)_{L^2(\RR^3_v)}\right | \hskip4cm\\
\hskip3cm \leq C_N \|f\|_{L^1_{\gamma^+}(\RR^3_v))} \Big(\|M
g\|_{L^2_{\gamma^+}(\RR^3_v)} + \| g \|_{H^{\lambda-N}_{\gamma^+}(\RR^3_v))} \Big)
\|h\|_{L^2(\RR_v^3)}.\nonumber
\end{eqnarray}

$(2)$ \,  If $1/2 < s <1$, for any
$N>0$ and any $\varepsilon >0$ there exists a $C_{N, \varepsilon}>0$ such that 
\begin{align}\label{10.8-3}
\left |(M(D_v) Q(f,\, g)-Q(f,\, M(D_v) g),\,\, h)_{L^2(\RR^3_v)} \right |  \hskip4cm \\
\leq C_{N, \varepsilon} \|f\|_{L^1_{(2s+ \gamma-1)^+}(\RR^3_v))} \Big(
\|M g\|_{H^{2s-1+\varepsilon} _{(2s+ \gamma-1)^+}(\RR^3_v))}  +
\|g\|_{H^{\lambda-N}_{\gamma^+}(\RR^3_v))}\Big)
 \|h\|_{L^2(\RR_v^3)}\, . \nonumber
\end{align}

$(3)$ \,
If $s = 1/2$, we have the same estimate  as (\ref{10.8-3}) with
$(2s+ \gamma-1)$ replaced by $(\gamma+ \kappa)$ for any small
$\kappa >0$.
\end{prop}

\smallbreak\noindent
{\bf Proof :}
Firstly,  set $\Phi_*(v)=\Phi(|v-v_*|)$ and write
\begin{eqnarray*}
&&\Big(M(D_v)Q(f,\, g),\,\, h\Big)-\Big(Q(f,\, M(D_v)g),\,\, h\Big)\\
&=&\displaystyle\int_{\RR^{6}}
\displaystyle\int_{\SS^{2}} B(|v-v_*|,\sigma)f(v_*)g(v)\Big(
\overline{\big(M\, h\big)(v')}- \overline{\big(M\, h\big)(v)}
\Big) d\sigma d v_* d v
\\
&&-\displaystyle\int_{\RR^{6}}
\displaystyle\int_{\SS^{2}} B(|v-v_*|,\sigma)f(v_*)\big(M\, g\big)(v)\Big(
 \overline{h(v')}-  \overline{h(v)}\Big) d\sigma d v_* d v\\
&=&\displaystyle\int_{\RR^{6}}
\displaystyle\int_{\SS^{2}} b(\cos \theta)f(v_*) \Big [ (\Phi_* g)(v)
\overline{(M h)(v')}-   \big\{M(\Phi_* g)\big\}(v)\overline{h(v')}
\Big] d\sigma d v_* d v
\\
&&+\displaystyle\int_{\RR^{6}}
\displaystyle\int_{\SS^{2}} b(\cos \theta) f(v_*) \big\{M(\Phi_* g)\big\}(v)\Big( \overline{h(v')}
- \overline{h(v)} 
\Big) d\sigma d v_* d v\\
&&-\displaystyle\int_{\RR^{6}}
\displaystyle\int_{\SS^{2}}b(\cos \theta) f(v_*)\big\{ \Phi_*\big(M\, g\big)\big\}(v)\Big(
\overline{ h(v')}-  \overline{h(v)}\Big) d\sigma d v_* d v\\
&=&\displaystyle\int_{\RR^{6}}
\displaystyle\int_{\SS^{2}} b(\cos \theta)f(v_*) \Big [ (\Phi_* g)(v)
\overline{(M h)(v')}-   \big\{M(\Phi_* g)\big\}(v)\overline{h(v')}
\Big] d\sigma d v_* d v
\\
&&+\displaystyle\int_{\RR^{6}}
\displaystyle\int_{\SS^{2}} b(\cos \theta)f(v_*)\Big( \big[M\,,\,\, \Phi_* \big]g\Big)(v)\big( \overline{h(v')}- \overline{h(v)}\big) d\sigma d v_* d v\\
&=&{\mathcal I}+{\mathcal {II}}.
\end{eqnarray*}
The above computation is justified with cutoff approximation, see the remark given after \eqref{Qfourier} and also 
\cite{HMUY}.
The first term ${\mathcal I}$  can be rewritten by using
Bobylev formula (see e.g. \cite{al-1}) as
\begin{align*}
&{\mathcal I}=
\displaystyle\int_{\RR^{6}}\displaystyle\int_{\SS^{2}}
b(\frac{\xi}{|\xi|}\cdot\sigma)f(v_*)
 \Big(M(\xi)-M(\xi^+)\Big) {\mathcal
F}(\Phi_* g)(\xi^{+})e^{-iv_{*}\xi^{-}}d\sigma d v_*
\overline{\hat{h}(\xi)}d \xi,
\end{align*}
where 
\begin{align*}
\qquad \xi^{\pm}=\frac{\xi\pm
|\xi|\sigma}{2}.
\end{align*}
Notice that in the case of
Maxwellian molecule type cross section with
$\gamma=0$ i.e. $\Phi(|v-v_*|)=1$,
${\mathcal {II}}\equiv0$.

Since ${\tilde M}'(|\xi|^2) = 2\xi \cdot \nabla M(\xi)/|\xi|^2$ and
$|\xi^+|\leq |\xi|\leq 2|\xi^+|$, it follows from
\eqref{equivalence} and \eqref{keep-form} that  
\begin{equation}\label{inegalite}
|M(\xi) - M(\xi^+)| \leq C \left| \sin
\frac{\theta}{2}\right|^2 M(\xi^+),
\end{equation}
and
$$
\int_{\SS^{2}}
b\Big(\frac{\xi}{|\xi|}\cdot\sigma\Big)\left|\sin\frac{\theta}{2}\right|^2
d\sigma\leq C<+\infty.
$$
Thus,
\begin{align*}
|{\mathcal I}|&\le C\displaystyle\int_{\R^{3}}\la v_*
\ra^{\gamma_+} |f(v_*)|\displaystyle\int_{\SS^{2}} \displaystyle
\int_{\R^{3}}b(\frac{\xi}{|\xi|}\cdot\sigma)\sin ^2\frac{\theta}{2}
M(\xi^+)|
{\mathcal F}(\Phi_* \la v_* \ra^{-\gamma_+} g)(\xi^{+})|\,\,
|\overline{\hat{h}(\xi)}|d \xi  d\sigma d v_*
\\
& \leq C\displaystyle \Big ( \int_{\R^{3}}| \la v_* \ra^{\gamma_+}
f(v_*)|\displaystyle\int_{\SS^{2}} \displaystyle
\int_{\R^{3}}b(\frac{\xi}{|\xi|}\cdot\sigma)\sin ^2\frac{\theta}{2}
|M(\xi^+) {\mathcal F}(\Phi_* \la v_* \ra^{-\gamma_+} g)(\xi^{+})| ^2 d \xi
d\sigma d v_* \Big)^{1/2}
\\
& \hskip0.5cm \times \Big ( \int_{\R^{3}}| \la v_* \ra^{\gamma_+}
f(v_*)|\displaystyle\int_{\SS^{2}} \displaystyle
\int_{\R^{3}}b(\frac{\xi}{|\xi|}\cdot\sigma)\sin
^2\frac{\theta}{2}|\overline{\hat{h}(\xi)}|^2 d \xi  d\sigma d
v_*\Big)^{1/2}
 \\
& \leq C\|f\|_{L^1_{\gamma_+}} \Big ( \sup_{v_*} \|  M(D_v)\Phi_* \la v_*
\ra^{-\gamma_+} g(v) \|_{L^2_{\gamma_+} } \Big ) \|h\|_{L^2 },
\end{align*}
where we have used Plancherel's equality, the change of variables
$\xi\to\xi^+$ for which $d\xi\sim d\xi^+$ uniformly with respect to $\sigma$,
the estimate $\Phi_* \la v_* \ra^{-\gamma_+} \leq \la v \ra^{\gamma_+} $.
Then by using the expansion formula of the pseudo-differential calculus
\begin{equation}\label{expansion-formula}
[M(D_v),\, \Phi_*(v)]\, g
 = \sum_{1 \leq |\alpha| < N_1} \frac{1}{\alpha!}
\Phi_{*(\alpha)} M^{(\alpha)}(D_v)g +  r_{N_1} (v,D_v;v_*)g\, ,
\end{equation}
with $N_1>\lambda$, and the condition \eqref{keep-form}, we obtain
\begin{equation}\label{major-estimate}
\sup_{v_*} \|  M(D_v) \Phi_* \la v_* \ra^{-\gamma_+} g(v)
\|_{L^2_{\gamma_+} }\leq C \Big (\|M g\|_{L^2_{\gamma_+} } +
||g||_{H^{\lambda-N_1}_{\gamma^+}} \Big) .
\end{equation}
Hence,
\begin{equation}\label{I11}
|{\mathcal I}|\leq C
\|f\|_{L^1_{\gamma_+}} \Big (\|M g\|_{L^2_{\gamma_+} } +
||g||_{H^{\lambda-N_1}_{\gamma^+}} \Big)\|h\|_{L^2} \, .
\end{equation}

We now turn to 
 the term ${\mathcal {II}}$. Firstly,
set
\begin{equation*}
F(v,v_*)=[M,\,\,\Phi_*]g(v),
\end{equation*}
and decompose
\begin{align*}
{\mathcal {II}} &=
\displaystyle\int_{\R^{6}}\displaystyle\int_{\SS^{2}} b(\cos
\theta)f(v_*)\Big\{F(v^\prime,v_*)h(v^\prime)- F(v,v_*)h(v)\Big\}d\sigma d
v_* d v
\\
&\hspace*{0.3cm}+\displaystyle\int_{\R^{6}}\int_{\SS^{2}} b(\cos
\theta)f(v_*)\Big(F(v,v_*)-F(v^\prime,v_*)\Big)h(v') d v_* d vd\sigma
\\
&=J_1+J_2.
\end{align*}
According to the cancellation lemma \cite{al-1}, we obtain
\begin{equation*}
\int_{\R^{3}}\displaystyle\int_{\SS^{2}} b(\cos
\theta)\Big\{F(v^\prime,\,v_*)h(v^\prime)- F(v,v_*)h(v) \Big\}d\sigma d
v=\Big(S*\Big\{F(\,\cdot\,,\,v_*)h\Big\}\Big)(v_*),
\end{equation*}
where the convolution product is in $v\in\RR^3$, and in this case,
\begin{equation*}
S= 2\pi\int_0^{\pi/2}\sin \theta
b(\cos\theta)\Big[\frac{1}{\cos^3(\theta/2)}-1\Big]d\theta
\end{equation*}
is a constant function. Consequently,
\begin{equation*}
J_1=\int_{\R^3}f(v_*)\Big(S*\{F(\cdot,v_*)h\}\Big)
(v_*)dv_*=S\int_{\R^{6}}f(v_*)
F(v,v_*)h(v)dvdv_*.
\end{equation*}
By (\ref{expansion-formula}) and (\ref{major-estimate}), we get
\begin{align}\label{J1}
|J_1|\le&
C\int_{\R^3}|f(v_*)|\|F(\cdot,v_*)\|_{L^2}\|h\|_{L^2}dv_*  \\
\le& C
\|f\|_{L^1_{\gamma_+}} \Big (\|M g\|_{L^2_{\gamma_+} } +
||g||_{H^{\lambda-N_1}_{\gamma^+}} \Big)\|h\|_{L^2} \,. \nonumber
\end{align}

To estimate the term $J_2$, we need to consider the following two
cases.

\smallbreak
\noindent{\bf Case  1: $0<s  <1/2$ .}  Since the mean
value theorem yields
\begin{equation*}
F(v,v_*)-F(v',v_*)=(v-v')\cdot
\int_0^1\nabla_v(F(v'+\tau(v-v'),v_*)d\tau,
\end{equation*}
by noticing that
\begin{equation*}
|v'-v|=|v-v_*|\sin(\theta/2)=|v'-v_*|\tan(\theta/2),
\end{equation*}
we have
\begin{align*}
&|J_{2}|
\le \int_0^1\Big(\int_{\R^{6}\times\SS^{2}} b(\cos\theta)|v'-v|
|f(v_*)||h(v')||(\nabla_vF)(v'+\tau(v-v'),v_*)|
 dvdv_*d\sigma\Big)d\tau.
 \\
 \le & C
\Big( \int_{\R^{6}\times\SS^{2}} b(\cos\theta)|\theta|  \la v_*
\ra^{\gamma_+} |f(v_*)||h(v')|^2
 dvdv_*d\sigma\Big)^{1/2}   \\
 & \times \int_0^1
\Big( \int_{\R^{6}\times\SS^{2}} b(\cos\theta)|\theta|
\la v_* \ra^{\gamma_+}|f(v_*)|
\left| \frac{|v-v_*|}{\la v_*
\ra^{\gamma_+}}(\nabla_vF)(v'+\tau(v-v'),v_*) \right |^2
 dvdv_*d\sigma\Big) ^{1/2 }d\tau
\\
 =&C \, J_{21} \times J_{22}.
 \end{align*}
By the change of variables $v\to v'$ for which $dv\sim dv'$
uniformly in $v_*\in\RR^3,\,\sigma\in\SS^2$ (see \cite{al-1}), we get
\begin{equation}\label{J21}
J_{21}^2 \le C\|f\|_{L^1_{\gamma_+}}\|h\|^2_{L^2}.
\end{equation}

To estimate $J_{22}$, we apply the change of variables \eqref{change-z-v}
and use \eqref{z-v-equiv}.
Setting
\begin{equation*}
\psi^*(v)=\displaystyle \frac{\langle v-v_* \rangle}{ \langle v_*
\rangle^{\gamma_+}},
\end{equation*}
   we get
\begin{align*}
J_{22}^2 &\le C \int_0^1\Big[\int_{\R^{2n}\times\SS^{2}}
b(\cos\theta)|\theta| \langle v_* \rangle^{\gamma_+}
|f(v_*)|\Big|\psi^*(z)(\nabla_vF)(z,v_*) \Big|^2
 dzdv_*d\sigma\Big]d\tau
\\
& \le C\|f\|_{L^1_{\gamma_+}}\sup_{v_*} \| \psi^*(\cdot)
 (\nabla_vF)(\cdot,v_*)\|^2_{L^2}.
\end{align*}

On the other hand, it follows from the expansion formula of
pseudo-differential operators that, with
$\Phi_*(v)=(1+|v-v_*|^2)^{\gamma/2}$ we have for any $N_1 \in \NN$
\begin{equation}\label{expansion-1}
 \begin{split}
&(\nabla_v F)(v,v_*) =  \nabla_v[M,\, \Phi_*]g (v)\\
 =& \sum_{1 \leq |\alpha| < N_1} \frac{1}{\alpha!}\left \{  (\nabla
\Phi_{*(\alpha)}) M^{(\alpha)}(D_v)g+
\Phi_{*(\alpha)} M^{(\alpha)}(D_v) \nabla_v g\right \}
+  \tilde{r}_{N_1} (v,D_v;v_*)g \\
 =& F_{N_1}(v,D_v;v_*)g(v) + \tilde{r}_{N_1}(v,D_v;v_*)g(v),
\end{split}
\end{equation}
where $\tilde{r}_{N_1}$ is a pseudo-differential operator with
symbol belonging to $S^{\lambda-N_1}_{1,0}$ uniformly with respect to $v_*\in\RR^3$
(cf.  \cite{Kuma}).
Since
\begin{equation*}
\big|\psi^* \Phi_{*(\alpha)}\big | \leq C_{\alpha} \frac{\langle
v-v_* \rangle}{ \langle v_* \rangle^{\gamma_+}}\langle v-v_*
\rangle^{\gamma- |\alpha|} \leq C_{\alpha} \langle
v\rangle^{\gamma_+},
\end{equation*}
by (\ref{keep-form}), we have for $\alpha \ne 0$ that,
$$
|M^{(\alpha)}(\xi) \,\,\xi|\leq C_\alpha M(\xi) \la \xi \ra^{-|\alpha|+1}\leq
 C_\alpha M(\xi).
$$
Hence
\begin{equation}\label{J22}
J_{22}^2 \leq
C
\|f\|_{L^1_{\gamma_+}} \Big (\|M g\|^2_{L^2_{\gamma_+} } +
||g||^2_{H^{\lambda-N_1}_{\gamma^+}} \Big).
\end{equation}
Now, it follows from \eqref{J1}, \eqref{J21}, and \eqref{J22} that
\begin{equation}\label{CU}
|{\mathcal {II}} |\le
 C\|f\|_{L^1_{\gamma_+}} \Big (\|M g\|_{L^2_{\gamma_+} } +
||g||_{H^{\lambda-N_1}_{\gamma^+}}\Big)||h||_{L^2}.,
\end{equation}
holds when  $0<s <1/2$.

\smallbreak
\noindent{\bf Case 2: $1/2<s
<1$.} We now decompose $J_2$  as follows:
\begin{align*}
J_{2}&= \int_0^1\Big(\int_{\R^{6}\times\SS^{2}} b(\cos\theta)
f(v_*)h(v')(v-v')\cdot (\nabla_vF)(v'+\tau(v-v'),v_*)
 dvdv_*d\sigma\Big)d\tau
 \\
 &
 = \int_{\R^{6}\times\SS^{2}} b(\cos\theta)
f(v_*)h(v')(v-v')\cdot (\nabla_vF)(v',v_*)
 dvdv_*d\sigma
\\
&+ \int_0^1\Big(\int_{\R^{6}\times\SS^{2}} b(\cos\theta)
f(v_*)h(v')\\
&\hskip2.2cm (v-v')\cdot \left\{(\nabla_vF)(v'+\tau(v-v'),v_*) -
(\nabla_vF)(v',v_*) \right \}
 dvdv_*d\sigma\Big)d\tau
\\
 &=J^0_2+J^1_2.
\end{align*}
The essential feature of this decomposition is that $J^0_2$ vanishes
by symmetry as in the proof of Lemma \ref{lemm3.4}.
Indeed, we have
\begin{align*}
J^0_2 &=  \int_{\R^{6}}
f(v_*)h(v')\\
&\enskip \left \{ \int_{\SS^{2}}
b\left(\frac{\psi_{\sigma}(v')-v_*}{ |\psi_{\sigma}(v')-v_*|}\cdot
\sigma\right) \Big|\frac{\partial (\psi_\sigma(v'))} {\partial
(v')}\Big|
(\psi_{\sigma}(v')-v') d \sigma \right\}\cdot (\nabla_vF)(v',v_*) dv'dv_* \\
&=0, \end{align*} because of the symmetry
in  $\sigma_1$ and $ \sigma_2$  in the sense that
$\psi_{\sigma_1}(v')-v'=-(\psi_{\sigma_2}(v')-v')$, cf.
Figure \ref{fig1}.

Now, by the change of variable $v \rightarrow z = v'+ \tau(v-v')$ defined by \eqref{change-z-v}, we
consider
\begin{equation*} 
J_2^1(\tau)=\int_{\R^{6}\times\SS^{2}} b f(v_*)h(v') (v-v')\cdot
\{(\nabla_v F)(z,v_*) - (\nabla_v F)(v',v_*)\}
 dvdv_*d\sigma.
\end{equation*}
By recalling the expansion formula \eqref{expansion-1} of $(\nabla_v
F)(v,v_*)$, we first consider
\begin{equation}
 \begin{split}
&J_2^1(\tau,\alpha)\\
=&\int_{\R^{6}\times\SS^{2}} b f(v_*)h(v') (v-v')\cdot
\{\Phi_{*(\alpha)}M^{(\alpha)} \nabla_v
g(z)-\Phi_{*(\alpha)} M^{(\alpha)} \nabla_v g(v')\}
 dvdv_*d\sigma  \\
=&\int_{\R^{6}\times\SS^{2}} b f(v_*)h(v')\left\{
\Phi_{*(\alpha)}(z) -
\Phi_{*(\alpha)}(v')\right\}(v-v')\cdot M^{(\alpha)} \nabla_v
g(z)
 dvdv_*d\sigma  \\
&+\int_{\R^{6}\times\SS^{2}} b f(v_*)h(v')\Phi_{*(\alpha)}(v')
(v-v')\cdot \{M^{(\alpha)}\nabla_v g(z)-
M^{(\alpha)} \nabla_v g(v')\}
 dvdv_*d\sigma \\
=& J_{2}^{1,0}(\tau,\alpha) +
\tilde{J}_2^1(\tau,\alpha).\label{decompo}
\end{split}
\end{equation}
Notice that the case when
 $|\alpha|=1$  is the most difficult one, in the sense that
$M^{(\alpha)}(D_v)\,\nabla_v$ is a pseudo-differential operator of
order $\lambda$ with symbol bounded by $C\,M(\xi)$ 
due to the assumption  (\ref{keep-form}).
By writing $(1)$ instead of $(\alpha)$ when $|\alpha|=1$, we
have
$$\left|\left \{\Phi_{*(1)}(z) - \Phi_{*(1)}(v')\right\} |v-v' |\right| \leq
C \la z-v_* \ra^\gamma \theta^2,$$ which gives
\begin{align}\label{J-2-10}
\left|J_{2}^{1,0}(\tau,(1)) \right| \leq&
\left(\int_{\R^{6}\times\SS^{2}} b \theta^2 | \la v_* \ra^{\gamma_+}
f(v_*)|
|h(v')|^2 d\sigma dvdv_* \right)^{1/2}\\
&\times \left( \int_{\R^{6}\times\SS^{2}} b \theta^2 | \la v_*
\ra^{\gamma_+} f(v_*)| \left | \frac{ \la z-v_* \ra^\gamma}{
\la v_* \ra^{\gamma_+} } M^{(\alpha)}\, \nabla_v g(z)
\right |^2
d\sigma dvdv_* \right )^{1/2}   \nonumber \\
\leq& C||f||_{L^1_{\gamma_+}} ||M\,g||_{L^2_{\gamma_+} }||h||_{L^2}.
\nonumber
\end{align}

In order to evaluate the term $\tilde{J}_2^1(\tau,(1)) $,  we take
the same Littlewood-Paley partition of unity $\{\psi_j(\xi)\}$ as in
the proof of Lemma \ref{lemm3.4} and write
\begin{align*}
&\tilde{J}_2^1(\tau,(1)) \\
=& \int_{\R^{6}\times\SS^{2}} b f(v_*)h(v')\Phi_{*(1)}(v')
(v-v')\cdot \{M^{(1)}\nabla_v g(z)- M^{(1)}
\nabla_v g(v')\}
 dvdv_*d\sigma \\
 =& \sum_{j=0}^\infty \int_{\R^{6}\times\SS^{2}}
 b f(v_*)h(v') \Phi_{*(1)}(v')
(v-v')\cdot\left( g_j(z) - g_j(v') \right) dvdv_*d\sigma\\
=& \sum_{j=0}^\infty \tilde{J}_{2,j}^1(\tau),
\end{align*}
where $g_j(v)= \psi_j(D_v)M^{(1)}(D_v)\nabla_v g(v)$. For each
$j$ we apply the 
following decomposition by using $\Omega_j$ introduced in the proof of
Lemma \ref{lemm3.4} to have
\begin{align*}
&\tilde{J}_{2,j}^1(\tau)  \\
=&\enskip \int_0^1 \Big(  \int_{\R^{6}} \Big(\int_{\Omega_j} b
f(v_*) h(v') \Phi_{*(1)}(v')(v-v')\\
&\hspace{3.5cm}\cdot (z - v') \nabla
g_j(v'+ s(z -v')) d \sigma \Big)  dvdv_* \Big ) ds \\
&\enskip + \int_{\R^{6}} \left(\int_{\Omega_j^c} b f(v_*) h(v')
\Phi_{*(1)}(v')(v-v')\cdot \left(
g_j(z) - g_j(v')\right) d \sigma \right)  dvdv_*\\
\enskip =& \tilde{J}_{2,j}^{1,1}(\tau) +
\tilde{J}_{2,j}^{1,2}(\tau).
\end{align*}
By setting $$v'_{\tau,s} =v'+ s(z -v'),$$ we have
\begin{align*}
|\tilde{J}_{2,j}^{1,1}(\tau)|&\le \int_0^1\Big(\int_{\R^{6}}
\left(\int_{\Omega_j}b|v'-v|^2 |f(v_*)||h(v')| |\Phi_{*(1)}(v')|
|\nabla g_j(v'_{\tau,s})|d\sigma\right) dvdv_*\Big)ds
 \\
 &
 \le
 C \int_0^1 \Big (\int_{\R^{6}} \Big(\int_{\Omega_j}
 b(\cos\theta)\theta^2\la v- v_* \ra^{2-2s}
\la v_* \ra^{(2s+ \gamma-1)_+} |f(v_*)|\\
&\hskip3cm \times
 |h(v')|
\left| \frac{\la v'_{\tau,s} -v_* \ra^{2s+ \gamma -1}}{\la v_*
\ra^{(2s+ \gamma-1)_+}} \nabla g_j(v'_{\tau,s}) \right|  d\sigma
\Big) dvdv_*\Big)ds
 \\
 &\leq
C   2^{-\varepsilon j}\Big (\int_{\R^{6}} \Big(\int_{\Omega_j}
b(\cos\theta)\theta^2 2^{j(2-2s)} \la v- v_* \ra^{2-2s} \la v_*
\ra^{(2s+ \gamma-1)_+} |f(v_*)|
 |h(v')|^2d\sigma   \Big) dvdv_* \Big)^{1/2}\\
 & \hskip1cm \times \int_0^1 \Big (
\Big (\int_{\R^{6}} \Big(\int_{\Omega_j} b(\cos\theta)\theta^2
2^{j(2-2s)} \la v- v_* \ra^{2-2s}
\la v_* \ra^{(2s+ \gamma-1)_+} |f(v_*)|\\
& \hskip4cm \left| \la v'_{\tau,s} \ra^{(2s+ \gamma-1)_+} 2^{j(2s-2
+ \varepsilon)}\nabla  g_j(v'_{\tau,s}) \right|^2 d\sigma   \Big)
dvdv_* \Big)^{1/2} \Big)ds
\\&
 =C 2^{-\varepsilon j}\tilde{J}_{2,j,1}^{1,1}(\tau)\times \tilde{J}_{2,j,2}^{1,1}(\tau) .
\end{align*}
By using the same change of variables as for $J_{2,1}$ in the
previous case, it follows from \eqref{first} that
\begin{equation}\label{J-1-1-1}
\tilde{J}_{2,j,1}^{1,1}(\tau)^2 \leq C ||f||_{L^1_{(2s
+\gamma-1)_+}}||h||_{L^2}^2.
\end{equation}
Similarly, by taking the 
 change of variables $v \to v'_{\tau,s}$ as in the
previous  case again, \eqref{first} leads to
\begin{equation}\label{J-1-1-2}
\tilde{J}_{2,j,2}^{1,1}(\tau)^2 \leq C ||f||_{L^1_{(2s
+\gamma-1)_+}} \Big(||M\,g||^2_{H^{2s-1+ \varepsilon} _{(2s + \gamma-1)_+} }
+||g||^2_{H^{\lambda-N_1+2s-1+ \varepsilon} _{(2s + \gamma-1)_+} }\Big),
\end{equation}
where we have used
\begin{equation*}
||2^{j(2s-2-\varepsilon)}\nabla g_j(v)||^2_{L^2_{(2s+ \gamma -1)_+}}
\leq C\Big(||M\,g||^2_{H^{2s-1+ \varepsilon} _{(2s + \gamma-1)_+} }
+||g||^2_{H^{\lambda-N_1+2s-1+ \varepsilon} _{(2s + \gamma-1)_+} }\Big).
\end{equation*}
Hence, it follows from \eqref{J-1-1-1} and \eqref{J-1-1-2} that, for $N_1>\lambda
+2s-1+ \varepsilon$, we have
\begin{equation}\label{J-2-j-1}
|\tilde{J}_{2,j}^{1,1}(\tau)| \leq C 2^{-\varepsilon j}
||f||_{L^1_{(2s +\gamma-1)_+}}  \Big(||M\,g||_{H^{2s-1+ \varepsilon} _{(2s + \gamma-1)_+} }
+||g||_{H^{\lambda-N}_{(2s + \gamma-1)_+} }\Big)|h||_{L^2}.
\end{equation}
On the other hand,  for $\tilde{J}_{2,j}^{1,2}(\tau)$,  note that 
\begin{equation*}
\tilde{J}_{2,j}^{1,2}(\tau) = \int_{\R^{6}} \left(\int_{\Omega_j^c}
b f(v_*) h(v') \Phi_{*(1)}(v')(v-v')\cdot g_j(z) d \sigma
\right)  dvdv_*,
\end{equation*}
by the symmetry in $\Omega_j^c$. We have
\begin{align*}
|\tilde{J}_{2,j}^{1,2}(\tau)| &\le C \int_{\R^{6}}
\Big(\int_{\Omega_j^c} b(\cos\theta)\theta \la v- v_* \ra^{1-2s}
\la v_* \ra^{(2s+ \gamma-1)_+} |f(v_*)|\\
&\hskip3cm \times
 |h(v')|
\left| \frac{\la z -v_* \ra^{2s+ \gamma -1}}{\la v_*
\ra^{(2s+ \gamma-1)_+}}
 g_j(z) \right|  d\sigma   \Big) dvdv_*
\\
 & \leq
C   2^{-\varepsilon j}\Big (\int_{\R^{6}} \Big(\int_{\Omega_j^c}
b(\cos\theta)\theta 2^{j(1-2s)} \la v- v_* \ra^{1-2s} \la v_*
\ra^{(2s+ \gamma-1)_+} |f(v_*)|
 |h(v')|^2d\sigma   \Big) dvdv_* \Big)^{1/2}\\
 & \hskip1cm \times
\Big (\int_{\R^{6}} \Big(\int_{\Omega_j^c} b(\cos\theta)\theta
2^{j(1-2s)} \la v- v_* \ra^{1-2s}
\la v_* \ra^{(2s+ \gamma-1)_+} |f(v_*)|\\
& \hskip4cm \left| \la z \ra^{(2s+ \gamma-1)_+} 2^{j(2s-1 +
\varepsilon)}  g_j(z) \right|^2
d\sigma   \Big) dvdv_* \Big)^{1/2}\\
&\leq C 2^{-\varepsilon j} ||f||_{L^1_{(2s +\gamma-1)_+}}|h||_{L^2}
 \Big(||M\,g||_{H^{2s-1+ \varepsilon} _{(2s + \gamma-1)_+} }
+||g||_{H^{\lambda-N}_{(2s + \gamma-1)_+} }\Big),
\end{align*}
because of \eqref{second}. This together with \eqref{J-2-10} and
\eqref{J-2-j-1} yield
\begin{equation}\label{estimate-strong}
|{J}_2^1(\tau,(1))| \leq  C ||f||_{L^1_{(2s +\gamma-1)_+}}
 \Big(||M\,g||_{H^{2s-1+ \varepsilon} _{(2s + \gamma-1)_+} }
+||g||_{H^{\lambda-N}_{(2s + \gamma-1)_+} }\Big)|h||_{L^2}.
\end{equation}
It is easy to see that all other terms coming from 
$F_{N_1}(v,D_v;v_*)g(v)$ in \eqref{expansion-1} have the same upper bound
estimates. Moreover,  all the
terms coming from $\tilde{r}_{N_1}(v,D_v;v_*)g(v)$ can be estimated by
\begin{equation*}
C||f||_{L^1_{(2s +\gamma-1)_+}}
||g||_{H^{\lambda-N}_{(2s + \gamma-1)_+} }|h||_{L^2}.
\end{equation*}
Therefore, we finally obtain
\begin{equation*}|J_2| =|J_2^1| \leq C
||f||_{L^1_{(2s +\gamma-1)_+}} \Big(||M\,g||_{H^{2s-1+ \varepsilon} _{(2s + \gamma-1)_+} }
+||g||_{H^{\lambda-N}_{(2s + \gamma-1)_+} }\Big)|h||_{L^2}.
\end{equation*}
In summary, when $1/2< s<1$ we obtain instead of \eqref{CU} that 
\begin{equation}\label{CU2}
|{\mathcal {II}} | \le C ||f||_{L^1_{(2s +\gamma-1)_+}}
\Big(||M\,g||_{H^{2s-1+ \varepsilon} _{(2s + \gamma-1)_+} }
+||g||_{H^{\lambda-N}_{(2s + \gamma-1)_+} }\Big)|h||_{L^2}.
\end{equation}
By combining \eqref{I11}, \eqref{CU} and \eqref{CU2}, the proof of Proposition
\ref{coro-comm-v} is completed.

\vskip0.5cm
\bigskip

The rest of this section is devoted to the proof \eqref{3.11+1*} of Lemma \ref{lemm3.4}.
\\

\noindent
{\bf Proof of \eqref{3.11+1*} of Lemma \ref{lemm3.4}}.For $m = 2s-1+\varepsilon>0$, we have with $\Lam = (1- \Delta_v)^{1/2}$
\begin{align*}
\Big(W_lQ(f,\, g) &- Q(f,\,W_l g),\, h\Big) = \Big((\Lambda^{-m} Q(f,\, g)-
Q(f,\, \Lambda^{-m}g)),\,W_l \Lambda^m  h \Big)\\
& + \Big((W_l Q(f,\, \Lambda^{-m}g)-
Q(f,\, W_l \Lambda^{-m} g)),\, \Lambda^m h \Big)\\
&+ \Big((Q(f,\,  \Lambda^{-m}\, W_l\,g)-
 \Lambda^{-m} Q(f,\, W_l g)),\, \Lambda^m h \Big)\\
 &+ \Big((\big[\Lambda^{-m},\, W_l\big]Q(f,\,g)-
 Q(f,\, \big[\Lambda^{-m},\, W_l\big] g)),\, \Lambda^m h \Big)\\
 &= (1) + (2) + (3) +(4).
 \end{align*}
It follows from \eqref{10.8-3} with $M(\xi)= \Lambda^{-m}$ that
\begin{align*}
&|(1)| \leq C \|f\|_{L^1_{(2s+\gamma-1)^+}}
\|g\|_{L^2_{(2s+\gamma-1)^+}}\|h\|_{H^m_l},\\
&|(3)| \leq C \|f\|_{L^1_{(2s+\gamma-1)^+}}
\|W_l g\|_{L^2_{(2s+\gamma-1)^+}}\|h\|_{H^m}.
\end{align*}
By means of \eqref{3.11+1}, we have
\begin{equation*}
|(2)| \leq C\|f\|_{L^1_{l+2s-1+\gamma^+}}\|g\|_{L^2_{l+2s-1+\gamma^+} }\|h\|_{H^m}.
\end{equation*}
To estimate (4), we first  note that
\begin{align*}
[\Lambda^{-m},\, W_l]&= \sum_{|\alpha|=1} \big( W_l \big)_{(\alpha)}
\big(\Lambda^{-m}\big)^{(\alpha)}
 + W_{l-1}R(v,D_v),
\end{align*}
where $R$ is a pseudo-differential operator
which  belongs to $S^{-m-2}_{1,\,0}$. Write
\begin{align*}
(4) = &\sum_{|\alpha|=1}
\Big( \Big\{ \big(\Lambda^{-m}\big)^{(\alpha)}Q(f,g) -Q(f,\big(\Lambda^{-m}\big)^{(\alpha)}g)
\Big\}
,\, \big( W_l \big)_{(\alpha)} \Lambda^m    h \Big)\\
&+ \sum_{|\alpha|=1}
\Big( \Big\{ \big( W_l \big)_{(\alpha)}Q(f,\big(\Lambda^{-m}\big)^{(\alpha)}g)
-
Q(f, \big(W_l \big)_{(\alpha)}\big(\Lambda^{-m}\big)^{(\alpha)}g)\Big\},
\Lambda^m h \Big)\\
&+\Big(R(v, D_v) Q(f,g), W_{l-1} \Lambda^m h\Big) +
\Big(Q(f, W_{l-1} R(v, D_v) g), \Lambda^m h\Big)\\
=& (a) +(b)+(c)+(d).
\end{align*}
It follows from \eqref{2.2+001} that
\begin{align*}
&|(c)| \leq C\|Q(f,g)\|_{H^{-2}} \|h\|_{H^m_{l-1}} \leq C \|f\|_{L^1_{(\gamma+2s)^+}}
\|g\|_{L^2_{(\gamma+2s)^+}}  \|h\|_{H^m_{l-1}}, \\
&|(d)| \leq C\|Q(f,W_{l-1} R g)\|_{L^2} \|h\|_{H^m}
\leq C
 \|f\|_{L^1_{(\gamma+2s)^+}}
\|g\|_{L^2_{l-1+(\gamma+2s)^+}}  \|h\|_{H^m} .
\end{align*}
By exactly the same method as the one for \eqref{3.11+1}, namely,
by replacing $W_l$ by $(W_l)^{(\alpha)}$ which is bounded by
$W_{l-|\alpha|}$,  we have
\begin{equation*}
|(b)|\leq C \|f\|_{L^1_{l-2+2s+\gamma^+}}\|\Big(\Lambda^{-m}\Big)^{(\alpha)}g\|
_{H^m_{l+2s-2+\gamma^+}}\|h\|_{H^m}\leq C
\|f\|_{L^1_{l-2+2s+\gamma^+}}\|g\|
_{L^2_{l+2s-2+\gamma^+}}\|h\|_{H^m}\,.
\end{equation*}
The estimation on $(a)$ is the same as the argument
in Proposition \ref{coro-comm-v}
by  replacing $M(D)$ by $(\Lambda^{-m})^{(\alpha)}$,
except for the term corresponding to  ${\mathcal I}$.
Notice that $D_{\xi}^\alpha (\left< \xi \right>^{-m}) := M^{(\alpha)}(\xi)$
is no longer a function of $|\xi|^2$.
Instead of \eqref{inegalite}, we only have
\begin{equation}\label{odd-property}
|M^{(\alpha)}(\xi) - M^{(\alpha)}(\xi^+)| \leq C \left| \sin \frac{\theta}{2}\right|
\la \xi^+ \ra^{-m-1}.
\end{equation}
Thus, we need to use the symmetry property as in the proof of Theorem \ref{theo2.1}.
The term corresponding to ${\mathcal I}$ is
\begin{align*}
&{\mathcal I}^\alpha=
\displaystyle\int_{\RR^{6}}\displaystyle\int_{\SS^{2}}
b(\frac{\xi}{|\xi|}\cdot\sigma)f(v_*)
\\
&\hspace*{1.5cm}\times \Big(M^{(\alpha)}(\xi) - M^{(\alpha)}(\xi^+)\Big) {\mathcal
F}(\Phi_* g)(\xi^{+})e^{-iv_{*}\xi^{-}}d\sigma d v_*
\overline{\hat{h_0}(\xi)}d \xi,
\end{align*}
where $h_0 = ( W_l )_{(\alpha)} \Lambda^m  h$.
By letting
\begin{equation*}
F(v,v_*) = \frac{\Phi(|v-v_*|)} { \langle v_*
\rangle^{\gamma^+}}g(v), \enskip \enskip \enskip h(v,v_*)=
\frac{h_0(v)}{\langle v_* \rangle},
\end{equation*}
we write
\begin{align*}
{\mathcal I}^{\alpha} =&\displaystyle\int_{\R^{3}}\langle
v_* \rangle^{1+\gamma^+} f(v_*)\\
&\hspace{-1cm}\times\Big\{\int_{\R^3}
\int_{\SS^{2}}b(\frac{\xi}{|\xi|}\cdot\sigma)
\Big(M^{(\alpha)}(\xi)-M^{(\alpha)}(\xi^+)\Big)
e^{iv_{*}\xi^{+}} \hat F(\xi^+,v_*)
\overline{e^{iv_{*}\xi} \hat h(\xi,v_*)} d \sigma d \xi \Big\} d v_* \\
=& \int_{\R^3} \langle v_* \rangle^{1+\gamma^+} f(v_*) {\mathcal
L}(v_*)dv_*.
\end{align*}
Set
\begin{equation*}
\widetilde{\hat F}(\xi,v_*)= e^{iv_{*}\xi} \hat F(\xi,v_*), \enskip \enskip
\widetilde{ \hat h}(\xi,v_*)=e^{iv_{*}\xi}\hat
h(\xi,v_*),\end{equation*}
 and write
\begin{align*}
{\mathcal L}(v_*) =& \int_{\R^3}
\int_{\SS^{2}}b(\frac{\xi}{|\xi|}\cdot\sigma)
\Big(M^{(\alpha)}(\xi)-M^{(\alpha)}(\xi^+)\Big)
\widetilde{\hat F}(\xi^+,v_*)\overline{\Big( \widetilde{
\hat h}(\xi,v_*)  -\widetilde{ \hat h}(\xi^+,v_*)\Big)} d \sigma d \xi\\
& + \frac{1}{2}\int_0^1(1-\tau)
\Big\{ \int_{\R_n}\int_{\SS^{2}} b\big(\frac{\xi}{|\xi|}\cdot \sigma\big) \\
&\enskip \hskip1cm \times (\nabla^2_{\xi}
M^{(\alpha)})(\xi^+ + \tau(\xi -\xi^+))(\xi^-)^2
\widetilde{\hat F}(\xi^+,v_*)\overline{
\widetilde{ \hat h}(\xi^+,v_*)}d \sigma d \xi \Big\} d \tau \\
=& {\mathcal L}^1(v_*)+{\mathcal L}^2(v_*).
\end{align*}
By the same symmetry property as shown in Figure
\ref{fig2} in the proof of Theorem \ref{theo2.1}, we have
\begin{equation*}
\int_{\R^3}\int_{\SS^{2}} b\big(\frac{\xi}{|\xi|}\cdot \sigma\big)
(\nabla_{\xi} M^{(\alpha)})(\xi^+)\cdot \xi^-(\sigma)
\widetilde{\hat F}(\xi^+,v_*)
\overline{\widetilde{ \hat h}(\xi^+,v_*)}d \sigma
d \xi=0.
\end{equation*}
Then it follows
from \eqref{odd-property}  that
\begin{equation*}
 \sup_{v_*}|{\mathcal L}^1(v_*)| \leq C \|g\|_{L^2_{\gamma^+}}\|h_0\|_{L^2_1}
\leq C\|g\|_{L^2_{\gamma^+}}\|h\|_{H^m_{l}},
\end{equation*}
and
\begin{equation*}
\sup_{v_*}|{\mathcal L}^2(v_*)| \leq C \|g\|_{L^2_{\gamma^+}}\|h_0\|_{L^2}
\leq C \|g\|_{L^2_{\gamma^+}}\|h\|_{H^m_{l-1}},
\end{equation*}
whence we obtain
$$
|{\mathcal I}^{\alpha}| \leq C\|f\|_{L^1_{1+\gamma^+}}
\|g\|_{L^2_{\gamma^+}}\|h\|_{H^m_{l}}.
$$
In summary,
we obtained the desired estimate \eqref{3.11+1*}.

\bigbreak


\section{Regularizing effect}\label{section3}
\smallbreak

In this section, we will prove the regularizing effect on solutions
to the non-cutoff Boltzmann equation starting from $f\in {\mathcal H}^{5}_l (]T_1,
T_2[\,\times\Omega\times\RR^3_{v}))$. 
Actually this will be 
proved by using
an induction argument in the following subsections.
In the first step, we will show the gain of regularity in the variable
$v$ mainly by using the singularity in the cross-section, that is, the
coercivity property in (\ref{3.5}). In the second step, we will apply the 
hypo-elliptic estimate obtained by a generalized version of the uncertainty principle
to show the gain of regularity in $(x,t)$ variables. Then an induction
argument will lead to at least one order higher 
regularity 
in  $(x,t)$ variables. By using the equation and an induction
argument again,  at least one order higher  regularity can
be obtained in  $v$ variable. Therefore, the solution is shown to be in
${\cH}^6_l (]T_1,\, T_2[\,
\times\Omega\times\RR^3_{v})$ which by induction leads to
${\cH}^\infty_l (]T_1,\, T_2[\,
\times\Omega\times\RR^3_{v})$.

Let $f\in {\mathcal H}^{5}_l (]T_1,
T_2[\,\times\Omega\times\RR^3_{v}))$,  for all $l \in\NN$, be a
(classical)  solution of the Boltzmann  equation (\ref{1.1}). 
We now want  to prove the full regularity of $ \varphi(t)\psi(x)f$ for
any smooth cutoff functions $\varphi\in C^\infty_0(]T_1, T_2[),\,
\psi\in C^\infty_0(\Omega)$.

\subsection{Initialization}\label{s3}
\setcounter{equation}{0} \smallbreak

Here and below,  $\phi$ denotes a cutoff  function satisfying $\phi\in
C^\infty_0$ and $0\leq \phi\leq 1$. Notation
$\phi_1\subset\subset\phi_2$ stands for two cutoff functions such
that $\phi_2=1$ on the support of $\phi_1$.

Take the smooth cutoff
functions $\varphi,\, \varphi_2, \varphi_3\in C^\infty_0(]T_1, T_2[)$ and $
\psi,\, \psi_2, \psi_3\in C^\infty_0(\Omega)$  such that $\varphi\subset\subset \varphi_2\subset\subset \varphi_3$ and
$\psi\subset\subset\psi_2\subset\subset\psi_3$. Set $f_1=\varphi(t)\psi(x) f$,
$f_2=\varphi_2(t) \psi_2(x) f$ and $f_3=\varphi_3(t) \psi_3(x) f$.
For $\alpha\in\NN^7, |\alpha|\leq 5$, define
$$
g=\partial^\alpha (\varphi(t)\psi(x) f)=\partial^\alpha_{t, x, v}
(\varphi(t)\psi(x) f)\in L^2_l
(\RR^7).
$$
Firstly, the translation invariance of the collision operator with respect to
the 
variable $v$ implies that (see \cite{desv-wen1,Grad,ukai} ), for the
translation operation $\tau_h$
in $v$ by $h$, that we have
$$
\tau_h G(f,\, g)=Q(\tau_h f, \tau_h g).
$$
Then the Leibniz formula with respect to the $t, x$ variables yields
 the following equation  in
a weak sense
\begin{equation}\label{3.1}
 g_t + v\,\cdot\,\partial_x {g } = Q(f_2,\,\, g)+ G,
\enskip (t, x, v) \in \RR^7,
\end{equation}
where
\begin{eqnarray}\label{3.1**}
G&=&\sum_{\alpha_1+\alpha_2=\alpha,\,\, 1\leq |\alpha_1|}C^{\alpha_1}_{\alpha_2}
Q\Big(\partial^{\alpha_1} {f_2},\,\,
\partial^{\alpha_2} f_1\Big)\\
&+&\partial^\alpha\Big( \varphi_t \psi(x) f+
v\,\cdot\,\psi_x(x) \varphi(t) f\Big)
+[\partial^\alpha,\,\,\,
v\,\cdot\,\partial_x](\varphi(t)\psi(x) f)\nonumber \\
&\equiv & (A)+(B)+(C).\nonumber
\end{eqnarray}

To prove the regularity of $g=\partial^\alpha (\varphi(t)\psi(x)
f)$, the natural idea would be to use $g$ as a test function for
equation (\ref{3.1}). But at this point, $g$ only  belongs to
$L^2_l (\RR^7)$ so that  it is only a weak solution to equation
(\ref{3.1}). By using the upper bound estimate on $Q$, we have
$Q(f_2, g)\in
L^2(\RR^4_{t, x};\, H^{-2s}(\RR^3_v))$. Thus, we need to choose the
test functions at least in the space
$L^2(\RR^4_{t, x};\, H^{2s}(\RR^3_v))$.
For this, we will  use a mollification of $g$ with respect to 
the variables $(x,\, v)$ as a test function.

\smallbreak
For this purpose, let $S\in C^\infty_0(\RR)$ satisfy $0\leq
S\leq 1$ and
$$
S(\tau)=1,\,\,\,\,|\tau|\leq 1;\,\,\,\,\,\,\,
S(\tau)=0,\,\,\,\,|\tau|\geq 2.
$$
Then
$$
S_N(D_x)S_N(D_v)=S(2^{-2N}|D_x|^2)S(2^{-2N}|D_v|^2)\, :\,\,
H^{-\infty}_l (\RR^6) \,\, \rightarrow\,\, H^{+\infty}_l  (\RR^6),
$$
is a regularization operator such that
$$
\|\big(S_N(D_x)S_N(D_v)f\big) -f\|_{L^2_l  (\RR^6)}\rightarrow 0,
\,\,\,\,\,\, \mbox{ as } N\rightarrow \infty.
$$
Choose another cutoff function $\psi\subset\subset\psi_1
\subset\subset\psi_2$ and set
$$
P_{N,\, l}=\psi_1(x) S_N(D_x)\, W_{l  }\,\, S_N(D_v).
$$
Then we can take
$$
\tilde{g}=P^\star_{N,\, l}\,\,( P_{N,\, l}\,\,g) \in C^1(\RR;
H^{+\infty}(\RR^6))
$$
as a test function for the equation (\ref{3.1}).


It follows by integration by parts on $\RR^7=\RR^1_t\times
\RR^3_x\times\RR^3_v$ that
\begin{eqnarray*}
&&\Big([S_N(D_v),\,\, v]\,\cdot\,\nabla_x S_N(D_x) g,\,\,
\psi_1(x)W_{l  } P_{N, l  }\,\,g\Big)_{L^2(\RR^7)}=
\\&&\,\,\,\,\,\,\,\, \Big(P_{N,
l  }\,Q({f_2},  g),\, P_{N, l  }\,\, g\Big)_{L^2(\RR^7)}+ \Big(G,\,
\tilde{g} \Big)_{L^2(\RR^7)},
\end{eqnarray*}
which implies that 
\begin{eqnarray}\label{3.5}
&&-\Big(Q({f_2}, P_{N,\, l}\, g),\,\, P_{N,\, l}\,
g\Big)_{L^2(\RR^7)}= -\Big([S_N(D_v),\,\, v]\,\cdot\,\nabla_x
S_N(D_x) g,\,\, \psi_1(x)W_{l  } P_{N, l  }\,\,g\Big)_{L^2(\RR^7)}
\\
&&\,\,\,\,\,\,\,\,\,\,\,+\Big(P_{N,\,l}\,Q(f_2,  g)-Q({f_2}, P_{N,\,
l}\, g),\, P_{N,\, l}\,\,g\Big)_{L^2(\RR^7)}+ \Big(G,\, \tilde{g}
\Big)_{L^2(\RR^7)}.\nonumber
\end{eqnarray}

By using (\ref{3.5}), we can deduce the regularity of $g$ from the coercivity 
property of the collision operator on
the left hand side and the upper bound estimate on the right hand side. The detailed argument will be given in the next subsection.

\subsection{Gain of regularity in $v$ }\label{s3++}
\setcounter{equation}{0} \smallbreak

In this subsection, we will prove
a  partial smoothing effect of the cross-section
 on the weak solution $g$ in
the velocity variable $v$ .

\begin{prop}\label{prop3.1}
Assume that $0<s<1,\,\, \gamma\in\RR$. Let $f\in {\cH}^5_l (]T_1,\, T_2[\,
\times\Omega\times\RR^3_{v})$ be a solution of the equation (\ref{1.1})
for all $l \in\NN$.
Assume furthermore that
\begin{equation}\label{1.4++}
f(t, x, v)\geq 0\,\,\,\,\mbox{and}\,\,\,\,\, \|f(t, x,
\cdot)\|_{L^1(\RR^3_v)}>0,
\end{equation}
for all $(t, x, v)\in ]T_1,\, T_2[\,\times\Omega\times\RR^3_{v}$. Then one has,
\begin{equation}\label{3.3}
\Lambda^s_v f_1\in H^5_l  (\RR^7),
\end{equation}
for any $l\in\NN$, where $f_1=\varphi(t)\psi(x)
f$ with $\varphi\in C^\infty_0(]T_1, T_2[), \psi\in C^\infty(\Omega)$.
\end{prop}

\noindent{Proof :} Firstly,  the local positive lower bound  assumption (\ref{1.4++}) implies that 
 $$
 \inf_{(t,  x)\in \Supp\,\varphi \times \Supp\,\psi_1}
 \|{f_2}(t, x, \cdot)\|_{L^1(\RR^3_v)}=c_0>0.
 $$
Thus,  the coercivity estimate (\ref{E-sub-estimate+}) in
Theorem \ref{E-theo-0.1} gives that for any  $\gamma \in\R,\,
0<s<1$,
\begin{eqnarray*}\label{3.7}
&&-\Big(Q({f_2},\,  P_{N, l  }\, g),\,\,
 P_{N, l  }\, g\Big)_{L^2(\RR^7)}=-\int_{t\in
 \Supp \varphi }\int_{x\in
 \Supp \psi_1 }\Big(Q({f_2},\,  P_{N, l  }\, g),\,\,
 P_{N, l  }\, g\Big)_{L^2(\RR^3_v)} dx dt
\\
&& \geq \int_{\RR_t}\int_{\RR^3_x}\Big( C_0\| W_{\gamma/2}
 P_{N, l  }\, g(t, x, \cdot)\|^2_{H^s(\RR^3_v)}\\
 &&
 \hskip4cm-C \|{f_2}(t, x, \cdot)\|_{L^1_{\max\{\gamma^+,\, 2-\gamma^+\}}(\RR^3_v)}
 \| P_{N, l  }\, g(t, x, \cdot)\|^2_{L^2_{\gamma^+/2}(\RR^3_v)}\Big)dx dt \nonumber\\
&&\geq C_0\|\Lambda^{s}_v W_{\gamma/2} P_{N, l }\,
g\|^2_{L^2(\RR^7)}- C \| {f_2}\|_{L^\infty(\RR^4_{t, x};\,\,
L^1_{\max\{\gamma^+,\, 2-\gamma^+\}}(\RR^3_v))}\|W_{l
}\,\,g\|^2_{L^2_{\gamma^+/2}(\RR^7)}, \nonumber
\end{eqnarray*}
where $C_0$  depends on $c_0, \sup_{t, x}\|f_2(t, x,
\cdot)\|_{L^1_1(\RR^3_v)}$ and $\, \sup_{t, x} \|f_2(t, x,
\cdot)\|_{L\log L(\RR^3_v)}$, see Remark \ref{coercivity}. 

\smallbreak

For the terms in (\ref{3.5}), note first of all that
\begin{equation}\label{bounded}
[S_N(D_v),\,\, v]\,\cdot\,\nabla_x\,\, S_N(D_x) =2^{-2N}\big(
S'\big)_N(D_v)\, D_v\,\cdot\,\nabla_x\, \,S_N(D_x)\, :\,
L^2(\RR^6_{x, v})\,\,\rightarrow\,\,
L^2(\RR^6_{x, v}),
\end{equation}
is a uniformly bounded operator so that 
\begin{equation*}\label{3.5+2}
\Big|\Big([S_N(D_v),\,\, v]\,\cdot\, \nabla_x\,\, S_N(D_x) g,\,\,
\psi_1(x)W_{l  } P_{N, l  }\,\,g\Big)_{L^2(\RR^7)}\Big|\leq C \|f_1\|^{2}_{H^5_{
l}(\RR^7)}.
\end{equation*}
Hence, by 
using (\ref{3.5}), we get, for $l>3/2+2$,
\begin{eqnarray}\label{3.6}
\|\Lambda^{s}_v W_{\gamma/2} P_{N, l  }\, g\|^2_{L^2(\RR^7)}&\leq& C
\Big\{\Big(1+\|f_2\|_{H^{2+\delta}_{l+\gamma^+}(\RR^6)}
\Big)\|f_1\|^{2}_{H^5_{l}(\RR^7)}
+\left|\Big(G,\, \tilde{g}\Big)_{L^2(\RR^7)}\right|\\
&&+\Big|\Big(P_{N, l }\, Q({f_2}, g)-Q({f_2},\, P_{N, l  }\, g), \,\,
P_{N, l }\,g\Big)_{L^2(\RR^7)}\Big| \Big\}\, .\nonumber
\end{eqnarray}
The above constants $C>0$ are independent of  $N$.

We complete the proof of Proposition \ref{prop3.1}
by estimating the  last two terms in \eqref{3.6} through the 
following three Lemmas.

\begin{lemm}\label{prop3.1+0}
Assume  $0<s<1, \gamma\in\RR$.  Let $f\in {\cH}^5_l (]T_1,\,
T_2[\, \times\Omega\times\RR^3_{v}),\, l \geq 3/2+2$. Then, for any
$\alpha\in\NN^7, |\alpha|\leq 5$, we have, for any $\varepsilon>0$,
\begin{equation}\label{3.5+1}
\left|\Big(G,\, \tilde{g}\Big)_{L^2(\RR^7)}\right|\leq
C_\varepsilon\|f_3\|^{4}_{H^{5}_{ l+4+|\gamma|}(\RR^7)}+\varepsilon
\|\Lambda^s_v W_{\gamma/2} P_{N,\, l}\, g\|^2_{L^2(\RR^7_{t, x,
v})}.
\end{equation}
\end{lemm}

\smallbreak \noindent {\bf Proof :} Firstly, we prove that 
\begin{equation}\label{3.2}
G\in L^2(\RR^4_{t, x};  H^{-(2s-1+\delta)^+}_{l}(\RR^3_v)),
\end{equation}
for any $l\in\NN$, where $(2s-1+\delta)^+=\max\{2s-1+\delta, 0\}$
and $\delta>0$ satisfying $2s-1+\delta<s$. 
By using the decomposition in \eqref{3.1**},
it is obvious that
$$
(B)=\partial^\alpha\Big( \varphi_t \psi(x) f+
v\,\cdot\,\psi_x(x) \varphi(t) f\Big)\in L^2_l (\RR^7),
$$
and
$$
\|(B)\|_{L^2_l(\RR^7)}\leq C \|f_2\|_{H^5_{l+1}(\RR^7)}.
$$
Since $[\partial^\alpha,\,\,\, v\,\cdot\,\partial_x]$ is a
differential operator of order $|\alpha|$, we have
$$
\|(C)\|_{L^2_l(\RR^7)}\leq C \|f_2\|_{H^{5}_{l}(\RR^7))}.
$$

For the term $(A)$, recall that $\alpha_1+\alpha_2=\alpha$,
$|\alpha|\leq 5$ and $|\alpha_2|<{5}$. In the following,
 we will apply Theorem
\ref{theo2.1} with $m=1-\delta-2s$. We separate the discussion into
two cases.\\

\noindent{\bf Case 1.} If $|\alpha_1|=1, 2 $, we have
\begin{eqnarray*}
&&\int_{\RR_t}\int_{\RR^3_x}\|Q(\partial^{\alpha_1}
{f_2},\,\,\partial^{\alpha_2} f_1)(t, x,
\cdot)\|^2_{H^{1-\delta-2s}_l(\RR^3_v)}dx dt\\
&\leq&
C\int_{\RR_t}\int_{\RR^3_x}\|\partial^{\alpha_1}
{f_2}(t, x, \cdot)\|^2_{L^{1}_{l+(2s+\gamma)^+}(\RR^3_v)} \|
\partial^{\alpha_2} f_1(t, x,
\cdot)\|^2_{H^{1-\delta}_{l+(2s+\gamma)^+}(\RR^3_v)}dx dt\\
&\leq& C\|\partial^{\alpha_1} {f_2}\|^2_{L^\infty(\RR^4_{t, x};\,
L^{1}_{l+(2s+\gamma)^+}(\RR^3_v))}
\int_{\RR_t}\int_{\RR^3_x} \|
\partial^{\alpha_2} f_1(t, x,
\cdot)\|^2_{H^{1}_{l+(2s+ \gamma)^+}(\RR^3_v)}dx dt\\
&\leq& C\|f_2\|^2_{H^{2+4/2+\delta}_{l+3/2+\delta+(2s+ \gamma)^+}(\RR^7)}
\|f_1\|^2_{H^{5}_{l+(2s+ \gamma)^+}(\RR^7)}.
\end{eqnarray*}
\noindent{\bf Case 2.} If $|\alpha_1|\geq 3$, then $|\alpha_2|\leq 2$, it follows that
\begin{eqnarray*}
&&\int_{\RR_t}\int_{\RR^3_x}\|\partial^{\alpha_1}
{f_2}(t, x, \cdot)\|^2_{L^{1}_{l+(2s+\gamma)^+}(\RR^3_v)} \|
\partial^{\alpha_2} f_1(t, x,
\cdot)\|^2_{H^{1-\delta}_{l+(2s+\gamma)^+}(\RR^3_v)}dx dt\\
&\leq& C \|\partial^{\alpha_2} {f_1}\|^2_{L^\infty(\RR^4_{t,x}; \,
H^{1-\delta}_{l+(2s+\gamma)^+}(\RR^3_v))} \int_{\RR_t}\int_{\RR^3_x}
\|
\partial^{\alpha_1} f_2(t, x,
\cdot)\|^2_{L^{2}_{l+3/2+\delta+(2s+\gamma)^+}(\RR^3_v)}dx dt\\
&\leq& C\|f_1\|^2_{H^{2+1-\delta+4/2+\delta/2}_{l+(2s+ \gamma)^+}(\RR^6)}
\|f_2\|^2_{H^{5}_{l+3/2+\delta+(2s+ \gamma)^+}(\RR^7)}.
\end{eqnarray*}
By combining these two cases,  we have proved (\ref{3.2}).

Now if $2s-1<0$, then (\ref{3.2}) implies that
\[
\left|\Big(G,\, \tilde{g}\Big)_{L^2(\RR^7)}\right|\leq C\|f_3\|^3_{H^5_{l+4+\gamma^+}(\RR^7)}.
\]
On the other hand, if $0\leq 2s-1$ and $\gamma<0$ (the case
$\gamma>0$ is easier), then (\ref{3.2}) implies that
\begin{eqnarray*}
&&\left|\Big(G,\, \tilde{g}\Big)_{L^2(\RR^7)}\right|\leq
\|G\|_{L^2(\RR^4_{t, x};
H^{1-2s-\delta}_{l+|\gamma|/2}(\RR^3_v)}\|W_{-|\gamma|/2}P_{N,\, l}\,
g\|_{L^2(\RR^4_{t, x};\, H^{2s-1+\delta}(\RR^3_v))}\\
&& \leq C\|f_3\|^2_{H^5_{l+4+|\gamma|}(\RR^7)}\,\, \|W_{-|\gamma|/2}P_{N,\, l}\,
g\|_{L^2(\RR^4_{t, x};\,
H^{2s-1+\delta}(\RR^3_v))},
\end{eqnarray*}
because $2s-1+\delta<s$. Therefore, the proof
of Lemma \ref{prop3.1+0} is completed.\\

We now  turn to  the estimates of commutators between the mollification operators and  the collision
operator, which are given in the  following two lemmas.

\begin{lemm}\label{lemm3.1}
For any  $\gamma\in\RR$, we have \\
{\rm (1)} If $0<s <1/2$,
then for any suitable functions $f$ and $g$ with the
following norms well defined, one has
\begin{equation}\label{3.8}
\|S_N(D_v) Q(f, g)-Q(f,\, S_N(D_v) g)\|_{L^2(\RR^3_v)} \leq C
\|f\|_{L^1_{\gamma^+}(\RR^3_v)} \|g\|_{L^2_{\gamma^+} (\RR^3_v)},
\end{equation}
for some constant $C$ independent of $N$.\\
{\rm (2)} If $1/2 < s <1$, then for any $\delta >0$ there exists
a constant $C_\delta >0$ such that
\begin{equation}\label{3.8-1}
\|S_N(D_v) Q(f, g)-Q(f,\, S_N(D_v) g)\|_{L^2(\RR^3_v)} \leq
C_\delta \|f\|_{L^1_{(2s+ \gamma-1)^+}(\RR^3_v)}
\|g\|_{H^{2s-1+\delta} _{(2s+ \gamma-1)^+}(\RR^3_v)}\, ,
\end{equation}
and
\begin{equation}\label{3.8-1*}
\|S_N(D_v) Q(f, g)-Q(f,\, S_N(D_v) g)\|_{H^{1-2s-\delta}
(\RR^3)} \leq C_\delta \|f\|_{L^1_{(2s+ \gamma-1)^+}(\RR^3_v)}
\|g\|_{L^2_{(2s+ \gamma-1)^+}(\RR^3_v)}.
\end{equation}
{\rm (3)} When $s = 1/2$, we have the same form of estimate as (\ref{3.8-1}) with $(2s+
\gamma-1)$ replaced by $(\gamma+ \kappa)$ for any small $\kappa >0$.
\vskip 0.5cm
\end{lemm}

Before giving the proof of this lemma, notice
that when $\gamma=0$ in the  Maxwellian molecule case, the following 
proof of Lemma \ref{lemm3.1} is similar to Lemma 3.1 in
\cite{MUXY-DCDS} (see also Lemma 5.1 in \cite{amuxy-nonlinear}) by
using the Fourier transformation of collision operator. However,  here we
 consider the case for  $\gamma\in\RR$.\\

\smallbreak\noindent {\bf Proof of Lemma \ref{lemm3.1} :} The proof
is a slight modification of  the proof for  Proposition
\ref{coro-comm-v}. Set
$$
M(|\xi|)=S_N(|\xi|)=S(2^{-2N}|\xi|^2).
$$
Then $S_N\in S^{0}_ {1,0}$ uniformly. 
Even though it does not satisfy (\ref{keep-form}),  we have
$$
|\partial^\alpha S_N(|\xi|)|\leq C_\alpha S_{N+1}(|\xi|)<\xi>^{-|\alpha|}
$$
with $C_\alpha$ independent of $N\in\NN$. Thus, (\ref{10.8-2}) implies
(\ref{3.8}) and (\ref{10.8-3}) implies (\ref{3.8-1}) respectively.

For (\ref{3.8-1*}), note that with $m= 2s-1+\delta$ we have
\begin{align*}
(S_N Q(f,g) - Q(f,S_N g), h) =& \left((\Lambda^{-m} Q(f,g)-
Q(f, \Lambda^{-m}g)), \Lambda^m S_N h \right)\\
& + \left((S_N Q(f, \Lambda^{-m}g)-
Q(f, \Lambda^{-m}S_N g)), \Lambda^m h \right)\\
&+ \left((Q(f, S_N \Lambda^{-m}g)-
 \Lambda^{-m} Q(f, S_N g)), \Lambda^m h \right)\\
 =& (I_1) + (I_2) + (I_3).
 \end{align*}
By applying \eqref{10.8-3} with $M (\xi) = \langle\xi \rangle^{-m}$
to $(I_1)$ and $(I_3)$, we obtain
\begin{equation*}
|(I_1)| +|(I_3)| \leq C \|f\|_{L^1_{(2s+\gamma-1)^+}}
\|g\|_{L^2_{(2s+\gamma-1)^+}}\|h\|_{H^m},
\end{equation*}
because $S_N \in S^0_{1,0}$ uniformly. The same bound on $(I_2$) follows
from \eqref{3.8-1}.

Notice  that the case of $s=1/2$ follows from the case of $s = 1/2 + \kappa$ 
for any positive $\kappa$
because the main concern here is the upper bound. And this completes the proof of the lemma.\\

The following lemma is on the
commutator of the collision opertor with mollifier in the $x$ variable.

\begin{lemm}\label{lemm3.3}
Let $0< s <1$ and $\gamma,  m\in\RR$. For any suitable functions $f$ and $h$
with the following norms well defined, one has
\begin{eqnarray}\label{3.10+1}
&&\|S_N(D_x) Q(f,\,\, h )-Q(f,\,\, S_N(D_x)\, h)\|_{L^2(\RR^4_{t, x},\, H^{m-2s}(\RR^3_v))} \\
&&\,\,\,\,\,\,\,\,\,\,\leq C 2^{-N} \|\nabla_x f\|_{L^\infty(\RR^4_{t,
x},\,\, L^1_{(2s+\gamma)^+}(\RR^3_v))}\|h\|_{L^2(\RR^4_{t, x},\,\,
H^m_{(2s+\gamma)^+}(\RR^3_v))}.\nonumber
\end{eqnarray}
for a constant $C$ independent of $N$.
\end{lemm}

\noindent {\bf Proof :} Let us introduce $\tilde{K}_N(z)=
2^{3N}\hat{S}(2^N z)2^N z$. Note that $\tilde{K}_N \in L^1(\RR^3)$
uniformly with respect to $N$. Then for any smooth function $\tilde{h}$, one
has
\begin{eqnarray*}
&&\Big(\big(S_N(D_x) \,\, Q(f,\,\,h)-Q(f,\,\,S_N(D_x)
h)\big),\,\, \tilde{h}\Big)_{L^2(\RR^7)}
= \int_0^1\Big\{ \int_{\RR_t} \int_{\RR_x^3 \times \RR_y^{3}}
\tilde{K}_N(x-y)\\
&&\,\,\,\,\,\,\,\,\,\,\,\times \Big(Q\big(\nabla_x f(t,y+\tau(x-y),\, \cdot\,),\,\,
2^{-N}h(t, y, \, \cdot\,)\big),\,\, \tilde{h}(t,x, \,
\cdot\,)\Big)_{L^2(\RR_v^3)}dtdxdy\Big\} d\tau .
\end{eqnarray*}
By applying Theorem \ref{theo2.1} with $m-2s$, the right hand side of
this equality can be estimated from above by
$$
C\Big\{\sup_{t,x} ||\nabla_x f(t, x,
\,\cdot\,)||_{L^1_{(2s+\gamma)^+}(\RR^3_v)}\Big \}\times \hskip7cm $$
$$ \int_{\RR_t}\int_{
\RR_x^3}\big(|\tilde{K}_N|*||2^{-N} h (t,\,\cdot\,)
||_{H^m_{(2s+\gamma)_+}(\RR^3_v)}\big)(x)||\tilde{h}(t, x, \cdot)||_{H^{2s-m}(\RR^3_v)}dxdt$$
$$\hskip1cm
\leq C 2^{-N}\|\nabla_x f\|_{L^\infty(\RR^4_{t, x}; L^1_{(2s+\gamma)^+}(\RR^3_v))}
\|h\|_{L^2(\RR_{t,x}^4;\, H^m_{(2s+\gamma)^+}(\RR_v^3))}||\tilde{h}||_{L^2(\RR_{t,x}^4;
\, H^{2s-m}(\RR_v^3))},
$$
which completes the proof of the lemma.\\

We now apply \eqref{3.10+1} with $h=S_N(D_v)g$ and $m=1$, we get 
\begin{eqnarray}\label{3.10+2}
&&\|S_N(D_x) Q(f,\,\, S_N(D_v)g)-Q(f,\,\, S_N(D_x)S_N(D_v)
g)\|_{L^2(\RR^4_{t, x},\, H^{1-2s}(\RR^3_v))} \\
&&\,\,\,\,\,\,\,\,\,\,\leq C \|\nabla_x f\|_{L^\infty(\RR^4_{t,
x},\,\, L^1_{(2s+\gamma)^+}(\RR^3_v))}\|g\|_{L^2(\RR^4_{t, x},\,\,
L^2_{(2s+\gamma)^+}(\RR^3_v))}.\nonumber
\end{eqnarray}
Here, we have used  the fact that
a mollification operator $S_N(D_v)$ in the $v$ variable
has the property that 
$$
\|2^{-N}S_N(D_v)g(t, x, \,\cdot\,) \|_{H^1_{(2s+\gamma)^+}(\RR^3_v)}\leq C
\|g(t, x, \,\cdot\,)\|_{L^2_{(2s+\gamma)^+}(\RR_v^3)},
$$
where $C$ is a constant independent on $N$.\\

Now we are ready to complete the proof of Proposition \ref{prop3.1}.\\

\smallbreak\noindent
{\bf Completion of proof of Proposition \ref{prop3.1}.}

\smallbreak

We study now the commutator terms in (\ref{3.6}). For this purpose, note that
\begin{eqnarray}\label{3.6++1}
&&\Big(P_{N, l }\, Q({f_2}, g)-Q({f_2},\, P_{N, l }\, g), \,\, P_{N,
l}\,g\Big)_{L^2(\RR^7)}\\
&=&\Big(S_{N}(D_v)\,Q({f_2}, g)-Q({f_2},\, S_{N}(D_v)\, g), \,\,
S^\star_N(D_x)\psi_1(x) W_l P_{N,
l }\,g\Big)_{L^2(\RR^7)}\nonumber\\
&+&\Big(S_{N} (D_x)\, Q({f_2}, S_{N}(D_v)\,g)-Q({f_2},\, S_{N}
(D_x)S_{N}(D_v)\, g), \,\,\psi_1(x) W_l  P_{N, l
}\,g\Big)_{L^2(\RR^7)}\nonumber\\
&+&\Big(\psi_1(x) W_l \, Q({f_2},  S_{N}
(D_x)S_{N}(D_v)\,g)-Q({f_2},\, P_{N, l }\, g), \,\, P_{N, l
}\,g\Big)_{L^2(\RR^7)}\,.\nonumber
\\
&=&(1)+(2)+(3).\nonumber
\end{eqnarray}
Note that $\Lambda^{s}_v [\psi_1(x),\,\,S_N(D_x)]S_N(D_v)$ is an
$L^2$ uniformly bounded operator with respect to the parameter $N$ 
for $0\leq s\leq 1$, and that $[W_l,\,\,S_N(D_v)]$ is also a uniformly
bounded operator from $L^2$ to $L^2_{l-1}$ with respect to the
parameter $N$. The discussion on (\ref{3.6++1}) can be divided into
the following two cases.

\smallbreak\noindent {\bf Case 1. $0 < s < 1/2$. }  In this
case, Lemma
\ref{lemm3.1} implies that , for $l>\max\{4, (\gamma+2s)^+\}$,
$$
|(1)|\leq C \|f_2\|_{L^\infty(\RR^4_{t, x}\, ,\,\, L^1_{\gamma^+
+2s}(\RR^3_v))}\|g\|_{L^2 (\RR^7)} \|g\|_{L^2_{2l}(\RR^7)}\leq
C\|f_3\|^3_{H^5_{2l}(\RR^7)}.
$$
And Lemma \ref{lemm3.3} implies that
$$
|(2)|\leq C \|\nabla_x f_2\|_{L^\infty(\RR^4_{t, x}\, ,\,\,
L^1_{\gamma^+ +2s}(\RR^3_v))}\|g\|_{L^2(\RR^4_{t, x}\, ,\,\, L^2_{
{\gamma^+}}(\RR^3_v))}\|g\|_{L^2_{2l}(\RR^7)}\leq
C\|f_3\|^3_{H^5_{2l}(\RR^7)}.
$$
As for the term (3), we use Lemma \ref{lemm3.4} to have 
$$
|(3)| \leq C \|f_2\|_{L^\infty(\RR^4_{t, x}\, ,\,\, L^1_{l+\gamma^+
+2s}(\RR^3_v))}\|g\|_{L^2_{l+\gamma^+ +2s}(\RR^7)} \|P_{N,\, l}\,
g\|_{L^2(\RR^7)}\leq C\|f_3\|^3_{H^5_{2l}(\RR^7)}.
$$

\smallbreak\noindent {\bf Case 2. $1/2\leq s<1$.}  By using
(\ref{3.8-1*}), we have
$$
|(1)|\leq C \|f_2\|_{L^\infty(\RR^4_{t, x}\, ,\,\, L^1_{l+\gamma^+
+2s-1}(\RR^3_v))}\|g\|_{L^2_{l+\gamma^+ +2s-1} (\RR^7)} \|W_{\gamma/2}P_{N,\,
l}\,g\|_{L^2(\RR^4_{t, x}\, ,\,\,
H^{2s-1+\delta}(\RR^3_v))}
$$
$$
\leq \varepsilon\|\Lambda^{s}_v W_{\gamma/2} P_{N, l }\,
g\|^2_{L^2(\RR^7)}+ C_\varepsilon\|f_3\|^{4}_{H^5_{l+4+\gamma^+}(\RR^7)}.
$$
We can use (\ref{3.10+2}) to show that
\begin{eqnarray*}\label{3.9+3}
|(2)|&\leq& C \|\nabla_x f_2\|_{L^\infty(\RR^4_{t, x}\, ,\,\,
L^1_{\gamma^+ +2s}(\RR^3_v))}\|g\|_{L^2(\RR^4_{t, x}\, ,\,\,
L^2_{\gamma^+ +2s}(\RR^3_v))}\|W_{l  }\, P_{N,\, l  }\,g\|_{L
^2(\RR^4_{t,
x}\, ,\,\, H^{2s-1}(\RR^3_v))}\\
&\leq& C_\varepsilon\|f_3\|^{\frac{4}{\theta}+2}_{H^5_{k l}(\RR^7)}
+\varepsilon \|\Lambda^s_v W_{\gamma/2} P_{N,\,
l}\, g\|^2_{L^2(\RR^7_{t, x, v})}.\nonumber
\end{eqnarray*}
Then, (\ref{3.11+1}) implies that
\begin{eqnarray*}\label{3.9+5}
|(3)| &\leq& C \|f_2\|_{L^\infty(\RR^4_{t, x}\, ,\,\,
L^1_{2s+l-1+\gamma_+}(\RR^3_v))}\|\psi_1(x)S_{N}(D_x)\,
S_{N}(D_v)\,g\|_{L^2(\RR_{t,x}^4\, ,\,
H^{2s-1+\delta}_{2s+l-1+\gamma_+}(\RR_v^3))} \|
P_{N,\, l }\,g\|_{L^2(\RR^7)}\\
&\leq& C_\varepsilon\|f_3\|^{\frac{4}{\theta}+2}_{H^5_{k l}(\RR^7)}
+\varepsilon \|\Lambda^s_v W_{\gamma/2} P_{N,\,
l}\, g\|^2_{L^2(\RR^7_{t, x, v})}.\nonumber
\end{eqnarray*}
In summary, we have obtained the following estimate for the second
term on the right hand side of (\ref{3.6})
\begin{eqnarray*}
&&\Big|\Big(P_{N, l  }\, Q({f}_2, g)-Q({f}_2,\, P_{N, l  }\,
g), \,\, P_{N, l  }\,g\Big)_{L^2(\RR^7)}\Big|\\
&&\leq C_\varepsilon\|f_3\|^{2\, k'}_{H^5_{k
l}(\RR^7)}+\varepsilon \|\Lambda^s_v W_{\gamma/2} P_{N,\, l}\,
g\|^2_{L^2(\RR^7_{t, x, v})}.
\end{eqnarray*}
Finally, it holds that 
\begin{equation}\label{3.12}
\|\Lambda^{s}_v W_{\gamma/2} P_{N, l  }\, g\|^2_{L^2(\RR^7)}\leq C
\|f_3\|^{2\, k'}_{H^5_{k l}(\RR^7)},
\end{equation}
where the constants $C$, $k$, and $k'$ are independent of  $N$.
Therefore, Proposition \ref{prop3.1} is proved by taking the limit
$N\rightarrow\infty$.

\vskip0.5cm
\subsection{Gain of regularity in $(t, x)$}
\label{s4+001} \setcounter{equation}{0}
\smallskip

First of all, let us consider a transport equation in the form of
\begin{equation}\label{2.1}
 f_t + v\cdot\nabla_x f = g \in D '({\RR}^{2n+1}) ,
\end{equation}
where $(t,x,v) \in \RR^{1 + n +n} = \RR^{2n +1}$. In
\cite{amuxy-nonlinear-b}, by using a generalized uncertainty
principle, we proved the following hypo-elliptic estimate.

\begin{lemm}\label{lemm2.1}
Assume that $ g \in H^{-s'} (\RR^{2n+1})$, for some $0\leq s' <1$.
Let $f\in L^2 (\RR^{2n+1}) $ be a weak solution of the transport
equation (\ref{2.1}) such that $\Lambda^s_v\, f \in L^2
(\RR^{2n+1})$ for some $0<s\leq 1$. Then it follows that
\[
\Lambda_x^{s (1-s')/(s +1)}f\in L^2_{-\frac{s s'}{s +1}}(\RR^{2n+1})
, \enskip \enskip \Lambda_t^{s (1-s')/(s +1)}f\in L^2_{-\frac{s}{s
+1}}(\RR^{2n+1}),
\]
where $\Lambda_\bullet=(1+|D_\bullet|^2)^{1/2}$.
\end{lemm}

As mentioned earlier, this hypo-elliptic estimate together
with Proposition \ref{prop3.1} are used to obtain
the partial regularity in the variable $(t,x)$. 
With this partial regularity in $(t,x)$, by
applying a Leibniz type formula for fractional derivatives, 
we will show some improved regularity in 
all variables, $v$ and $(t,x)$. Then the hypo-elliptic
estimate can be used again to get higher regularity in 
the variable $(t,x)$. This procedure can be continued to
obtain at least one order higher 
regularity
in $(t,x)$ variable.

For the details, we first recall a Leibniz type formula for fractional
derivatives with respect to variable $(t, x)$.

\begin{lemm}\label{lemm4.1}
Let $0<\lambda<1$. Then there exists a positive constant
$C_\lambda\neq 0$ such that for any $f\in \cS(\RR^n)$, one has
\begin{equation}\label{4.6}
|D_{y}|^{\lambda} f(y)=\cF^{-1}
\big(|\xi|^{\lambda}\hat{f}(\xi)\big)
=C_\lambda\int_{\RR^n}\frac{f(y)-f(y+h)}{|h|^{n+\lambda}} d h.
\end{equation}
\end{lemm}
Indeed, note that
$$
\int_{\RR^n}\frac{f(y)-f(y+h)}{|h|^{n+\lambda}} d h= \int_{\RR^n}
\hat{f}(\xi)e^{i y \cdot \xi}\int_{\RR^n}\frac{1-e^{i h \cdot
\xi}}{|h|^{n+\lambda}} d h d\xi,
$$
while
\[
  \int_{\RR^n}\frac{1-e^{-i\, h\cdot\xi}}{|
   h|^{n+\lambda}}\,d h=|\xi|^\lambda \int_{\RR^n}
  \frac{1-e^{-i\, u\cdot\frac{\xi}{|\xi|}}}{|
  {u}|^{n+\lambda}}\,d u ,
\]
so that (\ref{4.6}) follows from
$$
\int_{\RR^n}
  \frac{1-e^{-i\, u\cdot\frac{\xi}{|\xi|}}}{|
  {u}|^{n+\lambda}}\,d u \neq 0,
$$
which is a positive constant depending only on $\lambda$ and the
dimension $n$, but independent from $\xi$.

Using this Lemma, we have the following Leibniz type formula,
\begin{eqnarray}\label{4.7}
&&|D_{y}|^{\lambda} \big(f(y)g(y)\big)
=C_\lambda\int_{\RR^n}\frac{f(y)g(y)-f(y+h)g(y+h)}{|h|^{n+\lambda}} d h\\
&=& g(y)|D_{y}|^{\lambda} f(y) + f(y)|D_{y}|^{\lambda} g(y)+
C_\lambda\int_{\RR^n}\frac{\big(f(y)-f(y+h)\big)
\big(g(y+h)-g(y)\big)}{|h|^{n+\lambda}} d h. \nonumber
\end{eqnarray}

We now turn to the analysis of the fractional derivative with
respect to $(t, x)$ of the nonlinear collision operator.
Denote the difference with respect to $(t, x)$ by
$$
f_h(t, x, v)={f}(t, x, v)-{f}((t, x)+h\, , v),\,\,\,\,\,\,
h\in\RR^4_{t, x}\, .
$$
It follows that for the collision operator (where $n=1+3$),
\begin{eqnarray}\label{4.8}
\enskip |D_{t, x}|^{\lambda}Q\big({f},\,\, g\big)=
Q\big(|D_{t, x}|^{\lambda}{f},\,\, g\big)+ Q\big({f},\,\,
|D_{t, x}|^{\lambda}g\big) +C_\lambda\int_{\RR^4}
|h|^{-4-\lambda}Q\big(f_h,\,\, g_h\big) d h.
\end{eqnarray}

\bigbreak This kind of decomposition will be used extensively below
in order to get the partial regularity with respect to the $(t, x)$ variable.

\smallbreak
First of all, we have the following  proposition on the gain of regularity
in the variable $(t,x)$ through the uncertainty principle.

\begin{prop}\label{prop4.1}
Under the hypothesis of Theorem \ref{theo1}, one has
\begin{equation}\label{4.1}
\Lambda^{s_0}_{t, x}\, f_1\in H^5_l  (\RR^7),
\end{equation}
for any $l  \in\NN$ and $0<s_0=\frac{s(1-s)}{(s+1)}$.
\end{prop}

\noindent{\bf Proof:} In fact, for any $l\in\NN$, it follows from
Proposition \ref{prop3.1} that
$$
\Lambda^s_v W_l g \in L^2(\RR^7).
$$
Then the upper bound estimation given by Corollary \ref{lemm104}
with $m=-s$ implies that
$$
W_l Q( f_2,\,\, g)\in L^2(\RR^4_{t, x};\,\, H^{-s}(\RR^3_v)).
$$
On the other hand, Proposition \ref{prop3.1+0} and \eqref{3.2} gives
$$
W_l G\in L^2(\RR^4_{t, x};\,\, H^{-(2s-1+\delta)}(\RR^3_v)).
$$
By using (\ref{3.1}), it follows that
\begin{equation}\label{4.1+1-}
\partial_t(W_l g) + v\,\cdot\,\partial_x {(W_l g )} =W_l Q(f_2,\,\, g)+
W_l G \in H^{-s}(\RR^7).
\end{equation}
Finally, by using Lemma \ref{lemm2.1} with $s'=s$, we
can conclude (\ref{4.1}) and this completes the proof of the proposition.\\

Therefore,  under the hypothesis $f\in \cH^5_l (]T_1,
T_2[\times\Omega\times\RR^3_{v})$ for all $l  \in\NN$, it follows
that for any $ l \in\NN$ we have
\begin{equation}\label{4.2}
\Lambda^{s}_{v}(\varphi(t)\psi(x) f) \, \in H^5_l(\RR^7),\,\, \hskip
0.5cm \Lambda^{s_0}_{t, x}(\varphi(t)\psi(x) f) \, \in H^5_l
(\RR^7)\, .
\end{equation}

\smallbreak We now improve this partial regularity in $(t,
x)$ variable.

\begin{prop}\label{prop4.2}
Let $0<\lambda< 1$. Suppose that $f\in {\cH}^5_l (]T_1,\, T_2[\,
\times\Omega\times\RR^3_{v})$ is a solution of the equation
(\ref{1.1}) for all $l \in\NN$. Furthermore, assume that for any cutoff functions
$\varphi, \psi$, 
\begin{equation}\label{4.1+0}
\Lambda^{s}_{v}(\varphi(t)\psi(x) f) \, \in H^5_l(\RR^7),\,\, \hskip
0.5cm \Lambda^{\lambda}_{t, x}(\varphi(t)\psi(x) f)\in H^5_l(\RR^7).
\end{equation}
Then, one has
\begin{equation}\label{4.1+1}
\Lambda^s_v\Lambda^{\lambda}_{t, x}
(\varphi(t)\psi(x) f)\in H^5_l(\RR^7),
\end{equation}
for any $l  \in\NN$ and  any cutoff functions $\varphi, \psi$.
\end{prop}

\noindent{\bf Proof:} Set
$$
g_{N, l}=P_{N, l}\, g=\psi_1(x) S_N(D_x)\, W_{l  }\,\,
S_N(D_v)\partial^\alpha(\varphi(t)\psi(x) f),
$$
where $ \alpha\in\NN^7, |\alpha|\leq 5$ and $ l \in\NN$.
Then \eqref{4.1+0} yields
$$
 \|\Lambda^{s}_{v} g_{N, l}\|_{L^2(\RR^7)}\leq C
\|\Lambda^{s}_{v}\partial^\alpha (\varphi(t)\psi(x) f)
\|_{L^2_l(\RR^7)},
$$
and
$$
 \|\Lambda^{\lambda}_{t, x} g_{N, l}\|_{L^2(\RR^7)}\leq C
\|\Lambda^{\lambda}_{t, x}\partial^\alpha (\varphi(t)\psi(x) f)
\|_{L^2_l(\RR^7)},
$$
where the constant $C$ is independent of $N$.

\smallbreak
It follows that $g_{N,l}$ satisfies the equation
\begin{equation}\label{4.3}
\partial_t(g_{N, l}) + v\,\cdot\,\partial_x \,(g_{N, l}) = Q({f}_2,\,\, g_{N, l})
+ G_{N, l},
\end{equation}
where $G_{N, l}$ is given by
\begin{eqnarray*}
 G_{N, l}&=& \psi_1(x) W_{ l}\, \Big[S_N(D_v),\,\, v\Big]\,\cdot\,\nabla_x S_N(D_x) g
+\Big(P_{N, l}\,  Q\big({f}_2,\,\, g\big)-Q\big({f}_2,\,\,
P_{N, l}\,  g\big)\Big)\\
&&+\Big((v\,\cdot\,\nabla_x) \psi_1(x)\Big) W_{ l}\, S_N(D_x)
S_N(D_v) g + P_{N, l}\, G,\nonumber
\end{eqnarray*}
with $G$ defined in (\ref{3.1**}).

We now choose $ |D_{t, x}|^{\lambda} \psi_2^2(x)
|D_{t, x}|^{\lambda} g_{N, l} $ as a test function for equation
(\ref{4.3}). It follows that
\begin{eqnarray}\label{4.5}
&\Big( v\,\cdot\,\big(\partial_x \psi_2\big) {|D_{t, x}|^{\lambda}
g_{N, l}},\,\, \psi_2(x) |D_{t, x}|^{\lambda} g_{N, l}\Big)_{L^2(\RR^7)}&\\
& =\Big( \psi_2(x) |D_{t, x}|^{\lambda}\big\{Q({f}_2,\,\, g_{N,
l})+G_{N, l}\big\},\,\, \psi_2(x) |D_{t, x}|^{\lambda} g_{N,
l}\Big)_{L^2(\RR^7)}.& \nonumber
\end{eqnarray}

It is sufficient to prove that, for any
$l\in\NN$,
\begin{equation}\label{4.1+2}
\Lambda^{s}_v \Lambda^{\lambda}_{t, x} P_{N, l }\,g \in L^2(\RR^7),
\end{equation}
and is  uniformly  bounded with respect to $N$. In the rest of the
proof, we use
$C$ to denote a constant  independent of $N$.

\smallbreak

We first consider the linear terms in (\ref{4.5}).
On the left hand side of
(\ref{4.5}), the hypothesis (\ref{4.1+0}) implies that
\begin{eqnarray*}\label{4.1+3}
\Big\| v\,\cdot\,\big(\partial_x \psi_2\big) {|D_{t, x}|^{\lambda}
g_{N, l}}\Big\|_{L^2(\RR^7)} \leq C \|\,\, |\Lambda_{t, x}|^{\lambda}
\partial^\alpha (\varphi(t)\psi(x) f)
\|_{L^2_{l+1}(\RR^7)}.
\end{eqnarray*}
For the linear terms in $G_{N, l}$, by using \eqref{bounded}, one has
\begin{eqnarray*}\label{4.1+4}
 &\left\| \psi_2(x) |D_{t, x}|^{\lambda}\big\{\psi_1(x) W_{ l}\,
\big[S_N(D_v),\,\, v\big]\,\cdot\,\nabla_x S_N(D_x)
g\big\}\right\|_{L^2(\RR^7)}&\\
 &\leq C \| |\Lambda_{t, x}|^{\lambda}
\partial^\alpha (\varphi(t)\psi(x) f)
\|_{L^2_{l}(\RR^7)},\nonumber&
\end{eqnarray*}
and
\begin{eqnarray*}\label{4.1+5}
&& \left\| \psi_2(x) |D_{t, x}|^{\lambda}\big(v\,\cdot\,(\nabla_x\,
\psi_1)(x)\big) W_{ l}\, S_N(D_x) S_N(D_v) g\right\|_{L^2(\RR^7)}\\
&& \,\,\,\,\,\,\,\,\,\,\,\,\,\,\,\,\,\,\leq C \|
|\Lambda_{t, x}|^{\lambda}
\partial^\alpha (\varphi(t)\psi(x) f)
\|_{L^2_{l+1}(\RR^7)}.\nonumber
\end{eqnarray*}
Similarly, concerning the linear terms $(B)$ and $(C)$ in $G$, we have
\begin{eqnarray*}\label{4.1+6}
&& \left\| \psi_2(x) |D_{t, x}|^{\lambda}P_{N, l}\big((B)+
(C)\big)\right\|_{L^2(\RR^7)}
\leq C \|
|\Lambda_{t, x}|^{\lambda}
\partial^\alpha (\varphi(t)\psi(x) f)
\|_{L^2_{l+1}(\RR^7)}.
\end{eqnarray*}

For the nonlinear terms in (\ref{4.5}),  we shall
use  the formula (\ref{4.8}).
First of all, the coercivity estimate (\ref{E-sub-estimate+}) gives,
 as in (\ref{3.7}), that
\begin{eqnarray}\label{4.1+7}
&&-\Big(Q(f_2,\, \psi_1(x)|D_{t, x}|^{\lambda} g_{N, l}),\,\,
\psi_1(x)|D_{t, x}|^{\lambda} g_{N, l}\Big)_{L^2(\RR^7)}
\\
&&\geq C_0\|\Lambda^{s}_v W_{\gamma/2}\psi_1(x)|D_{t, x}|^{\lambda}
g_{N, l}\|^2_{L^2(\RR^7)}\nonumber \\
&&- C \| f_2\|_{L^\infty(\RR^4_{t, x};\,\,
L^1_{\max\{\gamma^+,\,
2-\gamma^+\}}(\RR^3_v))}\|\psi_1(x)|D_{t, x}|^{\lambda} g_{N,
l}\|^2_{L^2_{\gamma^+/2}(\RR^7)}. \nonumber
\end{eqnarray}

On the other hand, the upper estimate of Theorem \ref{theo2.1} with $m=-s$ and
$\a=-\gamma/2>0$ (the case $\gamma>0$ is easier) gives,
\begin{align*}
&\Big|\Big(Q(|D_{t, x}|^{\lambda}f_2,\, \psi_1(x) g_{N, l}),\,\,
\psi_1(x)|D_{t, x}|^{\lambda} g_{N, l}\Big)_{L^2(\RR^7)}\Big|
\\
&\leq C\||D_{t, x}|^{\lambda} f_2\|_{L^\infty(\RR^4_{t, x},\,
L^1_{|\gamma|/2+\gamma^+ +2s}(\RR^3_v))}
 { \|\psi_1(x)\Lambda^{s}_v g_{N, l}\|_{L^2_{|\gamma|/2+\gamma^+ +2s}(\RR^7)}}
\|\psi_1(x)|D_{t, x}|^{\lambda}\Lambda^{s}_v W_{\gamma/2}g_{N, l}\|_{L^2(\RR^7)}
\nonumber
\\
&\leq \varepsilon \|\psi_1(x)|D_{t, x}|^{\lambda}\Lambda^{s}_v W_{\gamma/2} g_{N,
l}\|^2_{L^2(\RR^7)} +C_\varepsilon \| |D_{t, x}|^{\lambda}
f_2\|^2_{L^\infty(\RR^4_{t, x},\, L^2_{|\gamma|/2+\gamma^+
+2s+4}(\RR^3_v))}
 { \|\Lambda^{s}_v g\|^2_{L^2_{|\gamma|/2+\gamma^+ +2s+l}(\RR^7)}} .\nonumber
\end{align*}
For  the cross term coming from the decomposition
(\ref{4.8}), by using again estimate (\ref{2.2+001}) with $m=-s$ and $\a=|\gamma|/2$, we
get
\begin{eqnarray*}
&&\int_{\RR^4} |h|^{-4-\lambda}\Big|\Big(Q((f_2)_h,\, (g_{N, l})_h),\,\,
\psi^2_1(x)|D_{x}|^{\lambda} g_{N, l}\Big)_{L^2(\RR^7)}d h\Big|
\\
&\leq& |C_\lambda| \|\psi_1(x)|D_{x}|^{\lambda}\Lambda^{s}_v W_{\gamma/2}g_{N,,
l}\|_{L^2(\RR^7)}\\
&&\,\,\,\,\,\,\,\,\,\times \int_{\RR^4} |h|^{-4-\lambda} \| (f_2)_h\|
_{L^\infty(\RR^4_{t, x},\, L^1_{|\gamma|/2+\gamma^+ +2s}(\RR^3_v))}
\|\Lambda^{s}_v (g_{N, l})_h\|_{L^2_{|\gamma|/2+\gamma^+ +2s}(\RR^7)}
 d h .
\end{eqnarray*}
Furthermore, 
\begin{eqnarray*}
&& \int_{\RR^4} |h|^{-4-\lambda} \| (f_2)_h \| _{L^\infty(\RR^4_{t, x},\,
L^1_{|\gamma|/2+\gamma^+ +2s}(\RR^3_v))} \|\Lambda^{s}_v (g_{N, l})_h
\|_{L^2_{|\gamma|/2+\gamma^+
+2s}(\RR^7)}
 d h\\
 &\leq& \int_{|h|<1} |h|^{-4-\lambda} \|(f_2)_h \|
_{L^\infty(\RR^4_{t, x},\, L^1_{|\gamma|/2+\gamma^+ +2s}(\RR^3_v))}
\|\Lambda^{s}_v (g_{N, l})_h\|_{L^2_{|\gamma|/2+\gamma^+ +2s}(\RR^7)} d h \\
&&+ 4\tilde{C}_\lambda\| f_2\| _{L^\infty(\RR^4_{t, x},\,
L^1_{|\gamma|/2+\gamma^+ +2s}(\RR^3_v))}
\|\Lambda^{s}_v g_{N, l}\|_{L^2_{|\gamma|/2+\gamma^+ +2s}(\RR^7)}\\
& \leq& 2\int_{|h|<1} |h|^{-4-\lambda+1} \| \nabla_{t, x} f_2\|
_{L^\infty(\RR^4_{t, x},\, L^1_{|\gamma|/2+\gamma^+ +2s}(\RR^3_v))}
\|\Lambda^{s}_v g_{N, l}\|_{L^2_{|\gamma|/2+\gamma^+ +2s}(\RR^7)} d h \\
&&+ 4\tilde{C}_\lambda\| f_2\| _{L^\infty(\RR^4_{t, x},\,
L^1_{|\gamma|/2+\gamma^+ +2s}(\RR^3_v))} \|\Lambda^{s}_v g_{N,
l}\|_{L^2_{|\gamma|/2+\gamma^+ +2s}(\RR^7)}.
\end{eqnarray*}
Thus
\begin{eqnarray*}
&&\int_{\RR^4} |h|^{-4-\lambda}\Big|\Big(Q((f_2)_h,\, (g_{N, l})_h),\,\,
\psi^2_1(x)|D_{x}|^{\lambda}W_{\gamma/2} g_{N, l}\Big)_{L^2(\RR^7)}d h\Big|
\\
&\leq& \varepsilon \|\psi_1(x)|D_{x}|^{\lambda}\Lambda^{s}_v g_{N,
l}\|^2_{L^2(\RR^7)} +C_\varepsilon \| \Lambda^1_{t, x}
f_2\|^2_{L^\infty(\RR^4_{t, x}; L^2_{|\gamma|/2+\gamma^+ +
2s+4}(\RR^3_v))} \|\Lambda^{s}_v  g_{N, l}\|^2_{L^2_{|\gamma|/2+\gamma^+
+2s}(\RR^7)} .\nonumber
\end{eqnarray*}
Hence, the  formula (\ref{4.8}) yields
\begin{align*}
&\Big|\Big(|D_{t, x}|^{\lambda}Q(f_2,\, \psi_1(x) g_{N, l})-
Q(|D_{t, x}|^{\lambda}f_2,\, \psi_1(x) g_{N, l}),\,\,
\psi_1(x)|D_{t, x}|^{\lambda} g_{N, l}\Big)_{L^2(\RR^7)}\Big|
\\&
\leq \varepsilon \|\psi_1(x)|D_{t, x}|^{\lambda}\Lambda^{s}_v W_{\gamma/2} g_{N,
l}\|^2_{L^2(\RR^7)}+C_\varepsilon \| \Lambda^1_{t, x}
f_2\|^2_{L^\infty(\RR^4_{t, x},\, L^2_{|\gamma|/2+\gamma^+
+2s+4}(\RR^3_v))} \|\Lambda^{s}_v g\|^2_{L^2_{|\gamma|/2+\gamma^+ +2s+l}(\RR^7)} .\nonumber
\end{align*}

In conclusion,  we get from coercivity property \eqref{4.1+7} that 
\begin{align}\label{4.1+10}
\|\Lambda^{s}_v W_{\gamma/2}&\psi_1(x)|D_{t, x}|^{\lambda} g_{N,
l}\|^2_{L^2(\RR^7)}
\\ \notag
\leq& C \| \Lambda^1_{t, x} f_2\|^2_{L^\infty(\RR^4_{t, x},\,
L^2_{|\gamma|/2+\gamma^+ +2s+4}(\RR^3_v))}
\\&
\hspace*{1cm}\notag
\times
\Big( \|\,\,|D_{t, x}|^{\lambda}
g\|^2_{L^2_{l+|\gamma|/2+\gamma^+ +2s}(\RR^7)}+\|\Lambda^{s}_v
g\|^2_{L^2_{l+|\gamma|/2+\gamma^+ +2s}(\RR^7)}
\Big)\nonumber
\\
& +\left|\left( |D_{t, x}|^{\lambda} \Big(P_{N, l}\,
Q\big(f_2,\,\, g\big)-Q\big(f_2,\,\, P_{N, l}\,
g\big)\Big),\,\, \psi^2_2(x)|D_{t, x}|^{\lambda}g_{N, l}
\right)_{L^2(\RR^7)}\right| \nonumber \\
&+\left|\left(|D_{t, x}|^{\lambda} P_{N, l}\,\, (A),\,\,
\psi^2_2(x)|D_{t, x}|^{\lambda} \, g_{N, l}\right)_{L^2(\RR^7)}\right|
\nonumber\\
&= {\rm (I)}+{\rm (II)}+{\rm (III)}\, . \nonumber
\end{align}
For the term (II),  since $[|D_{t, x}|^\lambda,\,\, \psi_1(x)]$ is a
bounded operator, we can replace $P_{N, l}$ by $\tilde{P}_{N, l}=
W_l\,\,S_N(D_x)S_N(D_v)$. Again, the formula (\ref{4.8}) yields
\begin{eqnarray*}
&&\Big(|D_{t, x}|^{\lambda} \big( \tilde{P}_{N, l} Q\big(f_2,\,\,
g\big)-Q\big(f_2,\,\, \tilde{P}_{N, l} g\big)\big) ,\,\,
\psi^2_2(x)|D_{t, x}|^{\lambda} g_{N, l}\Big)_{L^2(\RR^7)}\\
&&=\Big(\big( \tilde{P}_{N, l} Q\big(|D_{t, x}|^{\lambda}f_2,\,\,
g\big)-Q\big(|D_{t, x}|^{\lambda}f_2,\,\, \tilde{P}_{N, l}
g\big)\big) ,\,\,
\psi^2_2(x)|D_{t, x}|^{\lambda} g_{N, l}\Big)_{L^2(\RR^7)}\\
&&+\Big(\big( \tilde{P}_{N, l} Q\big(f_2,\,\,
|D_{t, x}|^{\lambda} g\big)-Q\big(f_2,\,\, \tilde{P}_{N, l}
|D_{x}|^{\lambda} g\big)\big) ,\,\, \psi^2_2(x)|D_{t, x}|^{\lambda}
g_{N, l}\Big)_{L^2(\RR^7)}
\\
&&+C_\lambda\int_{\RR^4} |h|^{-4-\lambda}\Big( \big( \tilde{P}_{N, l}
Q\big((f_2)_h,\,\, g_h\big)-Q\big((f_2)_h,\,\, \tilde{P}_{N,
l} g_h\big)\big) ,\,\, \psi^2_2(x)|D_{t, x}|^{\lambda} g_{N,
l}\Big)_{L^2(\RR^7)}d h.
\end{eqnarray*}
As for (\ref{3.12}), in the case 
when $1/2\leq
s<1$ (the other case when $0<s<1/2$ is similar and  easier to
handle), by applying Lemmas \ref{lemm3.4}, \ref{lemm3.1} and 
\ref{lemm3.3}, we have 
\begin{eqnarray}\label{4.1+11}
&&\Big|\Big(\big( \tilde{P}_{N, l}
Q\big(|D_{t, x}|^{\lambda}f_2,\,\,
g\big)-Q\big(|D_{t, x}|^{\lambda}f_2,\,\, \tilde{P}_{N, l}
g\big)\big) ,\,\,
\psi^2_2(x)|D_{t, x}|^{\lambda} g_{N, l}\Big)_{L^2(\RR^7)}\Big|\\
&&\leq C\|\Lambda_{t, x}^{1+\lambda}f_2\|_{L^\infty(\RR^4_{t, x}\,,
\,\, L^1_{(\gamma+2s-1)^+}(\RR^3_v))} \|g\|_{L^2(\RR^4_{t, x}\,,
\,\, H^{2s-1+\delta}_{(\gamma+2s-1)^+}(\RR^3_v))} \|\,|D_{t,
x}|^{\lambda} g\|_{L^2_{2l }(\RR^7)}. \nonumber
\end{eqnarray}

By using (\ref{3.8-1*}) of Lemma \ref{lemm3.1},  we can get, for
$2s-1+\delta<s$,
\begin{eqnarray*}\label{4.1+12+11}
&&\Big|\Big(\big( \tilde{P}_{N, l} Q\big(f_2,\,\,
|D_{t, x}|^{\lambda} g\big)-Q\big(f_2,\,\, \tilde{P}_{N, l}
|D_{t, x}|^{\lambda} g\big)\big) ,\,\,
\psi^2_2(x)|D_{t, x}|^{\lambda} g_{N, l}\Big)_{L^2(\RR^7)}\Big|\\
&&\leq C \|\Lambda_{t, x} f_2\|_{L^\infty(\RR^4_{t, x}\,, \,\,
L^1_{(\gamma+2s-1)^+}(\RR^3_v))} \|\,\,|D_{t, x}|^{\lambda}
g\|_{L^2(\RR^4_{t, x}\,, \,\, L^2_{l+(\gamma+2s-1)^+}(\RR^3_v))}\\
&&\hskip5cm \times\|\,|D_{t, x}|^{\lambda}\psi_1 g_{N,
l}\|_{L^2(\RR^4_{t, x}\, , \,\,
H^{2s-1+\delta}_{l+(\gamma+2s-1)^+}(\RR^3_v))} \nonumber\\
&&\leq \varepsilon \|\Lambda^{s}_v W_{\gamma/2}\psi_1(x)|D_{t,
x}|^{\lambda} g_{N, l}\|^2_{L^2(\RR^7)}\\
&& \hskip2cm+ C_\varepsilon \|\Lambda^1_{t,
x}f_2\|^{2k'}_{L^\infty(\RR^4_{t, x}\, ,\,\,
L^2_{l+3/2+\delta+(\gamma+2s-1)^+}(\RR^3_v))}\|\,\,|D_{t,
x}|^{\lambda} g\|^{2k'}_{L^2(\RR^4_{t, x}\,, \,\,
L^2_{kl}(\RR^3_v))}\, , \nonumber
\end{eqnarray*}
and
\begin{eqnarray*}
&&\Big|\int |h|^{-4-\lambda}\Big( \big( \tilde{P}_{N, l}
Q\big(f_{2, h},\,\, g_h\big)-Q\big(f_{2, h},\,\, \tilde{P}_{N,
l} g_h\big)\big) ,\,\, \psi^2_2(x)|D_{x}|^{\lambda} g_{N,
l}\Big)_{L^2(\RR^7)}d h\Big|\\
&&\,\,\,\,\,\,\,\,\,\leq C \|\Lambda_{t, x}\Lambda_x
f_2\|_{L^\infty(\RR^4_{t, x}; L^1_{l+(\gamma+2s-1)^+}(\RR^3_v))}\|
g\|_{L^2(\RR^4_{t, x}\,, \,\,
H^{2s-1+\delta}_{l+(\gamma+2s-1)^+}(\RR^3_v))}
\|\,|D_{t, x}|^{\lambda} g\|_{L^2_{l  }(\RR^7)} \\
&&\,\,\,\,\,\,\,\,\,\leq C \| f_2\|_{H^{2+4/2+\delta}_{l+3/2+\delta+
(\gamma+2s-1)^+}(\RR^7)}\|\Lambda_v^s g\|_{L^2_{l+(\gamma+2s-1)^+}(\RR^7)}
\|\,|D_{t, x}|^{\lambda} g\|_{L^2_{l  }(\RR^7)}  .
\end{eqnarray*}
Thus, we have
\begin{eqnarray*}
{\rm (II)}&\leq&\varepsilon\|\Lambda^{s}_v W_{\gamma/2}
\psi_1(x)|D_{t, x}|^{\lambda} g_{N,
l}\|^2_{L^2(\RR^7)}\\
&+& C_\varepsilon \| f_2\|^{2k'}_{H^2_{l+\gamma^+ +2s+4}(\RR^3_v))}
\Big( \|\,\,|D_{t, x}|^{\lambda} g\|^{2k'}_{L^2_{kl+\gamma^+ +2s}(\RR^7)}+\|\Lambda^{s}_v
g\|^2_{L^2_{l+\gamma^+ +2s}(\RR^7)} \Big)\nonumber.
\end{eqnarray*}

We now consider  the last term (III) of (\ref{4.1+10}). Recall
that  $(A)$
stands for the nonlinear terms from $G$ given in
(\ref{3.1**}). Precisely
\begin{eqnarray*}
(A)=\sum_{\alpha_1+\alpha_2=\alpha\,\, \alpha_1\neq 0}C^{\alpha_1}_{\alpha_2}
Q\Big(\partial^{\alpha_1} f_2,\,\,
\partial^{\alpha_2} f_1\Big).
\end{eqnarray*}
By using  (\ref{2.2+001}) (we consider also only
the case $1/2\leq s<1$) and formula \eqref{4.8}, we have 
\begin{eqnarray*}
&&\left|\left(|D_{t, x}|^{\lambda} \,
\Big(Q\Big(\partial^{\alpha_1} f_2,\,\,
\partial^{\alpha_2} f_1\Big)\Big),\,\,P_{N, l}
\psi^2_2(x)|D_{t, x}|^{\lambda} \, g_{N,
l}\right)_{L^2(\RR^7)}\right|\\
&\leq& C\|\Lambda^{-m}_v W_{\gamma/2} \psi_1(x)|D_{t, x}|^{\lambda}
g_{N, l}\|_{L^2(\RR^7)}\Big\{\Big\|Q\Big(|D_{t,
x}|^\lambda\partial^{\alpha_1} f_2,\,\,
\partial^{\alpha_2} f_1\Big)\Big\|_{L^2(\RR^4_{t, x}; H^{m}_{l+|\gamma|/2}(\RR^3_v))}\\
&&+\Big\|Q\Big(\partial^{\alpha_1} f_2,\,\, |D_{t, x}|^\lambda
\partial^{\alpha_2} f_1\Big)\Big\|_{L^2(\RR^4_{t, x};
H^{m}_{l+|\gamma|/2}(\RR^3_v))}\\
&&+\Big\|\int h^{-4-\lambda}Q\Big(\partial^{\alpha_1} (f_2)_h,\,\,
\partial^{\alpha_2}(f_1)_h\Big) dh\Big\|_{L^2(\RR^4_{t, x};
H^{m}_{l+|\gamma|/2}(\RR^3_v))}\Big\}.
\end{eqnarray*}

We divide the discussion into two cases.

\noindent{\bf Case 1.}  $|\alpha_1|=1, 2$. Take $m=-s$. We have
\begin{eqnarray*}
&&\Big\|Q\Big(|D_{t, x}|^\lambda\partial^{\alpha_1} f_2,\,\,
\partial^{\alpha_2} f_1\Big)\Big\|_{L^2(\RR^4_{t, x}; H^{-s}_{l+
|\gamma|/2}(\RR^3_v))}+\Big\|Q\Big(\partial^{\alpha_1} f_2,\,\,
|D_{t, x}|^\lambda \partial^{\alpha_2} f_1\Big)\Big\|_{L^2(\RR^4_{t,
x};
H^{-s}_{l+|\gamma|/2}(\RR^3_v))}\\
&&\leq C\|\,\Lambda_{t, x}^{\lambda}
\partial^{\alpha_1} f_2\|_{L^\infty(\RR^4_{t, x};
L^1_{l+\gamma^+ +2s}(\RR^3_v))}\|\,\Lambda_v^{s}\Lambda_{x}^{\lambda}
\partial^{\alpha_2} f_1\|_{L^2_{l +\gamma^+ +2s}(\RR^7)}
\\
&&\leq C\|f_2\|_{H^{\lambda+2+4/2+\delta}_{l+3/2+\delta+\gamma^+
+2s}(\RR^7)}\|\Lambda_v^s f_1\|_{H^{4+\lambda}_{l +\gamma^+ +2s}(\RR^7)},
\end{eqnarray*}
and
\begin{eqnarray*}
&&\Big\|\int_{\RR^4}h^{-4-\lambda}Q\Big(\partial^{\alpha_1} (f_2)_h,\,\,
\partial^{\alpha_2} (f_1)_h\Big) dh\Big\|_{L^2(\RR^4_{t, x};
H^{-s}_{l+|\gamma|/2}(\RR^3_v))}\\
&&\leq C\int |h|^{-4-\lambda}\|\,\partial^{\alpha_1} (f_2)_h\|_{L^\infty(\RR^4_{t, x};
L^1_{l+\gamma^+ +2s}(\RR^3_v))}\|\,\Lambda_v^{s}
\partial^{\alpha_2} (f_1)_h\|_{L^2_{l +\gamma^+ +2s}(\RR^7)}dh\\
&&\leq C\|\,\partial^{\alpha_1} f_2\|_{L^\infty(\RR^4_{t, x};
L^1_{l+\gamma^+ +2s}(\RR^3_v))}\|\,\Lambda_v^{s}
\partial^{\alpha_2} \nabla_{t, x}f_1)\|_{L^2_{l +\gamma^+ +2s}(\RR^7)}
\\
&&\leq C\|f_2\|_{H^{2+4/2+\delta}_{l+3/2+\delta+\gamma^+
+2s}(\RR^7)}\|\Lambda_v^s f_1\|_{H^{5}_{l +\gamma^+ +2s}(\RR^7)}.
\end{eqnarray*}

\noindent{\bf Case 2.}  $|\alpha_1|\geq 3$. By the same argument
as above, one has
\begin{eqnarray*}
&&\Big\|Q\Big(|D_{t, x}|^\lambda\partial^{\alpha_1} f_2,\,\,
\partial^{\alpha_2} f_1\Big)\Big\|_{L^2(\RR^4_{t, x}; H^{-s}_{l+
|\gamma|/2}(\RR^3_v))}+\Big\|Q\Big(\partial^{\alpha_1} f_2,\,\,
|D_{t, x}|^\lambda \partial^{\alpha_2} f_1\Big)\Big\|_{L^2(\RR^4_{t,
x};
H^{-s}_{l+|\gamma|/2}(\RR^3_v))}\\
&&\leq C\|\,\Lambda_{t, x}^{\lambda}
\partial^{\alpha_1} f_2\|_{L^2(\RR^4_{t, x};
L^1_{l+\gamma^+ +2s}(\RR^3_v))}\|\,\Lambda_v^{s}\Lambda_{t, x}^{\lambda}
\partial^{\alpha_2} f_1\|_{L^\infty(\RR^4_{t, x}; L^2_{l +\gamma^+ +2s}(\RR^3_v))}
\\
&&\leq C\|\Lambda_{t, x}^{\lambda}f_2\|_{H^{5}_{l+3/2+\delta+\gamma^+
+2s}(\RR^7)}\|\Lambda_v^s f_1\|_{H^{2+4/2+\lambda+\delta}_{l +\gamma^+ +2s}(\RR^7)}.
\end{eqnarray*}
When $|\alpha_1|= 3, 4$, we have
\begin{eqnarray*}
&&\Big\|\int_{\RR^4}h^{-4-\lambda}Q\Big(\partial^{\alpha_1} (f_2)_h,\,\,
\partial^{\alpha_2} (f_1)_h\Big) dh\Big\|_{L^2(\RR^4_{t, x}; H^{-s}_{l+|\gamma|/2}(\RR^3_v))}\\
&&\leq C\int |h|^{-4-\lambda}\|\,\partial^{\alpha_1} (f_2)_h\|_{L^2(\RR^4_{t, x};
L^1_{l+\gamma^+ +2s}(\RR^3_v))}\|\,\Lambda_v^{s}
\partial^{\alpha_2} (f_1)_h\|_{L^\infty(\RR^4_{t, x}; L^2_{l +\gamma^+ +2s}(\RR^3_v))}dh\\
&&\leq C\|\,\nabla_{t, x}\partial^{\alpha_1} f_2\|_{L^2(\RR^4_{t, x};
L^1_{l+\gamma^+ +2s}(\RR^3_v))}\|\,\Lambda_v^{s}
\partial^{\alpha_2} f_1)\|_{L^\infty(\RR^4_{t, x}; L^2_{l +\gamma^+ +2s}(\RR^3_v))}
\\
&&\leq C\|f_2\|_{H^{5}_{l+3/2+\delta+\gamma^+ +2s}(\RR^7)}
\|f_1\|_{H^{2+4/2+s+\delta}_{l +\gamma^+ +2s}(\RR^7)},
\end{eqnarray*}
while when $|\alpha_1|= |\alpha|=5$, we have
\begin{eqnarray*}
&&\Big\|\int_{\RR^4}h^{-4-\lambda}Q\Big(\partial^{\alpha} (f_2)_h,\,\,
(f_1)_h\Big) dh\Big\|_{L^2(\RR^4_{t, x}; H^{-s}_{l+|\gamma|/2}(\RR^3_v))}\\
&&\leq C\int |h|^{-4-\lambda}\|\,\partial^{\alpha} (f_2)_h\|_{L^2(\RR^4_{t, x};
L^1_{l+\gamma^+ +2s}(\RR^3_v))}\|\,\Lambda_v^{s}
(f_1)_h\|_{L^\infty(\RR^4_{t, x}; L^2_{l +\gamma^+ +2s}(\RR^3_v))}dh\\
&&\leq C\|\,\partial^{\alpha} f_2\|_{L^2(\RR^4_{t, x};
L^1_{l+\gamma^+ +2s}(\RR^3_v))}\|\,\Lambda_v^{s}
\nabla_{t, x} f_1)\|_{L^\infty(\RR^4_{t, x}; L^2_{l +\gamma^+ +2s}(\RR^3_v))}
\\
&&\leq C\|f_2\|_{H^{5}_{l+3/2+\delta+\gamma^+
+2s}(\RR^7)}\|f_1\|_{H^{1+4/2+s+\delta}_{l +\gamma^+ +2s}(\RR^7)},
\end{eqnarray*}

Thus, by the Cauchy-Schwarz inequality, we obtain
\begin{eqnarray*}
{\rm (III)}\leq\varepsilon\|\Lambda^{s}_v W_{\gamma/2}
\psi_1(x)|D_{t, x}|^{\lambda} g_{N, l}\|^2_{L^2(\RR^7)}+ C_\varepsilon
\Big(\|\,\Lambda_{t, x}^{\lambda} f_3\|^{4}_{H^5_{2l+\gamma^+ +7 }(\RR^7)}
+\|\,\Lambda^{s} f_3\|^{4}_{H^5_{2l+\gamma^+ +7 }(\RR^7)}\Big).
\end{eqnarray*}
Finally, we get from (\ref{4.1+10}) that 
\begin{eqnarray*}\label{4.1+16}
&&\|\Lambda^{s}_v W_{\gamma/2} \psi_1(x)|D_{t, x}|^{\lambda} g_{N,
l}\|^2_{L^2(\RR^7)}\leq C\Big(\|\,\Lambda_{t, x}^{\lambda} f_3\|^{k'4}_{H^{5}_{kl
+\gamma^+ +7}(\RR^7)}+\|\,\Lambda_{v}^{s} f_3\|^4_{H^5_{2l+\gamma^+ +7
}(\RR^7)}\Big).
\end{eqnarray*}
Therefore, we complete the proof for
Proposition \ref{prop4.2}.\\

\smallbreak
We are now ready to prove the following regularity result on the solution
with respect to the $(t, x)$ variable.
\begin{prop}\label{prop4.3}
Under the hypothesis of Theorem \ref{theo1}, one has
\begin{equation}\label{4.1+19}
\Lambda^{1+\varepsilon}_{t, x}\,
(\varphi(t)\psi(x) f)\in H^5_l  (\RR^7),
\end{equation}
for any $l  \in\NN$ and some $\varepsilon>0$.
\end{prop}

\noindent{\bf Proof:}  Fix $s_0=\frac{s(1-s)}{(s+1)}$. Then
\eqref{4.2} and Proposition \ref{prop4.2} with $\lambda=s_0$ imply
$$
\Lambda^{s}_{v}\Lambda^{s_0}_{t, x} g\in H^{5}_l(\RR^7).
$$
It follows that,
\begin{equation*}\label{4.1+17}
 (\Lambda^{s_0}_{t, x} g)_t + v\,\cdot\,\partial_x {(\Lambda^{s_0}_{t, x} g)}
  =\Lambda^{s_0}_{t, x} Q(f_2,\,\, g)+ \Lambda^{s_0}_{t, x}  G\in
H^{-s}_l(\RR^7).
\end{equation*}
By applying Lemma \ref{lemm2.1} with $s'=s$, we can deduce that
\begin{equation*}\label{4.1+18}
\Lambda^{s_0+s_0}_{t, x}
(\varphi(t)\psi(x) f)\in H^5_l  (\RR^7),
\end{equation*}
for any $l  \in\NN$. If $2s_0<1$, by using Proposition \ref{prop4.2}
with $\lambda=2s_0$ and Lemma \ref{lemm2.1} with $s'=s$, we have
$$
\Lambda^{s}_{v} (\varphi(t)\psi(x) f),\,\,\Lambda^{2s_0}_{t, x}
(\varphi(t)\psi(x) f)\in H^5_l (\RR^7)\,\,\Rightarrow\,
\Lambda^{3s_0}_{t, x} (\varphi(t)\psi(x) f)\in H^5_l (\RR^7).
$$
Choose $k_0\in\NN$ such that
$$
k_0 s_0<1,\,\,\,\,\,\,\,\,\, (k_0+1) s_0=1+\varepsilon>1.
$$
Finally, (\ref{4.1+19}) follows from (\ref{4.1}) and Proposition
\ref{prop4.2} with $\lambda=k_0s_0$ by induction. And this completes
the proof of the proposition.

\vskip0.5cm
\subsection{Proof of Theorem \ref{theo1}}\label{s5}
\setcounter{equation}{0}
\smallskip
In this subsection, we give the proof of Theorem 1.1 with the above
preparations. The proof is also based on an induction argument.

\smallbreak From Propositions \ref{prop3.1} and 
\ref{prop4.3}, it follows that for any $l\in\NN$,
\begin{equation}\label{5.1}
\Lambda^{s}_{v}\, (\varphi(t)\psi(x) f),\,\,\,\,\,\,\, \nabla_{t, x}
\, (\varphi(t)\psi(x) f) \in H^5_l(\RR^7).
\end{equation}

These facts will be used to get the high order regularity
with respect to the variable $v$.

\begin{prop}\label{prop5.1}
Let $0<\lambda< 1$. Suppose that, for any cutoff functions $\varphi\in
C^\infty_0(]T_1, T_2[), \psi\in C^\infty_0(\Omega)$ and all $l\in\NN$,
\begin{equation}\label{5.2+0}
\Lambda^{\lambda}_{v}\,
(\varphi(t)\psi(x) f),\,\,\,\,\,
\nabla_{x}\, (\varphi(t)\psi(x) f)
\in H^5_l  (\RR^7).
\end{equation}
Then, for any cutoff function and any $l\in\NN$,
\begin{equation}\label{5.2}
\Lambda^{\lambda+s}_v\, (\varphi(t)\psi(x)
f)\in H^5_l(\RR^7).
\end{equation}
\end{prop}

\noindent {\bf Proof : } Recall that $g=\partial^\alpha
(\varphi(t)\psi(x) f)$ with $|\alpha|\leq 5$ and
$$
g_{N, l}=P_{N, l}\, g=\psi_1(x) S_N(D_x)\, W_{l  }\,\, S_N(D_v) g.
$$
Choose $ \Lambda^{2 \, \lambda}_v g_{N, l} $ as a test
function for equation (\ref{4.3}). Then, one has
\begin{equation}\label{5.2+}
\Big( \big[\Lambda^{\lambda}_v, \,\, v\big]\,\cdot\,\partial_x \,
g_{N, l} ,\,\,  \Lambda^{\lambda}_v g_{N, l}\Big)_{L^2(\RR^7)}
=\Big( \Lambda^{\lambda}_v \big\{Q({f}_2,\,\, g_{N, l})+G_{N,
l}\big\},\,\, \Lambda^{\lambda}_v g_{N, l}\Big)_{L^2(\RR^7)}.
\end{equation}
Since
$$
\big[\Lambda^{\lambda}_v, \,\, v\big]\,\cdot\,\partial_x =\lambda
\Lambda^{\lambda-2}_v \, \partial_v\,\cdot\,\partial_x,
$$
and $\Lambda^{\lambda-2}_v \, \partial_v$ are bounded operators in
$L^2$, for any $0<\lambda< 1$, we have, by using the hypothesis
(\ref{5.2+0}) that
\begin{equation}\label{5.3}
\Big|\Big( \big[\Lambda^{\, \lambda}_v, \,\,
v\big]\,\cdot\,\partial_x g_{N, l},\,\,\Lambda^{\, \lambda}_v g_{N,
l}\Big)_{L^2(\RR^7)}\Big| \leq C
\|\Lambda^{\lambda}_{v}\,g\|_{L^2_l(\RR^7)}\|\nabla_x
\,g\|_{L^2_l(\RR^7)} ,
\end{equation}
and when $1/2\leq s<1$.
\begin{equation}\label{5.4}
\Big|\Big(\Lambda^{\lambda}_v G_{N, l},\,\, \Lambda^{\lambda}_v
g_{N, l}\Big)_{L^2(\RR^7)}\Big| \leq C \|f_2\|_{H^5_7(\RR^7)}
\|\Lambda^{\lambda}_v g\|_{L^2_{l+\gamma^+ +2s}(\RR^7)}
\|\Lambda^{\lambda+2s-1+\delta}_v g_{N, l}\|_{L^2(\RR^7)}
\end{equation}
$$
\leq \varepsilon \|\Lambda^{s}_vW_{\gamma/2}\Lambda^{\lambda}_v g_{N,
l}\|^2_{L^2(\RR^7)}+ C_\varepsilon \|f_2\|^2_{H^5_7(\RR^7)}
\|\Lambda^{\lambda}_{v}\,g\|^{2k}_{L^2_{k'l}(\RR^7)}.
$$

By setting $M = \Lambda_v^\lambda$ in Proposition \ref{coro-comm-v},
we have
\begin{eqnarray}\label{5.5}
&&\Big|\Big(\Lambda^{\lambda}_v Q(\tilde{f},\,\, g_{N, l})-
Q(\tilde{f},\,\,\Lambda^{\lambda}_v g_{N, l}),\,\,
\Lambda^{\lambda}_v g_{N, l}\Big)_{L^2(\RR^7)}\Big|\\
&\leq& C \|f_2\|_{L^\infty(\RR^4_{t, x}; L^1_{\gamma_+}(\RR^3_v))}
\Big( \|\Lambda^{\lambda}_v g_{N, l}
\|^2_{L^2(\RR_{t,x}^4 ; L^2_{\gamma+}(\RR_v^3))}
+ \| g_{N, l}\|^2_{L^2(\RR^7)}\Big)
\|\Lambda^{\lambda}_vg_{N, l}\|^2_{L^2(\RR^7)} \nonumber\\
&\leq& C\|f_3\|_{H^5_7(\RR^7)}
\|\Lambda^{\lambda}_{v}\,g\|^2_{L^2_{l+1}(\RR^7)} , \nonumber
\end{eqnarray}
when $0<s<1/2$. Moreover when $1/2 \leq s <1$, we have
\begin{align}\label{5.5+0}
\Big|&\Big(\Lambda^{\lambda}_v Q(f_2,\,\, g_{N, l})-
Q(f_2,\,\,\Lambda^{\lambda}_v g_{N, l}),\,\,
\Lambda^{\lambda}_v g_{N, l}\Big)_{L^2(\RR^7)}\Big|
\\ \notag
&\leq C \|f_2\|_{L^\infty(\RR^4_{t, x};
L^1_{(2s + \gamma -1)_+}(\RR^3_v))}
\\&\hspace*{1cm}\times \nonumber
\Big( \|\Lambda^{\lambda}_v g_{N, l}\|^2_{L^2(\RR_{t,x}^4 ;
L^2_{(2s+\gamma-1)_+}(\RR_v^3))}
+ \| g_{N, l}\|^2_{L^2(\RR^7)}\Big)
\|\Lambda^{\lambda+2s-1+\delta}_vg_{N, l}\|^2_{L^2(\RR^7)} \nonumber
\\
&\leq  \varepsilon \|\Lambda^{s}_vW_{\gamma/2}\Lambda^{\lambda}_v
g_{N, l}\|^2_{L^2(\RR^7)}+ C_\varepsilon \|f_3\|^2_{
H^5_7(\RR^7)} \|\Lambda^{\lambda}_{v}\,g\|^{2k}_{L^2_{k'l}(\RR^7)}.
\nonumber
\end{align}

Now the coercivity estimate (\ref{E-sub-estimate+}) gives,
\begin{eqnarray}\label{5.6}
&&-\Big(Q(f_2,\, \Lambda^{\lambda}_v g_{N,
l}),\,\,\Lambda^{\lambda}_v g_{N, l}\Big)_{L^2(\RR^7)} \geq
C_0\|\Lambda^{s}_vW_{\gamma/2}\Lambda^{\lambda}_v g_{N,
l}\|^2_{L^2(\RR^7)}
\\
&&\hskip 4cm - C \| f_2\|_{L^\infty(\RR^4_{t, x};
L^1_{\max\{\gamma^+,\, 2-\gamma^+\}}(\RR^3_v))}\|\Lambda^{\lambda}_v
g_{N, l}\|^2_{L^2_{\gamma^+/2}(\RR^7)}. \nonumber
\end{eqnarray}
Thus, Proposition \ref{prop5.1} is proved by the following estimate
\begin{equation}\label{5.7}
\|\Lambda^{s}_vW_{\gamma/2}\Lambda^{\lambda}_v g_{N,
l}\|^2_{L^2(\RR^7)}\leq C \Big(\|f_3\|^2_{H^5_7(\RR^7)}+
\|\Lambda^{\lambda}_{v}\,g\|^{2k}_{L^2_{k'l}(\RR^7)}\Big),
\end{equation}
where $C$ is independent on $N$.\\


We can now conclude  the following regularity result with
respect to the variable $v$.
\begin{prop}\label{prop5.2}
Under the hypothesis of Theorem \ref{theo1}, one has
\begin{equation}\label{5.8}
\Lambda^{1+\varepsilon}_{v}
(\varphi(t)\psi(x) f)\in H^5_l  (\RR^7),
\end{equation}
for any $l  \in\NN$ and some $\varepsilon>0$.
\end{prop}
Again, this result follows by induction.  Indeed, notice that there
exists $k_0\in\NN$ such that
$$
k_0 s<1,\,\,\,\,\,\,\,\,\, (k_0+1) s=1+\varepsilon>1.
$$
Then we get (\ref{5.8}) from (\ref{3.3}), Proposition
\ref{prop5.1} with $\lambda=k_0s$ and (\ref{5.7}), by induction.

\bigbreak\noindent {\bf High order regularity by iterations}

\smallbreak From Proposition \ref{prop4.3} (more precisely
(\ref{4.1+19})) and Proposition \ref{prop5.2}, we can now deduce
that, for any $l\in\NN$, and any cutoff functions $\varphi(t)$
and $\psi(x)$, 
$$
\varphi(t)\psi(x) f\in H^6_l(\RR^7).
$$
The proof of Theorem \ref{theo1} is then completed by induction.

Indeed, if $f$ is a solution of Boltzmann equation satisfying the
assumptions of Theorem \ref{theo1}, then, when
$m\geq 5$,  we have
$$
f\in \cH^m_l(]T_1, T_2[\times\Omega\times\RR^3_v),\,\,\forall
l\in\NN\,\,\Longrightarrow \,\, f\in
\cH^{m+1}_l(]T_1, T_2[\times\Omega\times\RR^3_v),\,\,\forall l\in\NN .
$$
Thus, the full regularity of Theorem \ref{theo1} is obtained by
induction from $m=5$.

\bigbreak

\section{Existence and uniqueness of local solutions}\label{section1}
\smallskip

The local existence of solutions to the spatially inhomogeneous
Boltzmann equation without angular cutoff
is so far not well studied. The strategy of proving the existence
in this section
is to approximate the non-cutoff cross-section by a family of
cutoff cross-sections and approximate the Boltzmann equation
by a sequence of iterative linear equations. Then by proving the existence of solutions to these approximate linear equations and by obtaining a
uniform estimate on these solutions
 with respect to the cutoff parameter
in some suitable weighted Sobolev space, the compactness will
lead to the convergence of the approximate
solutions to the desired solution for the original problem.
One of the techniques used here is to introduce a 
transformation defined by the time dependent Maxwellian developed previously in
\cite{ukai-2}. The purpose of this transformation is to get
an extra gain of one order higher weight in the velocity variable at the
expense of the loss of the decay in the time dependent Maxwellian.
Moreover, the uniqueness of the solution is also proved in some
function space.

\subsection{Modified Cauchy Problem}\label{E-s0} \setcounter{equation}{0}
\smallskip

By taking $\kappa, \rho>0$, we
set, for $ 0\leq t\leq T_0=\rho/(2\kappa)$,
$$
\mu_\kappa(t)=\mu(t,v)=e^{-(\rho-\kappa t)(1+ |v|^2)},
$$
and
$$
f=\mu_\kappa(t) g, \,\,\, \quad \Gamma^t(g, g)=\mu_\kappa(t)^{-1}Q(\mu_\kappa(t) g,\,
\mu_\kappa(t) g).
$$
Then the Cauchy problem (\ref{1.1b}) is reduced to 
\begin{equation}\label{E-Cauchy-B}
\left\{\begin{array}{l}
g_t+v\cdot\nabla_x g\  +\kappa (1+ |v|^2) g=\Gamma^t(g, g),\\
g|_{t=0}=g_0.
\end{array}\right.
\end{equation}

Our existence theorem can be stated as follows

\begin{theo}\label{E-theo-0.2}
Assume that $0<s<1/2,\, \gamma+2s<1$ and $\kappa, \rho > 0$. Let $
g_0 \in H^k_l ({\mathbb R}^6 ),\, g_0\geq 0$ for some $l\geq 3$ and $k\geq
{4}$. Then there exists $T_* \in ]0, T_0]$ such that the Cauchy problem
(\ref{E-Cauchy-B}) admits a unique  non-negative solution
$$
g \in C^0 ([0,T_*];\,\, H^k_l({\mathbb R}^6))\bigcap\,\,
 L^2(]0, T_*[;H^k_{l+1}({\mathbb R}^6))\,. 
$$
\end{theo}

 We shall prove Theorem \ref{E-theo-0.2} by  cutoff approximations.
For simplicity of  notations, we will denote $\mu_\kappa(t)$ by $\mu(t)$
without any confusion.

Recall that the cross-section is of the form of $B(|v-v_*|, \cos
\theta)=\Phi(|v-v_*|) b(\cos\theta)$ which satisfies (\ref{E-soft-p})
and (\ref{1.2}). For $0<\varepsilon<\,<1$,
 we approximate (cutoff)  the cross-section by
$$
b_\varepsilon(\cos\theta)=\left\{\begin{array}{l} b(\cos\theta),
\,\,\,\,\,\mbox{if}\,\,\,|\theta|\geq 2\varepsilon,\\
b(\cos\varepsilon), \,\,\,\,\,\mbox{if}\,\,\,|\theta|\leq
2\varepsilon.
\end{array}\right.
$$
Denote by $\Gamma^t_\varepsilon(g,\,g)$ the collision operator
corresponding to the above cutoff cross-section
$B_\varepsilon=\Phi(v-v_*)b_\varepsilon(\cos\theta)$.

By using the collisional energy conservation,
$$
|v'_*|^2+|v'|^2=|v_*|^2+|v|^2,
$$
we have $\mu_*(t)=\mu^{-1}(t)\,\mu'_*(t)\,\mu'(t)$. Then for
some suitable functions $U, V$, it holds that
\begin{eqnarray}
\Gamma_\varepsilon^t(U,\,
V)(v)&=&\mu^{-1}(t,v)\iint_{\RR^3_{v_*}\times\mathbb S^{2}_{\sigma}}
B_\varepsilon(v-v_*,\, \sigma) \big(\mu'_*(t) U'_* \mu'(t) V'-
\mu_*(t) U_*
\mu(t) V\big) d v_* d \sigma\nonumber\\
&=&\iint_{\RR^3_{v_*}\times\mathbb S^{2}_{\sigma}}
B_\varepsilon(v-v_*,\, \sigma)\mu_*(t)\, \big( U'_*  V'-  U_* V\big)
d v_* d \sigma
={\cT}_\varepsilon(U,\, V,\, \mu(t))\label{E-3.105}\\
&=&Q_\varepsilon(\mu(t) U,\, V)+ \iint_{\RR^3_{v_*}\times\mathbb
S^{2}_{\sigma}} B_\varepsilon(v-v_*,\, \sigma) (\mu_*(t)
-\mu'_*(t))U'_* V' d v_* d \sigma.\nonumber
\end{eqnarray}
Then we have the following formula coming from the Leibniz formula
in the $x$ variable and the translation invariance property in the $v$
variable. For any $\alpha, \beta\in\NN^3$,
\begin{eqnarray}
&&\partial^\alpha_x\partial^\beta_v\Gamma_\varepsilon^t(U,\,\,
V)\nonumber\\
&=&\sum_{\alpha_1+\alpha_2=\alpha;\,
\beta_1+\beta_2+\beta_3=\beta} C_{\alpha_1, \alpha_2, \beta_1,
\beta_2,
\beta_3}{\cT}_\varepsilon(\partial^{\alpha_1}_x\partial^{\beta_1}_v
U,\,\,\partial^{\alpha_2}_x\partial^{\beta_2}_v V,\,\,
\partial^{\beta_3}_v\mu(t)
)\nonumber\\
&=&Q_\varepsilon(\mu(t) U,\, \partial^\alpha_x\partial^\beta_v V)+
\iint_{\RR^3_{v_*}\times\mathbb S^{2}_{\sigma}}
B_\varepsilon(v-v_*,\, \sigma) (\mu_*(t) -\mu'_*(t))U'_*
(\partial^\alpha_x\partial^\beta_v V)'
d v_* d \sigma\nonumber\\
&&\,\,\,\,\,\,\,+\sum_{|\alpha_2|+|\beta_2|\leq|\alpha+\beta|-1}
C_{\alpha_1, \alpha_2, \beta_1, \beta_2,
\beta_3}{\cT}_\varepsilon(\partial^{\alpha_1}_x\partial^{\beta_1}_v
U,\,\,\partial^{\alpha_2}_x\partial^{\beta_2}_v V,\,\,
\partial^{\beta_3}_v\mu(t) )\nonumber\\
&=&A_1+A_2+A_3 \,. \label{E-3.106}
\end{eqnarray}

Firstly, we give the following
 upper weighted estimate on the nonlinear collision operator
with cutoff.

\begin{lemm}\label{E-lemm3.111}
Let $\, \gamma\in\RR $.
Then for any $\varepsilon>0,\, k\geq 4,\, l\geq 0$,  there exists
$C>0$ depending on $\varepsilon,\, k,\, l$ such that for any $U,V$
belonging to $H^k_l(\RR^6)$
\begin{equation}\label{E-uper-estimate}
\|\Gamma_\varepsilon^t(U,\, V)\|_{H^k_l(\RR^6)}\leq C
\|U\|_{H^k_{l+\gamma^+}(\RR^6)} \|V\|_{H^k_{l+\gamma^+}(\RR^6)},
\quad 0\le t\le T_0=\frac{\rho}{2\kappa}.
\end{equation}
\end{lemm}
{\bf Proof.} 
To prove \eqref{E-uper-estimate}, put
\begin{align*}
&g_1=\partial^{\alpha_1}_x\partial^{\beta_1}_v U,\,\, \qquad
h_2=\partial^{\alpha_2}_x\partial_v^{\beta_2} V,\,\, \qquad
\mu_3(t)=
\partial^{\beta_3}_v\mu(t),
\\&
{\cT}_\varepsilon(g_1,h_2,\mu_3(t))={\cT}_\varepsilon^+
-{\cT}_\varepsilon^-.
\end{align*}
Throughout this section, the estimates
$$
\mu(t,v),\quad |\mu_3(t)|= |\partial^{\beta_3}_v\mu(t,v)|\le
C_{\rho,\, k} \,\,e^{-\rho \langle v\rangle ^2/4}, \qquad t\in[0,
T_0], \quad v\in\RR^3,
$$
will often be used.

Firstly, we compute ${\cT}_\varepsilon^+$ as follows.
\begin{align*}
|W_l {\cT}_\varepsilon^+|&\le C \iint \langle |v-v_*|\rangle^\gamma
|\mu_3(t,v_*)|\frac{W_l}{(W_l)'_*(W_l)'}|(W_lg_1)'_*|| (
W_lh_2)'|dv_*d\sigma
\\&
\le C \Big[ \iint
\Big|\mu_3(t,v_*)\frac{W_l}{(W_l)'_*(W_l)'}\Big|^2dv_*d\sigma
\Big]^{1/2} \Big[ \iint\langle v'-v'_*\rangle^{2\gamma}|(W_lg_1)'_*
( W_lh_2)'|^2dv_*d\sigma \Big]^{1/2}
\\&
\le C_{\varepsilon} \Big[ \iint|(W_{l+\gamma^+}g_1)'_* (
W_{l+\gamma^+} h_2)'|^2dv_*d\sigma\Big]^{1/2},
\end{align*}
where we have used $|v-v_*|=|v'-v'_*|$ and $
\frac{W_l}{(W_l)'_*(W_l)'}\leq 1.$ 
Since the change of variables
\begin{equation}\label{change}
(v,v_*,\sigma) \to (v',v'_*,\sigma'),\qquad \sigma'=(v-v_*)/|v-v_*|,
\end{equation}
has a unit Jacobian, 
we get
\begin{align*}
\|W_l {\cT}_\varepsilon^+\|_{L^2(\RR^6)}^2 &\le
C\iiiint|(W_{l+\gamma^+}g_1)'_*
 ( W_{l+\gamma^+} h_2)'|^2dv_*d\sigma dvdx
 \\
&=C\iiiint|(W_{l+\gamma^+}g_1)'_*
 ( W_{l+\gamma^+} h_2)'|^2dv_*'d\sigma' dv'dx
 \\
&
\le C\int\|(W_{l+\gamma^+}g_1)\|^2_{L^2(\RR^3_v)} \| (
W_{l+\gamma^+} h_2)\|^2_{L^2(\RR^3_v)}dx.
\end{align*}
If $|\a_1 + \b_1| \leq k /2$,  then we have
\begin{align*}
\|W_l{\cT}_\varepsilon^+\|_{L^2(\RR^6)}
&\le C \|(W_{l+\gamma^+}g_1)\|_{L^\infty(\RR_x^3 ; L^2(\RR^3_v))} \| (
W_{l+\gamma^+} h_2)\|_{L^2(\RR^6_{x,v})}\\
&\leq C
\|U\|_{H^k_{l+\gamma^+}(\RR^6)} \|V\|_{H^k_{l+\gamma^+}(\RR^6)},
\end{align*}
because of the  Sobolev embedding theorem and the fact $k/2 +3/2 < k$ when $k \ge 4$.
When $|\a_2 + \b_2| \leq k /2$, the proof is similar.
This completes  the proof of the lemma.


\subsection{Cutoff approximations}\label{E-s3} \setcounter{equation}{0}
\smallskip

\smallbreak

We now study the following Cauchy problem for the cutoff Boltzmann
equation 

\begin{equation}\label{E-Cauchy-cut-off}
\left\{\begin{array}{l}
g_t+v\cdot\nabla_x g +\kappa \la v \ra^2 g=\Gamma^t_\varepsilon(g, g),\\
g|_{t=0}=g_0\,,
\end{array}\right.
\end{equation}
for which we shall obtain uniform estimates in weighted
Sobolev spaces.

\bigskip
We first prove the existence of weak solutions to this 
cutoff
Boltzmann equation.

\begin{theo}\label{E-exist-cut-off}
Assume that $\gamma\le 1$. Let $k\ge  4,\, l \ge 0$, $\varepsilon>0$
and $D_0>0$. Then, there exists $T_\varepsilon \in ]0, T_0]$
such that for any  non-negative initial
data $g_0$ satisfying
$$
g_0\in H^k_{l}({\mathbb R}^6), \qquad \|g_0\|_{H^k_{l}({\mathbb
R}^6)}\le D_0,
$$
the Cauchy problem (\ref{E-Cauchy-cut-off}) admits a unique non-negative
solution $g^\varepsilon$ having the property
$$
 g^\varepsilon\in C^0(]0, T_\varepsilon[;\,\, H^k_l({\mathbb R}^6)) ,
\qquad \|g^\varepsilon\|_{L^\infty(]0, T_\varepsilon[;\,\,
H^k_l({\mathbb R}^6))}\le 2D_0.  
$$
Moreover, this solution enjoys a moment gain in the sense that 
\begin{equation}\label{momentgain}
g^\varepsilon\in L^2(]0,T_\varepsilon[; H^k_{l+1}(\RR^6)).
\end{equation}
\end{theo}
\begin{rema}\label{Tepsilon}
{\rm (1) } Notice that we do not assume $g_0\in H^k_{l+1}(\RR^6)$
and the gain of the moment will be  essentially used  below in the proof of 
uniform estimates to compensate the singularity in the cross-section.\\
{\rm (2) } The regularity of $g^\varepsilon $ with respect to $t$
variable follows directly from the equation \eqref{E-Cauchy-cut-off}.
\\
{\rm (3) } Fix $\gamma,\, k,\, l$ as in the theorem. Then
$T_\varepsilon$ is a function of $\varepsilon$ and $ D_0$. In the following, when
we need to emphasize this dependency, we shall write
\[
T = T_\varepsilon(D_0).
\]
{\rm (4) } If $\gamma \leq 0$, we may take $\kappa=0$. In this case,
we do not have the moment gain (\ref{momentgain}), which is anyway not needed.
\end{rema}

\noindent {\bf Proof of Theorem \ref{E-exist-cut-off}. } We prove
the existence of non-negative solutions by successive approximation
that preserves the  non-negativity, which is defined by using the usual
splitting of the collision operator \eqref{E-3.105} into the the
gain (+) and loss (-) terms,
\begin{align*}
\Gamma^{t,+}_\varepsilon(g,h)&=\iint_{\RR^3_{v_*}\times\mathbb
S^{2}_{\sigma}} B_\varepsilon(v-v_*,\, \sigma)\mu_*(t)\,  g'_*  h' d
v_* d \sigma,
\\
\Gamma^{t,-}_\varepsilon(g,h)& = hL_\varepsilon(g),
\\
L_\varepsilon(g)&=\iint_{\RR^3_{v_*}\times\mathbb S^{2}_{\sigma}}
B_\varepsilon(v-v_*,\, \sigma)\mu(t,v_*)\,  g_* d v_* d \sigma.
\end{align*}
Evidently, Lemma \ref{E-lemm3.111} applies to
$\Gamma^{t,\pm}_\varepsilon$, and in view of \eqref{E-soft-p}, the
linear operator $L_\varepsilon$ satisfies
\begin{equation}\label{L}
|\partial^{\alpha}_x\partial^{\beta}_vL_\varepsilon(g)(t,x,v)|\le
 C\langle v
\rangle^{\gamma-|\beta|}
\|\partial^{\alpha}_xg\|_{L^2(\R^3_v)},\quad t\in[0,T_0],
\end{equation}
for a constant $C>0$ depending on $\varepsilon$, because $|\mu(t,v_*)
\partial_v^\beta \la v - v_* \ra^\gamma| \leq C \la v \ra^{\gamma-
|\beta|}$.

We now define a sequence of approximate solutions
$\{g^n\}_{n\in\NN}$ by
\begin{equation}\label{iteration}
\left\{\begin{array}{l}
g^0=g_0\, ;\\
\partial_t g^{n+1}+v\cdot\nabla_x g^{n+1}+\kappa\langle |v|
\rangle^2 g^{n+1}=\Gamma_\varepsilon^{t,+} (g^n,\,
g^n)-\Gamma_\varepsilon^{t,-} (g^n,\, g^{n+1}),\qquad
 \\
g^{n+1}|_{t=0}=g_0.
\end{array}\right.
\end{equation}
Actually, in view of \eqref{L} we consider the mild form
\begin{align}\label{mild}
g^{n+1}(t,\, x,\, v)=&e^{-\kappa\langle |v|\rangle^2 t-V^n(t,\,
0)}g_0(x-tv,\, v)
\\ \notag
&+\int^t_0
 e^{-\kappa\langle |v|\rangle^2(t-s)-V^n(t,\, s)}
\Gamma_\varepsilon^{s,\, +} (g^n,\, g^n)(s, x-(t-s)v, \, v) ds,
\end{align}
where
\[
V^n(t,\, s)=\int_s^t L_{\varepsilon}(g^n)(s, x-(t-s)v,\, v)ds.
\]

First, we note from Lemma \ref{E-lemm3.111} that for any $T\in\,]0,
T_0]$, $T_0=\rho/(2\kappa) $, $g_0\ge0$, and
\[
g^n\in L^\infty(]0, T[;\,\, H^k_l(\RR^6)),\qquad g^n\ge 0,
\]
 the mild form \eqref{mild} determines
$g^{n+1}$ in the function class
\begin{equation}\label{loss}
g^{n+1}\in L^\infty(]0, T[;\,\, H^k_{l-\gamma^+}(\RR^6)),\qquad
g^{n+1}\ge 0,
\end{equation}
and solves \eqref{iteration}. Thus $g^{n+1}$ exists and is
non-negative, but appears to have a loss of weight in the velocity
variable. We shall now show that
the term $\kappa\langle v \rangle^2 g^{n+1}$ in \eqref{iteration}
not only recovers this weight loss but also creates a higher moment.
More precisely, we have the following lemma. Introduce the space and
norm by
$$
\begin{array}{ll}
&X=L^\infty(]0, T[;\,\, H^k_l(\RR^6))\cap L^2(]0, T[;\,\,
H^k_{l+1}(\RR^6)),
\\[0.3cm]&
|||g|||^2=\|g\|^2_{L^\infty(]0, T[;\,\, H^k_l(\RR^6))}
+\kappa\|g\|^2_{L^2(]0, T[;\,\, H^k_{l+1}(\RR^6))}\,.
\end{array}
$$
This norm depends on $k,l,T,\kappa$, but we omit this dependence in the
notation for
simplicity.

\begin{lemm}\label{estimate1}
Assume that $\gamma\le 1$ and let $k\ge  4, l \ge 0, \varepsilon>0$.
Then, there exist positive numbers $ C_1,C_2$ such that if
$\rho>0,\, \kappa >0$ and if
\begin{equation}\label{initial}
g_0\in H^k_l(\RR^6), \quad 
g^n\in L^\infty(]0, T[;\,\, H^k_l(\RR^6)), \quad 
\end{equation}
with some $T\le T_0$, the function $g^{n+1}$ given by \eqref{mild}
enjoys the properties
\begin{equation}\label{X}\begin{array}{l}
g^{n+1}\in X, 
\\
|||g^{n+1}|||^2\le e^{C_1 K_n T}
\left(\|g_0\|^2_{H^k_l(\RR^6)}+\frac{C_2}{\kappa} ||g^n||^4_{L^4(]0,
T[;\,\, H^k_{l}(\RR^6))}\right),
\end{array}
\end{equation}
where $K_n$ is a positive constant depending on
$\|g^n\|_{L^\infty(]0, T[;\,\, H^k_l({\mathbb R}^6))}$ and $\kappa$.
\end{lemm}

\noindent {\bf Proof.} 
 Put
\[
h^n =
h^n_{\alpha}=\partial^\alpha g^n.
\]
Differentiation of equation \eqref{iteration} yields
\begin{align*}
\partial_t h^{n+1}& +v\cdot\nabla_x h^{n+1}+
\kappa \langle  v \rangle ^2 h^{n+1} = G_1^+ -G_1^- +G_2 +G_3,
\\
&G_1^{+}=\partial^\alpha \Gamma_\varepsilon^{t,+}(g^n,\,\,
g^n), \qquad G_1^{-}=\partial^\alpha \Gamma_\varepsilon^{t,-}(g^n,\,\,
g^{n+1}), \\
&G_2=-[\partial^\a ,\,\, v\cdot\nabla_x] g^{n+1} ,
\\
&G_3=-
\kappa \sum_{|\tilde \beta|=1,2}C_{\tilde \beta}
 \partial_v^{\tilde \beta}\langle v \rangle ^2\partial^{\alpha - (0, \tilde \b)}
 g^{n+1}.
\end{align*}

Let
$\chi_j\in C^\infty_0(\RR^3), \,\, j\in\NN \,, $ be the cutoff function
\[
\chi_j(v)=\left\{\begin{array}{ll} 1\, ,\,\,\,\,\,\,& |v|\le j\,,\\
0\,, & |v|\ge j+1\, . \end{array}\right.
\]
We remark that
\eqref{loss} does not necessarily imply $W_{l+1}
h^{n+1}(t)\in L^2(\RR^6)$, but $\chi_jW_{l+1} h^{n+1}(t)$ $\in
L^2(\RR^6)$ for all $j\in \NN$.
Hence, we can use  $\chi_j^2W^2_l S_N^2(D_x) h^{n+1}$ as a test function to get
\begin{align}\label{diff1}
&\frac{1}{2}\frac{d}{dt}\|S_N(D_x)\chi_jW_l h^{n+1}\|^2+\kappa\|S_N(D_x)\chi_jW_{l+1}
h^{n+1}\|^2 \\ \notag
&\hskip2cm =(G_1^+ -G_1^-+G_2+G_3, S_N(D_x)^2 \chi_j^2W^2_l h^{n+1}).
\end{align}
Here and in what follows,
the norm
$\| \ \|$ and inner product $( \ , \ )$ are those of
$L^2(\RR^6_{x,v})$ unless otherwise stated.
We shall evaluate the inner products on the right hand side. Observe
that Lemma \ref{E-lemm3.111} gives, for $ t\in[0,T]$,
\begin{align*}
\Big|(G_1^+, S_N^2\chi_j^2W^2_l h^{n+1})\Big|&= \Big|(S_N\chi_jW_{l-1}G_1^+,
S_N\chi_jW_{l+1} h^{n+1})\Big| \le C \|W_{l-1}G_1^+\| \,
\|S_N \chi_jW_{l+1} h^{n+1}\|
\\&
\le C \|\Gamma^{t,+}_\varepsilon(g^n,g^n)\|_{H^k_{l-1}(\RR^6)}
 \, \|S_N\chi_jW_{l+1} h^{n+1}\|
\\&
\le C \|g^n\|^2_{H^k_{l}(\RR^6)}\| \, \|S_N\chi_jW_{l+1} h^{n+1}\|
\\&
\le \frac{C}{ \kappa}\|g^n\|^4_{H^k_{l}(\RR^6)}+\frac{\kappa}{4}
\|S_N\chi_jW_{l+1} h^{n+1}\|^2.
\end{align*}
On the other hand, Lemma \ref{E-lemm3.111} is not enough to evaluate
 $G_1^-$
because $G_1^-$ contains $g^{n+1}$ which is not known, at this
moment, to have moments required by Lemma \ref{E-lemm3.111}. However, this
obstacle is only superficial. Observe that
\[
G_1^- =\sum_{(\alpha_1, \beta_1) + \a_2 = \a}C_{\alpha_1,
\beta_1,\a_2 }
\Big(\partial^{\alpha_2}
g^{n+1}\Big)
\Big(\partial^{\beta_1}_vL(\partial^{\alpha_1}_xg^n)\Big).
\]
Define,
\[
H_{j,l}(g)
=\sum_{|\alpha|\le
k}\|\chi_jW_l\partial^\alpha  g\|^2,
\]
and write $H_{j,l}^n= H_{j,l}^n(t) = H_{j,l}(g_n(t))$.
By recalling \eqref{L}, we get
\begin{align*}
\Big|(G_1^-, S_N^2\chi_j^2W^2_l h^{n+1})\Big|& \le
\sum_{(\alpha_1, \beta_1) + \a_2 = \a}C_{\alpha_1,
\beta_1,\a_2 }
\|\chi_j \langle v\rangle^{\gamma-|\beta_1|}
W_{l-1}\partial^{\a_2}
g^{n+1}\| \, \|\partial^{\alpha_1}_xg^{n}\| \, \|S_N\chi_jW_{l+1}
h^{n+1}\|
\\&
\le C \|g^n\|_{H^k_{l}(\RR^6)}\| \, (H^{n+1}_{j,l})^{1/2} \,
\|S_N\chi_jW_{l+1} h^{n+1}\|
\\&
\le \frac{C'}{ \kappa}\|g^n\|^2_{H^k_{l}(\RR^6)}
H^{n+1}_{j,l}+\frac{\kappa}{4} \|S_N\chi_jW_{l+1} h^{n+1}\|^2.
\end{align*}
Here  $C, C'$ are positive constants independent of $\kappa$.

The estimate on the remaining two inner products are more straightforward
and can be given as follows.
\begin{align*}
\Big|(G_2+G_3,& \, S_N^2\chi_j^2W^2_l  h^{n+1})\Big| \le C  \|\chi_jW_{l-1}
(G_2+G_3)\| \,\|S_N\chi_jW_{l+1} h^{n+1}\|
\\&
\le C(\kappa+1)\Big(H^{n+1}_{j,l}\Big)^{1/2}\,\|S_N\chi_jW_{l+1}
h^{n+1}\|
\le C''\frac{(\kappa+1)^2}{\kappa} H^{n+1}_{j,l}+\frac{\kappa}{4} \|S_N\chi_j W_{l+1}
h^{n+1}\|^2.
\end{align*}
The constants $C,C''$ are independent of $\varepsilon$ and
$\kappa$.

Putting together all the estimates obtained above in \eqref{diff1}
yields
\[
\frac 12\frac{d}{dt}\|S_N\chi_jW_l
h^{n+1}\|^2+\frac{\kappa}{4}\|S_N\chi_jW_{l+1} h^{n+1}\|^2 \le
C'''\left\{\kappa+\frac{1}{ \kappa}(1+\|g^n\|^2_{H^k_{l}(\RR^6)})\right\}
H^{n+1}_{j,l} +\frac{C}{ \kappa}\|g^n\|^4_{H^k_{l}(\RR^6)}.
\]
Summing up estimates for $h^{n+1} = h^{n+1}_{\a}$ over $|\alpha|\le k$ then yields,
$$
\frac{d}{dt}H_{j,l}(S_Ng^{n+1})+\kappa H_{j,l+1}(S_N g^{n+1}) \le
C_1 K_n
H_{j,l}(g^{n+1}) +\frac{C_2}{ \kappa}\|g^n\|^4_{H^k_{l}(\RR^6)},
$$
where
\[
K_n=\kappa+\frac{1}{\kappa}\Big(\|g^n\|^2_{L^\infty(]0,T[;H^k_{l}(\RR^6))}+1 \Big),
\]
and
  $C_1>0$ is a constant independent of $\varepsilon, \kappa$
while $C_2$ is independent of $ \kappa$ but depends on
$\varepsilon$.
By integrating the above estimate over $[0,t]$ and taking the limit $N \rightarrow \infty$, we get
\begin{align*}
H^{n+1}_{j,l}(t)&+ \kappa\int_0^t H^{n+1}_{j,l+1}( \tau )d \tau
\\ &\notag
\le  H^{n+1}_{j,l}(0)+ C_1 K_n \int_0^t H^{n+1}_{j,l}( \tau)d\tau
 +
\frac{C_2}{ \kappa} \int_0^t \|g^n( \tau )\|^4_{H^k_{l}(\RR^6)}d\tau, \qquad t\in[0,T],
\end{align*}
which gives a Gronwall type inequality
\begin{align}\label{uniform1}
H^{n+1}_{j,l}(t)&+ \kappa\int_0^t e^{C_1 K_n (t-\tau )}H^{n+1}_{j,l+1}(\tau )d  \tau
\\ &\notag
\le  e^{C_1 K_n t}H^{n+1}_{j,l}(0)
 +
\frac{C_2}{ \kappa} \int_0^te^{C_1 K_n
(t-\tau)}\|g^n(\tau )\|^4_{H^k_{l}(\RR^6)}d\tau , \qquad t\in[0,T],
\end{align}
for all $j\in\NN$.
Since
\[
H^{n+1}_{j,l}(0)\le \|g_0\|_{H^k_{l}}^2,
\]
and $1\le e^{C_1 K_n (t-\tau )}\le e^{C_1 K_n t}$, \eqref{uniform1}
gives
\begin{align*}
H^{n+1}_{j,l}(t)+& \kappa\int_0^t H^{n+1}_{j,l+1}(\tau )d \tau
\le  e^{C_1 K_n t}\Big\{\|g_0\|_{H^k_{l}}^2
 +
\frac{C_2}{ \kappa} \int_0^t\|g^n(\tau )\|^4_{H^k_{l}(\RR^6)}d \tau\Big\},
\qquad t\in[0,T].
\end{align*}
Since the right hand side is independent of $j$, we see that
$\{\chi_j \pa^\a g^{n+1}\}_{j\in\NN}$, $|\a| \leq k$ is  weakly* compact in
$L^\infty(]0,T[;L^2_l(\RR^6))$ and weakly compact in
$L^2(]0,T[;L^2_{l+1}(\RR^6))$. Take a convergent subsequence.
Apparently, its limit is $ h^{n+1}(t)$. This is true for all
$|\alpha|\le k$ so that we can now conclude that
\[
g^{n+1}\in X\:=L^\infty(]0,T[;H^k_{l}(\RR^6))\cap
L^2(]0,T[;H^k_{l+1}(\RR^6)),
\]
and by Fatou's theorem,
\begin{align*}
|||g^{n+1}&|||^2\le \liminf_{j\to\infty}
\|H^{n+1}_{j,l}\|_{L^\infty(]0,T[)} +\kappa \liminf_{j\to\infty}
\|H^{n+1}_{j,l+1}\|_{L^1(]0,T[)}
\\&
 \le  e^{C_1K_n T}\Big(\|g_0\|_{H^k_{l}}^2
 +
\frac{C_2}{ \kappa} \|g^n\|^4_{L^4(]0,T[;H^k_{l}(\RR^6)}\Big).
\end{align*}
Now the proof of Lemma \ref{estimate1} is completed.\\

\smallbreak We are now ready to prove the convergence
of $\{g^n\}_{n\in\NN}$. Fix ${ \kappa>0} $, let $D_0, g_0$ be as in
Theorem \ref{E-exist-cut-off} and introduce an induction hypothesis
\begin{equation}\label{E-3.101}
\|g^{n}\|_{L^\infty(]0, T[;\,\, H^k_l(\RR^6))}\leq 2D_0.
\end{equation}
for some $T\in\,]0,T_0]$. Notice that the factor $2$ can be  any
number $>1$.

\eqref{E-3.101} is true for $n=0$ due to \eqref{initial}. Suppose
that this is true for some $n>0$. We shall determine $T$ independent
of $n$. A possible choice is given by
\begin{align}\label{Tchoice}
e^{C_1K_0 T}= 2, \quad \frac{2^4C_2}{\kappa}TD_0^2= 1 \qquad
\text{where\qquad}K_0=\kappa+\frac{2D_0+1}{\kappa },
\end{align}
or
\[
T=\min\left\{\frac{\log 2}{ C_1K_0}, \,
\frac{\kappa}{2^4C_2D_0^2}\right\}.
\]
In fact, \eqref{X} and \eqref{E-3.101} yield that $g^{n+1}\in X$ and
\begin{align*}
|||g^{n+1}|||^2&\le e^{C_1K_0 T}
\Big(\|g_0\|^2_{H^k_l(\RR^6)}+\frac{C_2}{\kappa}
T||g^n||^4_{L^\infty(]0, T[;\,\, H^k_{l}(\RR^6))}\Big)
\\&
\le e^{C_1K_0 T} \Big(D_0^2+\frac{C_2}{\kappa} T2^4D_0^4\Big)\le
4D_0^2.
\end{align*}
That is, the induction hypothesis \eqref{E-3.101} is fulfilled for
$n+1$, and hence  holds for all $n$.

For the convergence, set $w^n=g^{n}(t)-g^{n-1}(t)$, for which
\eqref{iteration} leads to
\begin{equation*}
\left\{\begin{array}{ll}
\partial_t w^{n+1}+v\cdot\nabla_x w^{n+1}+\kappa\langle |v|
\rangle^2 w^{n+1}=&\Gamma_\varepsilon^{t,+} (w^n,\,
g^n)+\Gamma_\varepsilon^{t,+} (g^{n-1},\, w^{n}),
\\&-
\Gamma_\varepsilon^{t,-} (w^n,\, g^{n+1})-\Gamma_\varepsilon^{t,-}
(g^{n-1},\, w^{n+1}),
 \\
w^{n+1}|_{t=0}=0.&
\end{array}\right.
\end{equation*}
By the same computation as used for \eqref{diff1}, but more
directly since we can now use test functions as
$S_N(D_x)^2 W_l^2 \partial^\a \, w^{n+1}$,  we get
\begin{align*}
|||w^{n+1}|||^2 \leq \frac 12 C_2e^{C_1K_0 T}\frac{1}{\kappa}T
 \Big\{\| g^{n+1}&\|^2_{L^\infty(]0, T[;\,\,H^k_l(\RR^6))}
 +\|g^n\|^2_{L^\infty(]0, T[;\,\,
H^k_l(\RR^6))}
\\&
+\|g^{n-1}\|^2_{L^\infty(]0, T[;\,\, H^k_l(\RR^6))}\Big\}
\|w^{n}\|^2_{L^\infty(]0, T[;\,\, H^k_l(\RR^6))},
\end{align*}
with the same constants $C_1, C_2$ and $K_0$ as above. Then, \eqref{E-3.101} and
\eqref{Tchoice} give
\[
|||g^{n+1}-g^n|||^2 \\
\le 2^4C_2D_0^2\kappa^{-1}T
 \|g^{n}-g^{n-1}\|^2_{L^\infty(]0, T[;\,\, H^k_l(\RR^6))}.
\]
Finally,  choose  $T$ smaller if necessary so that
\[
2^4C_2D_0^2\kappa^{-1}T \leq \frac 14.
\]
Then, we have proved that for any $n\geq 1,$
\begin{equation}\label{E-3.102}
|||g^{n+1}-g^n|||\leq \frac 12\,\, |||g^{n}-g^{n-1}|||.
\end{equation}
Consequently, $\{g^{n}\}$ is a convergence sequence in $X$, and the limit
\[
g^\varepsilon\in X,
\]
is therefore a non-negative solution of the Cauchy problem
(\ref{E-Cauchy-cut-off}). The estimate (\ref{E-3.102}) also implie 
the uniqueness of solutions.

By means of the mild form \eqref{mild}, it can be proved also that for each $n$,
\[
g^{n}\in C^0([0,T]; H^k_l(\RR^6))
\]
and hence so is the limit $g^\varepsilon$.
The non-negativity of $g^\varepsilon$ follows because $g^n \geq 0$.
Now the  proof of Theorem \ref{E-exist-cut-off} is completed.

\subsection{Uniform estimate}\label{s4}
\setcounter{equation}{0}
\smallskip

We now prove the existence of solutions for the Cauchy problem
(\ref{E-Cauchy-B}) by the convergence of approximation sequence
$\{g^\varepsilon\}$ as $\varepsilon\rightarrow0$. The first step is to
prove the uniform boundedness of this approximation sequence. Below, the constant $C$ are various constants independent
of $\varepsilon>0$.

\begin{theo}\label{E-uniform-estimate}
Assume that $0<s<1/2,\,\, \gamma +2s< 1$. Let $g_0\in H^k_l(\RR^6),
g_0\geq 0$ for some $k\geq 4,\,\, l\geq 3$. Then there exists $T_* \in ]0,T_0]$
depending only on $\|g_0\|_{H^k_l}$ and independent of $\varepsilon$
satisfying the following property:
 If for some $0<T\leq T_0$\, ,
\begin{equation}\label{assume}
g^\varepsilon\in C^0(]0, T];\,\, H^k_l(\RR^6))\cap L^2(]0, T[;\,\,
H^k_{l+1}(\RR^6)),
\end{equation}
is a non-negative solution of the Cauchy problem
(\ref{E-Cauchy-cut-off}) and if 
$T_{**}=\min\{T,\,\,T_*\}$,
then it holds that
\begin{equation}\label{E-3.103}
\|g^\varepsilon\|_{L^\infty(]0, T_{**}[;\,\, H^k_l(\RR^6))}\leq 2
\|g_0\|_{H^k_l(\RR^6)}.
\end{equation}
 \end{theo}
\begin{rema}
The case $T_*\le T$ gives a uniform estimate 
of local solutions on 
the fixed time interval $[0,T_*]$
while the case $T<T_*$ gives an a priori estimate on the 
existence time interval $[0,T]$ of local solutions. The latter is used for the continuation argument of
local solutions, in Subsection 4.4 below.
\end{rema}
In the following, $\rho>0,\,\kappa>0$ are fixed. Furthermore, recall
$T_0=\rho/(2\kappa)$. We start with  a solution $g^\varepsilon$ subject
to \eqref{assume} for some  $T\in\,]0,T_0]$. For $\alpha \in\NN^6,
|\alpha| \leq k$, the differentiation of the equation
(\ref{E-Cauchy-cut-off}) implies
\begin{equation}\label{E-equ-1}
\partial_t (\partial^\alpha g^\varepsilon) +v\,\cdot\,\nabla_x (\partial^\alpha g^\varepsilon)
+\kappa \la v \ra^2  (\partial^\alpha g^\varepsilon)
=\partial^\alpha\Gamma_\varepsilon^t(g^\varepsilon,\,\,
g^\varepsilon) -[\partial^\alpha,\,\, v\cdot\nabla_x] g^\varepsilon
-\kappa [\partial^\alpha,\, \langle v\rangle ^2] g^\varepsilon.
\end{equation}

Since  $\partial^\alpha g^\varepsilon$ only belongs to $L^2_l$, now
as in Section \ref{section3}, we take,
$$
P^\star_{N,\, l}P_{N,\, l} (\partial^\alpha g^\varepsilon)
$$
as a test function in (\ref{E-equ-1}), where $l\geq 3$ and  $P_{N,\,
l}= S_N(D_x)S_N(D_v) \, W_l$ (we do not need the cutoff functions
$\varphi, \psi$ here). Then we have
\begin{align}\label{E-equ-2}
\frac 1 2 \frac{d}{dt}\|P_{N,\, l}(\partial^\alpha
g^\varepsilon)(t)\|^2_{L^2(\RR^6)} &
 + \kappa
\| W_1\,P_{N,\, l}(\partial^\alpha g^\varepsilon)(t)\|^2_{L^2(\RR^6)}
\\ \notag
& + \kappa \Big ([S_N(D_v),\, \la v \ra^2] W_l\,(\partial^\alpha
g^\varepsilon), \,\,S_N(D_x) P_{N,\, l} (\partial^\alpha
g^\varepsilon) \Big )_{L^2(\RR^6)}\\
= &\left(A_1+A_2+A_3+A_4 +A_5,\,\, P^\star_{N,\, l}P_{N,\,
l}(\partial^\alpha g^\varepsilon)\right)_{L^2(\RR^6)},\notag
\end{align}
where $A_1,\, A_2,\, A_3$ are defined in (\ref{E-3.106}) with
$U=V=g$ and
$$
 A_4=-[\partial^\a,\,\,
 v\cdot\nabla_x]g^\varepsilon,\qquad
 A_5= -\kappa \sum_{|\tilde \beta|=1,2}C_{\tilde \beta}
 \partial_v^{\tilde \beta}\langle v \rangle ^2\partial^{\alpha - (0, \tilde \b)} g^\varepsilon.
$$
We have firstly,
\begin{equation}\label{E-3.105+1}
\left|\left(A_4,\,\, P^\star_{N,\, l}P_{N,\, l}(\partial^\alpha
g^\varepsilon)\right)_{L^2(\RR^6)}\right| \leq
 C\|g^\varepsilon(t)\|^2_{H^k_l(\RR^6)},
\end{equation}
and
\begin{equation} \label{E-3.105+1-2}
\left|\left(A_5,\,\,P^\star_{N,\, l}P_{N,\, l}(\partial^\alpha
g^\varepsilon)\right)_{L^2(\RR^6)}\right| \le C\kappa
\|g^\varepsilon(t)\|^2_{H^k_{l}(\RR^6)}+ \frac{\kappa}{
4}\|g^\varepsilon(t)\|^2_{H^k_{l+1}(\RR^6)}.
\end{equation}
We also have
\begin{align}\label{additional-term}
\left|\kappa \Big ([S_N(D_v),\, \la v \ra^2] W_l\, (\partial^\alpha
g^\varepsilon), \,\,S_N(D_x) P_{N,\, l} (\partial^\alpha
g^\varepsilon) \Big )_{L^2(\RR^6)}\right| \\  \notag \leq C\kappa
\|g^\varepsilon(t)\|^2_{H^k_{l}(\RR^6)}+ \frac{\kappa}{ 4}\|
g^\varepsilon(t)\|^2_{H^k_{l+1}(\RR^6)}.
\end{align}
We now study the term $A_1$  by using the non-negativity of
$g^\varepsilon$ and the coercivity of collision operators.

\begin{prop}\label{E-prop4.1}
Assume that  $0<s<1/2,\, \gamma\in\RR $. There exists $C>0$
independent of $\varepsilon$ such that for any $\alpha \in\NN^6,
|\alpha|\leq k$, $ k\geq  4, \, l\geq 3$,
\begin{equation}\label{E-3.105+4}
\left(A_1,\,\, P^\star_{N,\, l}P_{N,\, l}(\partial^\alpha
g^\varepsilon)\right)_{L^2(\RR^6)}\leq C
\|g^\varepsilon(t)\|^2_{H^k_l(\RR^6)}
\|g^\varepsilon(t)\|_{H^k_{l+\gamma^+}(\RR^6)},
\end{equation}
for any $0\leq t\leq T\leq T_0$.
\end{prop}

\noindent {\bf Proof : } By setting $h= \partial^\alpha
g^\varepsilon$, we have,
\begin{eqnarray*}
&&\hskip-1cm \Big (A_1,\,\, P^\star_{N,\, l}P_{N,\, l} h\Big
)_{L^2(\RR^6)} =\Big( P_{N,\, l} Q_\varepsilon(\mu
g^\varepsilon,\,h),\,\, (P_{N,\, l} h)\Big)_{L^2(\RR^6)}\\
&&=\Big(Q_\varepsilon\big(\mu(t) g^\varepsilon,\,(P_{N,\, l}
h)\big),\,\,
(P_{N,\, l} h)\Big)_{L^2(\RR^6)} \\
&&+ \Big( P_{N,\, l} Q_\varepsilon\big(\mu(t)
g^\varepsilon,\,h\big)-Q_\varepsilon\big(\mu(t)
g^\varepsilon,\,(P_{N,\, l} h)\big),\,\,
(P_{N,\, l} h)\Big)_{L^2(\RR^6)}\\
&&=B_1+B_2.
\end{eqnarray*}
Since $\mu(t)\,g^\varepsilon(t, \, x,\, v)\geq 0$, we have, in the
same way as  Theorem \ref{E-theo-0.1} with the cancellation lemma,
\begin{eqnarray*}
B_1&=&-\frac 1 2 \iiiint_{\RR^3_{x}\times\RR^3_{v}
\times\RR^3_{v_*}\times\SS^2_{\sigma}} B_\varepsilon(v-v_*,\,
\sigma)\,(\mu(t)\,g^\varepsilon)_*\, \Big((P_{N,\, l} h)'-
(P_{N,\, l} h)\Big)^2\, d v_* d \sigma dv dx\\
&&+ \frac 1 2 \iiiint_{\RR^3_{x}\times\RR^3_{v}
\times\RR^3_{v_*}\times\SS^2_{\sigma}} B_\varepsilon(v-v_*,\,
\sigma)\,(\mu(t)\,g^\varepsilon)_*\, \Big\{\Big((P_{N,\,
l}h)'\Big)^2- \big(P_{N,\, l} h\big)^2\Big\}\, d v_* d \sigma dv
dx\\
& \leq& \frac 1 2 \iiiint_{\RR^3_{x}\times\RR^3_{v}
\times\RR^3_{v_*}\times\SS^2_{\sigma}} B_\varepsilon(v-v_*,\,
\sigma)\,(\mu(t)\,g^\varepsilon)_*\, \Big\{\Big((P_{N,\, l}
h)'\Big)^2-
\big(P_{N,\, l} h\big)^2\Big\}\, d v_* d \sigma dv dx\\
 &\leq &
 C\iiint_{\RR^3_{x}\times\RR^3_{v_*}\times\RR^3_{v}} (\mu(t)
\,g^\varepsilon)_*\, \langle v-v_* \rangle^{\gamma^+} (P_{N,\, l}
h)^2dvdv_*dx
\\
& \le &C \|\mu W_{\gamma^+}  g^\varepsilon(t) \|_{L^\infty(
\RR^3_x;L^1(\RR^3_v))}  \|W_{l} h(t) \|_{L^2(\RR^6_{x, \, v}) }
\|W_{l+\gamma^+} h(t) \|_{L^2(\RR^6_{x, \, v}) }
 \\&
\leq& C  \|g^\varepsilon(t)\|_{H^{3/2+\delta}(\RR^6_{x, \, v})}
\|g^\varepsilon(t)\|_{H^k_{l}(\RR^6_{x, \, v})}
\|g^\varepsilon(t)\|_{H^k_{l+\gamma^+}(\RR^6_{x, \, v})}, \quad
t\in[0,\, T],
\end{eqnarray*}
where we used the fact that $ b_\varepsilon(\cos\theta)\leq
b(\cos\theta)$.

By putting $S_N = S_N(D_x) $ $ \tilde S_N = S_N(D_v)$, we decompose
\begin{align*}
B_2&=\Big( S_N \tilde S_N \Big \{ W_lQ_\varepsilon\big(\mu(t)
g^\varepsilon,\,h\big)-Q_\varepsilon\big(\mu(t)
g^\varepsilon,\,(W_l h)\big) \Big \} ,\,\,
(P_{N, \, l}h )\Big)_{L^2(\RR^6)} \\
&+
\Big(S_N \Big \{ \tilde S_N Q_\varepsilon\big(\mu(t)
g^\varepsilon,\, (W_l h) \big)-Q_\varepsilon\big(\mu(t) g^\varepsilon,\,
\tilde S_N (W_l h)\big) \Big \},\,\,
(P_{N,\,l} h)\Big)_{L^2(\RR^6)}\\
&+\Big(S_NQ_\varepsilon\big(\mu(t)
g^\varepsilon,\,( \tilde S_N W_l h )\big)-Q_\varepsilon\big(\mu(t)
g^\varepsilon,\, S_N (  \tilde S_N W_l h)\big),\,\,  (P_{N,\,l} h) \, \Big)_{L^2(\RR^6)}
\\
&=B_{21}+B_{22}+B_{23}.
\end{align*}
By Lemma \ref{lemm3.4}, we get
\begin{align*}
|B_{21}|&= \left |\Big( \Big \{W_lQ_\varepsilon\big(\mu(t)
g^\varepsilon,\,h\big)-Q_\varepsilon\big(\mu(t)
g^\varepsilon,\,(W_l h)\big) \Big \},\,\, (\tilde S_N S_N P_{N,\,l}
h)\Big)_{L^2(\RR^6)} \right |
\\
& \leq C \|\mu
(t)g^\varepsilon(t)\|_{L^\infty(\RR^3_x;L^1_{l+\gamma^+}(\RR^3_v))}
\int_{\RR^3_x}\|W_{l+\gamma^+} h\|_{L^2(\RR^3_v)}\| P_{N,\,l}
h\|_{L^2(\RR^3_v)}dx
\\&
\leq C \|g^\varepsilon(t)\|_{L^\infty(\RR^3_x;L^2(\RR^3_v))}
\|g^\varepsilon(t)\|_{H^k_{l}(\RR^6)}
\|g^\varepsilon(t)\|_{H^k_{l+\gamma^+}(\RR^6)},
\\&
\le C\|g^\varepsilon(t)\|^2_{H^k_{l}(\RR^6)}
\|g^\varepsilon(t)\|_{H^k_{l+\gamma^+}(\RR^6)}, \quad t\in[0,\, T].
\end{align*}
It follows from Lemma \ref{lemm3.1} that
\begin{align*}
|B_{22}|&\le \left(\int_{\RR^3_x} \|\tilde S_N Q_\varepsilon(\mu(t) g^\varepsilon,
(W_l h) )-Q_\varepsilon(\mu(t) g^\varepsilon,\, \tilde S_N
(W_l h) )\|_{L^2(\RR^3_v)}^2 dx \right)^{1/2} \, \|P_{N,\,l} h \|_{L^2(\RR^6)}
\\&
\leq C \|\mu
(t)g^\varepsilon(t)\|_{L^\infty(\RR^3_x;L^1_{\gamma^+}(\RR^3_v))}
\|W_{l+\gamma^+} h\|_{L^2(\RR^6_v)}\|g^\varepsilon(t)\|_{H^k_{l}(\RR^6)} \\
& \le C\|g^\varepsilon(t)\|^2_{H^k_{l}(\RR^6)}
\|g^\varepsilon(t)\|_{H^k_{l+\gamma^+}(\RR^6)}, \quad t\in[0,\, T].
\end{align*}
Lemma \ref{lemm3.3} with $m =2s$ yields
\begin{align*}
|B_{23}|&\le C
 \|S_N Q(\mu(t)g^\varepsilon,\,\,  (\tilde S_N W_lh) \, )-
Q(\mu(t)g^\varepsilon,\,\, S_N\, ( \tilde S_N W_l h) \, )\|_{L^2(\RR^3_{x},\,
L^2(\RR^3_v))} \|P_{N,\,l} h \|_{L^2(\RR^6) }
\\
& \leq C \|\mu(t)\nabla_x  g^\varepsilon\|_{L^\infty(\RR^3_{ x},\,\,
L^1_{(2s+\gamma)^+}(\RR^3_v))}\| ( 2^{-N} \tilde S_N (W_l h) \|_{L^2(\RR^3_{
x},\,\, H^{2s}_{(2s+\gamma)^+}(\RR^3_v))}\|P_{N,\,l} h \|_{L^2(\RR^6) }
\\
&\le C\|g^\varepsilon(t)\|^2_{H^k_{l}(\RR^6)}
\|g^\varepsilon(t)\|_{H^k_{l+\gamma^+}(\RR^6)}, \quad t\in[0,\, T].
\end{align*}
Combining the above estimates proves Proposition \ref{E-prop4.1}.\\

For the term $A_2$ and $A_3$, we prove the following proposition.

\begin{prop}\label{E-prop4.2}
Assume that $0<s<1/2,\, \gamma+2s< 1$. Then, for any $\delta>0$,
there exists $C>0$ independent of $\,\,\varepsilon>0$ such that for
any $\alpha \in\NN^6, |\alpha|\leq k$, $k\geq 4, \, l\geq 3$,
\begin{equation}\label{E-3.105+5}
\left(A_2+A_3,\,\, P^\star_{N,\, l}P_{N,\, l} (\partial^\alpha
g^\varepsilon)\right) _{L^2(\RR^6)} \leq C
\|g^\varepsilon(t)\|^2_{H^k_l(\RR^6)}
\|g^\varepsilon(t)\|_{H^k_{l+(\gamma+2s+\delta)^+}(\RR^6)},
\end{equation}
for $t\in[0,\, T]$.
\end{prop}

\noindent {\bf Proof. } By putting $h =\partial^\alpha g^\varepsilon$
and $\tilde h = W_{l}^{-2} P^\star_{N,\, l}P_{N,\, l}
(\partial^\alpha g^\varepsilon) $ , we get
\begin{eqnarray*}
&&\left|\left(A_2,\,\, W^2_{l} \tilde
h\right)_{L^2(\RR^6)}\right|=\left|\iiiint_{\RR^3_{x}\times\RR^3_{v}
\times\RR^3_{v_*}\times\SS^2_{\sigma}} B_\varepsilon(v-v_*,\,
\sigma) (\mu_*(t) -\mu'_*(t))(g^\varepsilon)'_*  h'
(W^2_{l} \tilde h) d v_* d \sigma dv dx\right|\\
&&\leq \iiiint_{\RR^3_{x}\times\RR^3_{v}
\times\RR^3_{v_*}\times\SS^2_{\sigma}} B(v-v_*,\, \sigma) |\mu_*(t)
-\mu'_*(t)|\,|(g^\varepsilon)'_*|\, |(W_l h)'
(W_{l} \tilde h) |d v_* d \sigma dv dx\\
&&+\iiiint_{\RR^3_{x}\times\RR^3_{v}
\times\RR^3_{v_*}\times\SS^2_{\sigma}} B(v-v_*,\, \sigma) |(\mu_*(t)
-\mu'_*(t))|\, |(g^\varepsilon)'_*|\, \big|W_l-W'_l\big|\,\,| h'
 (W_{l} \tilde h)| d v_* d \sigma dv dx \\
&&= I_1+I_2.
\end{eqnarray*}
To estimate $I_1$, we notice that
\begin{equation}\label{mudiff}
|\mu(t,v_*)-\mu(t,v'_*)| \le C|v_*-v_*'|^\lambda\le C\theta^\lambda
|v-v_*|^\lambda \le C\theta^\lambda |v'-v'_*|^\lambda,\quad
\lambda\in[0,1], \ t\in[0,T_0],
\end{equation}
which is elementary for $\lambda=0,1$ and is obtained for general
$\lambda\in(0,1)$ by interpolation. Since $\gamma+2s<1$ is assumed
in the proposition, there is $\lambda\in(0,1)$ such that
$\lambda>2s, \gamma+\lambda\le 1$. By the manipulation on the primed
and non-primed variables ( see \eqref{change} ) we have
\begin{align*}
 I_1 &\le C\iiiint_{\RR^3_{x}\times\RR^3_{v}
 \times\RR^3_{v_*}\times\SS^2_{\sigma}} \langle v'-v'_* \rangle^{(\gamma+\lambda)^+}
 \theta^{-2-2s+\lambda}|(g^\varepsilon)'_*|\,\, | (W_l h)'|\,\, | (W_{l} \tilde h)| d v_* d
 \sigma dv dx
\\
 &\le C\iiint_{\RR^3_{x}
 \times\RR^3_{v_*}\times\SS^2_{\sigma}}
 \theta^{-2-2s+\lambda}|(W_{(\gamma+\lambda)^+}g^\varepsilon)_*|
\Big\{\int_{\RR^3_{v}}| (W_{l+(\gamma+\lambda)^+} h)(W_{l} \tilde
h)'|dv\Bigr\} d v_* d
 \sigma  dx
\\
& \leq C \|g^\varepsilon (t)\|_{L^\infty(\RR^3_x;
L^1_{(\gamma+\lambda)^+}(\RR^3_v))}\|g^\varepsilon (t)\|_{
H^k_{l}(\RR^6)} \|g^\varepsilon (t)\|_{
H^k_{l+(\gamma+\lambda)^+}(\RR^6)}
\\&
\le C\|g^\varepsilon (t)\|_{ H^k_{l}(\RR^6)}^2 \|g^\varepsilon
(t)\|_{ H^k_{l+(\gamma+\lambda)^+}(\RR^6)},
\end{align*}
for $l>(\gamma+\lambda)^++3/2$. In the third inequality we have used
again the fact that the Jacobian of changing of variable
$v\,\rightarrow\,v'$ is bounded.

Using \eqref{E-3.105+3} gives
\begin{align*}
I_2 &\le C \iiiint_{\RR^3_{x}\times\RR^3_{v}
\times\RR^3_{v_*}\times\SS^2_{\sigma}} \langle v'-v'_*
\rangle^{\gamma} \theta^{-1-2s}
(\mu_*(t)+\mu'_*(t))|(g^\varepsilon)'_* \big(W_l'W'_{l,*}\big) h'
  (W_{l} \tilde h) |d v_* d \sigma dv dx
  \\&
  =C(J_1+J_2).
\end{align*}
By the Schwarz inequality and  the Sobolev inclusion, 
we have
\begin{align*}
J_1&\le C \iiiint_{\RR^3_{x}\times\RR^3_{v}
\times\RR^3_{v_*}\times\SS^2_{\sigma}} \theta^{-1-2s}
\mu_*(t)|(W_{l+\gamma^+}g^\varepsilon)'_* (W_{l+\gamma^+} h)'
  (W_{l} \tilde h)| d v_* d \sigma dv dx \\
&\le  C \int_{\RR^3_x}  
  \Big( \iiint_{\RR^3_{v}
\times\RR^3_{v_*}\times\SS^2_{\sigma}}
  \theta^{-1-2s}
\mu_*(t)^2 
  |(W_{l} \tilde h)|^2 dv d v_* d \sigma  \Big)^{1/2}\\
  & \hskip1cm  \times 
  \Big( \iiint_{\RR^3_{v}
\times\RR^3_{v_*}\times\SS^2_{\sigma}}
  \theta^{-1-2s}
|(W_{l+\gamma^+}g^\varepsilon)'_* (W_{l+\gamma^+} h)'|^2 dv d v_*  d \sigma  \Big)^{1/2}dx
  \\&
\le  C \|\mu\|_{L^2(\RR_v^3)}
\int_{\RR^3_x}  
\|W_{l} \tilde h(x)\|_{L^2(\RR^3_v)}  
\|W_{l+\gamma^+}g^\varepsilon(x)\|_{L^2(\RR_v^3)}
\|W_{l+\gamma^+}h(x) \|_{L^2(\RR^3_v)}
dx
\\ & \le 
C\|g^\varepsilon\|_{L^\infty(\RR_x^3; L^2(\RR_v^3)}
   \|W_{l+\gamma^+}h\|_{L^2(\RR^6)}\|W_{l} \tilde h\|_{L^2(\RR^6)}
   \\
   &
\le C\|g^\varepsilon\|^2_{H^k_{l+\gamma^+}(\RR^6)}
   \|g^\varepsilon\|_{H^k_{l}(\RR^6)}.
\end{align*}
On the other hand, again by the  manipulation on the primed
and non-primed variables,
\begin{align*}
J_2&\le C \iiiint_{\RR^3_{x}\times\RR^3_{v}
\times\RR^3_{v_*}\times\SS^2_{\sigma}} \theta^{-1-2s}   | (\mu(t)
W_{l+\gamma^+}g^\varepsilon)'_* (W_{l+\gamma^+} h)'
  (W_{l} \tilde h)| d v_* d \sigma dv dx
  \\&
\le  C \iiiint_{\RR^3_{x} \times\RR^3_{v_*}\times\SS^2_{\sigma}}
\theta^{-1-2s}    |(\mu(t) W_{l+\gamma^+}g^\varepsilon)_*|
\Big\{\int_{\RR^3_v}|W_{l+\gamma^+} h| |  (W_{l} \tilde h)'|dv\Big\}
d v_* d \sigma dx
  \\&
  \le
    C\|\mu(t) W_{l+\gamma^+}g^\varepsilon\|_{L^\infty(\RR^3_x;L^1(\RR^3_v))}
   \|W_{l+\gamma^+}h\|_{L^2(\RR^6)}\|W_{l} \tilde h\|_{L^2(\RR^6)}
   \\&
   \le C\|g^\varepsilon\|^2_{H^k_{l}(\RR^6)}
   \|g^\varepsilon\|_{H^k_{l+\gamma^+}(\RR^6)}.
\end{align*}
Here, we have used $W_{l+\gamma^+}\mu^{1/2}(t)\le C$.

We consider now the term $A_3$.  For any $\alpha \in\NN^6, |\alpha|
\leq k$, $k\geq 4, \, l\geq 3$, denote
$$
h_1=\partial^{\alpha_1} g^\varepsilon,\,\,\, h_2=\partial^{\alpha_2}
g^\varepsilon,\,\,
$$
where
$$
\alpha_1+ \alpha_2 \leq \alpha,\,\,\,\,  \alpha_2 <\alpha \,.
$$
We shall compute
\begin{align*}
\left({\cT}_\varepsilon(h_1,\,  h_2,\,\, \tilde{\mu}),\,\, W^2_{l}
\tilde h\right)_{L^2(\RR^6)}
 =&\iiiint_{\RR^3_{x}\times\RR^3_{v}
\times\RR^3_{v_*}\times\SS^2_{\sigma}} B_\varepsilon( \tilde{\mu}_*(t)
-\tilde{\mu}'_*(t))(h_1)'_*  h_2' (W^2_{l} \tilde h) d v_* d \sigma dv dx
\\
&+ \iiiint_{\RR^3_{x}\times\RR^3_{v}
\times\RR^3_{v_*}\times\SS^2_{\sigma}} B_\varepsilon \, ( \tilde \mu
h_1)'_* \, \big (W_l-W'_l\big )\,\, h_2 '
 (W_{l} \tilde h)  d v_* d \sigma dv dx \\
&+ \left( {Q}_\varepsilon( \tilde{\mu} h_1,\, (W_l h_2) \,
),W_{l}\tilde h\right)_{L^2(\RR^6)} .
\end{align*}
For the last term,
by Theorem \ref{theo2.1} 
with $m=2s<1$, there exists $C>0$ independent of $\varepsilon$ such
that for $|\alpha_2|\leq|\alpha|-1$
and $\delta>0$,
\begin{eqnarray*}
&&\|{Q}_\varepsilon( \tilde \mu h_1,\, (W_l h_2)\,
)\|^2_{L^2(\RR^6)} \leq C \int_{\RR^3} \| \tilde \mu h_1(t,\, x,\,
\cdot\,)\|^2_{L^1_{(\gamma+2s)^+}(\RR^3_v)} \|(W_l h_2)(t,\,
x,\,\cdot\,)\|^2_{H^{2s}_{(\gamma+2s)^+}(\RR^3_v)}dx\\
&\leq&\left\{
\begin{array}{ll}C \|\tilde \mu h_1(t)\|^2_{L^\infty(\RR^3_x,\,\,
L^2_{3/2 +(\gamma+2s)^++\delta} (\RR^3_v))} \|(W_l
h_2)(t)\|^2_{L^2(\RR^3_x;\,\,
H^{2s}_{(\gamma+2s)^+}(\RR^3_v))}\,,\,\,\,\,\,\, &|\alpha_1|\leq 2,\\
C\|\tilde \mu h_1(t)\|^2_{L^2(\RR^3_x;\,\, L^2_{3/2
+(\gamma+2s)^++\delta}(\RR^3_v))} \|(W_l
h_2)(t)\|^2_{L^\infty(\RR^3_x,\,\,
H^{2s}_{(\gamma+2s)^+}(\RR^3_v))}\,,\,\,\,\,\,\, &|\alpha_1|>2,
\end{array}\right.\nonumber
\\
&\leq&\left\{
\begin{array}{ll}C \|h_1(t)\|^2_{H^{3/2+\delta}(\RR^3_x;\,\,
L^2(\RR^3_v))} \|(W_l h_2)(t)\|^2_{L^2(\RR^3_x;\,\,
H^{2s}_{(\gamma+2s)^+}(\RR^3_v))}\,,\,\,\,\,\,\, &|\alpha_1|\leq 2,\\
C\|h_1(t)\|^2_{L^2(\RR^3_x;\,\, L^2(\RR^3_v))} \|(W_l
h_2)(t)\|^2_{H^{3/2+\delta}(\RR^3_x;\,\,
H^{2s}_{(\gamma+2s)^+}(\RR^3_v))}\,,\,\,\,\,\,\, &|\alpha_1|>2,
\end{array}\right.\nonumber
\\
&\leq& C \|g^\varepsilon(t)\|^2_{H^k_{l}(\RR^6)}
\|g^\varepsilon(t)\|^2_{H^k_{l+(\gamma+2s)^+}(\RR^6)}, \qquad k\geq
4>3+2s,\,\,l>(\gamma+2s)^++3/2. \nonumber
\end{eqnarray*}
The estimation on the first term is similar to $(A_2, W_l^2 \tilde
h)_{L^2(\RR^6)}$ by taking into account the same manipulation
concerning $\a_2$.  The estimation for the second term is also
similar to the part $J_2$ of $I_2$ as above. Hence, we have obtained
\begin{equation}\label{E-3.105+8}
\left|\left(A_3,\,\, W^2_{l} \tilde h
\right)_{L^2(\RR^6)}\right|\leq C
\|g^\varepsilon(t)\|^2_{H^k_{l}(\RR^6)}
\|g^\varepsilon(t)\|_{H^k_{l+\gamma+2s}(\RR^6)}.
\end{equation}
This completes the proof of Proposition \ref{E-prop4.2}.\\

If (\ref{E-3.105+1}), \eqref{E-3.105+1-2}, \eqref{additional-term},
(\ref{E-3.105+4}) and (\ref{E-3.105+5}) are combined, then it follows
from \eqref{E-equ-2} that
\begin{align*}
\frac 1 2 \frac{d}{dt}\|P_{N,\, l}(\partial^\alpha
g^\varepsilon)(t)\|^2_{L^2(\RR^6)} &
 + \kappa
\|W_1\,P_{N,\, l}(\partial^\alpha g^\varepsilon)(t)\|^2_{L^2(\RR^6)} -
\frac{\kappa}{2}\|g^\varepsilon(t)\|^2_{H^k_{l+1}(\RR^6)}
\\
\le & C \left( \|g^\varepsilon(t)\|^2_{H^k_l(\RR^6)}+C_2
\|g^\varepsilon(t)\|^2_{H^k_l(\RR^6)}\|g^\varepsilon(t)
\|_{H^k_{l+(\gamma+2s+\delta)^+}(\RR^6)}  \right) .
\end{align*}
Take the sum over $|\a| \le k$, integrate from $0$ to $t \in [0,T]$
and make $N \rightarrow \infty$. Then there exists $C_1, C_2>0$
independent of $\varepsilon>0$ such that , for any $\delta>0$ and
$t\in [0, T]$,
\begin{align}\label{E-energy-estimate}
&\|g^\varepsilon(t)\|^2_{H^k_l(\RR^6)} + \kappa \int_0^t
\|g^\varepsilon(\tau )\|^2_{H^k_{l+1}(\RR^6)} d \tau
\,\,\,\,\,\,\,\,\,\,
\\
&\notag \le  \|g^\varepsilon(0)\|^2_{H^k_l(\RR^6)} + C_1 \int_0^t
\|g^\varepsilon(\tau )\|^2_{H^k_l(\RR^6)} d\tau  +C_2 \int_0^t
\|g^\varepsilon(\tau )\|^2_{H^k_l(\RR^6)}\|g^\varepsilon(\tau )
\|_{H^k_{l+(\gamma+2s+\delta)^+}(\RR^6)} d \tau .
\end{align}

\begin{rema}\label{comment-soft}
We give here some technical reasons about the choice of the time
dependent distribution $\mu(t)$ as moment control in the
equation (\ref{E-Cauchy-B}). If we take $\kappa=0$ in the definition
of Maxwellian distribution $\mu(t)$, the above computation gives also
(\ref{E-energy-estimate}) without the second term on the left hand
side because $\kappa=0$. But the upper bound estimate, by using
Theorem \ref{theo2.1}, always gives the last term in
(\ref{E-energy-estimate}) with the factor $\|g^\varepsilon(t)
\|_{H^k_{l+(\gamma+2s)^++\delta}(\RR^6)}$. If $\gamma+2s<0$, there
is no loss of moment, we can get (\ref{E-3.104+1}) with $\kappa=0$.
If $0\leq\gamma+2s<1$, we choose $\delta$ such that $\gamma + 2s +
\delta \leq 1$ so the second term on left hand side 
absorbs the last term in (\ref{E-energy-estimate}) because
$$
\|g^\varepsilon(t) \|_{H^k_{l+(\gamma+2s+\delta)^+}(\RR^6)}
\leq\|g^\varepsilon(t) \|_{H^k_{l+1}(\RR^6)}.
$$
In conclusion, the choice  of $\mu(t)$ is mainly for the hard
potential.
\end{rema}

\noindent {\bf Completion of proof of  Theorem \ref{E-uniform-estimate}. }
Set $X(t) = \|g^\varepsilon(t)\|^2_{H^k_l(\RR^6)} $ and $F(t) =
\int_0^t X(\tau)(1+X(\tau)) d\tau$. Since $\gamma+2s<1$, by
\eqref{E-energy-estimate} there exists a $C>0$ independent of
$\varepsilon>0$ such that
\begin{equation}\label{E-3.104+1}
X(t)  + \frac{\kappa}{2} \int_0^t \|g^\varepsilon(\tau
)\|^2_{H^k_{l+1}(\RR^6)} d \tau
             \leq X(0) + C F(t).
\end{equation}
Noticing that $F'(t) \leq \big (X(0) + CF(t) \big )\big (1+X(0)+CF(t) \big
)$, we have
$$
\|g^\varepsilon(t) \|^2_{H^k_l(\RR^6))}\leq
\frac{\|g_0\|^2_{H^k_l(\RR^6)}e^{{C} t}} {1-\big(e^{{C}
t}-1\big)\|g_0\|^2_{H^k_l(\RR^6)}}\,, \enskip 
$$
as long as the denominator remains  positive. 
We choose $T_*>0$ small enough such that
$$
\frac{e^{{C} T_*}} {1-\big(e^{{C}
T_*}-1\big)\|g_0\|^2_{H^k_l(\RR^6)}}=4.
$$
Then
$$
T_*=\frac{1}{C}\log\Big(1+\frac{3}{1+4\|g_0\|^2_{H^k_l(\RR^6)}}\Big),
$$
is independent of $\varepsilon>0$, but depends on
$\|g_0\|_{H^k_l(\RR^6)}$ and the constant $C$ which depends on
$\rho, \kappa, k$ and $l$. Now  we
have  (\ref{E-3.103}) for $T_{**} = \min(T,T_*)$.

From (\ref{E-3.103}) and (\ref{E-3.104+1}), we get also, for
$\kappa>0$,
\begin{equation}\label{L2-estimate}
\kappa  \|g^\varepsilon\|^2_{L^2(]0, T_{**}[;\, H^k_{l+1}(\RR^6))}
\leq 2 \|g_0\|^2_{H^k_l(\RR^6) } \Big ( 1+ 2C T_{*} (1 + 2
\|g_0\|^2_{H^k_l(\RR^6)}) \Big)  .
\end{equation}
We have proved Theorem \ref{E-uniform-estimate}.

\subsection{Convergence and uniqueness}\label{s4+++}
\setcounter{equation}{0}
\smallskip

\smallbreak

The second step is to prove that, for any $0<\varepsilon<1$,  we can
extend the approximation solution $g^\varepsilon$, obtained by
Theorem \ref{E-Cauchy-cut-off}, to a fixed interval $]0, T_*[$ with
$T_*>0$ determined in Theorem \ref{E-uniform-estimate} which is
independent on $\varepsilon>0$. Then this sequence is convergent.

\begin{theo}\label{E-unif-existence}
Assume that $0<s<1/2,\,\, \gamma+2s<1$,  $g_0\geq 0, \, g_0\in
H^k_l(\RR^6)$ for some $k\geq 4,\,\, l\geq 3$. Let $T_* >0$ be given
in Theorem \ref{E-uniform-estimate}. Then the Cauchy problem
(\ref{E-Cauchy-cut-off}) admits a unique non-negative solution up to
$T_*$ satisfying
$$
g^\varepsilon\in L^\infty(]0, T_*[;\,\, H^k_l(\RR^6))\cap L^2(]0,
T_*[;\,\, H^k_{l+1}(\RR^6)).
$$
\end{theo}

\smallskip
\noindent {\bf Proof:} We recall the notation $T=T_\varepsilon(D_0)$
from Remark \ref{Tepsilon}. Then Theorem \ref{E-exist-cut-off}
asserts that the Cauchy problem (\ref{E-Cauchy-cut-off}) with
initial data $g_0$ admits a unique
 non-negative solution
$$
g^\varepsilon_1\in C^0([0,\,\, 2T_{1, \varepsilon}];\,\,
H^k_l(\RR^6))\cap L^2(]0, 2T_{1, \varepsilon}[;\,\,
H^k_{l+1}(\RR^6)), \quad T_{1,\varepsilon}=\frac 12
T_\varepsilon(\|g_0\|_{H^k_l(\RR^6)}).
$$
If $T_{1,\varepsilon} \geq  T_*$, then the proof is completed. If $T_{1,
\varepsilon}<T_*$, then Theorem \ref{E-uniform-estimate} implies
$$
\|g^\varepsilon_1(T_{1, \varepsilon})\|_{H^k_l(\RR^6)}\leq 2
\,\|g_0\|_{H^k_l(\RR^6)}.
$$
We now consider the Cauchy problem (\ref{E-Cauchy-cut-off}) with
initial data $g^\varepsilon(T_{1, \varepsilon})$. Again Theorem
\ref{E-exist-cut-off} asserts that there exists
\[
T_{2,\varepsilon}=\frac 12 T_\varepsilon(2\|g_0\|_{H^k_l(\RR^6)}),
\]
such that the Cauchy problem (\ref{E-Cauchy-cut-off}) admits a
unique non-negative solution
$$
g^\varepsilon_2\in C^0([T_{1, \varepsilon},\,\, T_{1,
\varepsilon}+2T_{2, \varepsilon}];\,\, H^k_l(\RR^6))\bigcap
L^2(]T_{1, \varepsilon}, T_{1, \varepsilon}+2T_{2,
\varepsilon}[;\,\, H^k_{l+1}(\RR^6)).
$$
By uniqueness of solution, we obtain a  non-negative solution of the
Cauchy problem (\ref{E-Cauchy-cut-off}),
$$
g^\varepsilon\in C^0([0,\,\, T_{1, \varepsilon}+2T_{2,
\varepsilon}];\,\, H^k_l(\RR^6))\bigcap L^2(]0, T_{1,
\varepsilon}+2T_{2, \varepsilon}[;\,\, H^k_{l+1}(\RR^6)).
$$
If $T_{1, \varepsilon}+2T_{2, \varepsilon}\geq T_*$, we finish the
proof. If $T_{1, \varepsilon}+2T_{2, \varepsilon}< T_*$, we consider
again the Cauchy problem (\ref{E-Cauchy-cut-off}) with initial date
$g^\varepsilon(T_{1, \varepsilon}+T_{2, \varepsilon})$. Since
Theorem \ref{E-uniform-estimate} gives again
$$
\|g^\varepsilon(T_{1, \varepsilon}+T_{2,
\varepsilon})\|_{H^k_l(\RR^6)}\leq 2 \,\|g_0\|_{H^k_l(\RR^6)},
$$
the interval of the existence of solution is the same, that
is, $2T_{2,\varepsilon}$,
 so that we can
extend the solution to
$$
g^\varepsilon\in L^\infty(]0,\,\, T_{1, \varepsilon}+3T_{2,
\varepsilon}[;\,\, H^k_l(\RR^6))\bigcap L^2(]0, T_{1,
\varepsilon}+3T_{2, \varepsilon}[;\,\, H^k_{l+1}(\RR^6)).
$$
By iteration, there exists $m\in\NN$ such that
$$
T_{1, \varepsilon}+ m T_{2, \varepsilon}<T_*,\,\,T_{1, \varepsilon}+
(m+1) T_{2, \varepsilon}\geq T_*,
$$
and we extend the solution up to
$$
g^\varepsilon\in C^0([0,\,\, T_{1, \varepsilon}+(m+1)T_{2,
\varepsilon}];\,\, H^k_l(\RR^6))\bigcap L^2(]0, T_{1,
\varepsilon}+(m+1)T_{2, \varepsilon}[;\,\, H^k_{l+1}(\RR^6)).
$$
We have proved Theorem \ref{E-unif-existence}.\\

\bigbreak Theorem \ref{E-unif-existence} asserts the existence of an
approximation solution sequence
$$
\Big\{g^\varepsilon\Big\}_{\varepsilon>0}\subset C^0([0,\,\,
T_*];\,\, H^k_l(\RR^6))\bigcap L^2(]0, T_*;\,\, H^k_{l+1}(\RR^6)),
$$
and
$$
\|g^\varepsilon\|_{L^\infty (]0, T_*[\,;\,\,H^k_l(\RR^6))}\leq 2
\,\|g_0\|_{H^k_l(\RR^6)}.
$$
This implies that it is a weakly* compact set of $L^\infty(]0,\,\,
T_*[;\,\, H^k_l(\RR^6))$. Let
\[
g\in L^\infty (]0, T_*[\,;\,\,H^k_l(\RR^6)),
\]
be a limit of a subsequence of
$\big\{g^\varepsilon\big\}_{\varepsilon>0}$.

On the other hand, by using the equation (\ref{E-Cauchy-cut-off}) and
Theorem \ref{theo2.1}, we obtain
\begin{eqnarray*}
\|\partial_t g^\varepsilon\|_{L^\infty (]0,
T_*[\,;\,\,H^{k-1}_{l-1}(\RR^6))}&\leq&
C\Big(\|g^\varepsilon\|_{L^\infty (]0, T_*[\,;\,\,H^k_l(\RR^6))}+\,
\|g^\varepsilon\|^2_{L^\infty
(]0, T_*[\,;\,\,H^k_l(\RR^6))}\Big)\\
&\leq& 2C \,(1+2\|g_0\|_{H^k_l(\RR^6)})\|g_0\|_{H^k_l(\RR^6)}.
\end{eqnarray*}
Thus, $\big\{g^\varepsilon\big\}_{\varepsilon>0}$ is a compact subset
in
\[
C^{1-\delta} (]0, T_*[\,;\,\,H^{k-1-\delta}_{l-1}(\Omega
\times\RR^3_v)),
\]
for any compact bounded open set $\Omega \subset\RR_x^3$ and for any
$\delta>0$. For the variable $v$, we have the weight $W_{l-1}$ with
$l-1>3/2$. Then, we can take the limit in the equation
(\ref{E-Cauchy-cut-off})
 and also in the mild form \eqref{mild}. Then $g$ is  a solution
 of the Cauchy problem (\ref{E-Cauchy-B}). The limit $g$ belongs to
 $ L^2(]0, T_*[;\,\, H^k_{l+1}(\RR^6))$  deduced from \eqref{L2-estimate}.
Now if $g_0\geq 0$, Theorem \ref{E-exist-cut-off} implies that
$g^\varepsilon\geq 0$, so that the limit $g$ is also non-negative on
$]0, T_*[$. We have completed the proof for the local existence of  solutions
stated in Theorem \ref{E-theo-0.2}.

\smallskip
It remains to prove the uniqueness of solutions in Theorem
\ref{E-theo-0.2}. We state it more precisely
 as follows.

\begin{prop}\label{prop-unique}
Assume that $0<s<1/2,\,\, \gamma+2s<1,\, 0<T\leq T_0,\, m>3$ and
$g_0\geq 0, \, g_0\in H^m_3(\RR^6)$. Suppose that the Cauchy problem
\eqref{E-Cauchy-B} admits two (non-negative) solutions
$$
g_1,\, g_2\in C^0([0, T];\,\, H^m_{4}(\RR^6)).
$$
Then $g_1\,\equiv\, g_2$.
\end{prop}

Set $f=g_1-g_2$, by using \eqref{E-Cauchy-B}, we have
\begin{equation}\label{E-Cauchy-U}
\left\{\begin{array}{l}
f_t+v\cdot\nabla_x f\  +\kappa (1+ |v|^2) f=\Gamma^t(g_1,\, f)+\Gamma^t(f,\ g_2)\, ,\\
f|_{t=0}=0.
\end{array}\right.
\end{equation}
We can now take  $W_3 f$ as a test function to get
\begin{equation}\label{2.4.1}
\frac 12 \frac{d}{d t} \|W_3 f(t)\|^2_{L^2(\RR^6)}+\kappa\|W_{4}
f(t)\|^2_{L^2(\RR^6)}=\Big(W_3\Gamma^t(g_1,\, f)+W_3\Gamma^t(f,\
g_2)\, ,W_3 f\Big)_{L^2(\RR^6)}.
\end{equation}
Recall that
$$
\Gamma^t(g,\, h)=Q(\mu(t) g,\, h)+\int_{\RR^3_{v_*}\times\SS^2} B
\big(\mu(t)_*-\mu(t)'_*\big) g'_* h' d v_* d\sigma.
$$

We estimate the last two terms of \eqref{2.4.1} in the following
lemma.

\begin{lemm}\label{lemm2.4.1}
Assume that $g_1\geq 0$. Then for any $\varepsilon>0$, there exist constants
$C_\varepsilon>0$ and  $K(\varepsilon, \|g_2\|_{L^\infty(]0,\, T[;
H^m_{4}(\RR^6))})>0$ such that
\begin{equation}\label{2.4.2}
\Big(W_3\Gamma^t(g_1,\, f)\, ,W_3 f\Big)_{L^2(\RR^6)}\leq \varepsilon
\|W_{4} f(t)\|^2_{L^2(\RR^6)}+C_\varepsilon\|g_1\|^2_{L^\infty(]0,\,
T[; H^m_{4}(\RR^6))}
 \|W_{3}f(t)\|^2_{L^2(\RR^6)},
\end{equation}
\begin{equation}\label{2.4.3}
\left|\Big(W_3\Gamma^t(f,\, g_2)\, ,W_3
f\Big)_{L^2(\RR^6)}\right|\leq \varepsilon \|W_{4}
f(t)\|^2_{L^2(\RR^6)}+K(\varepsilon, \|g_2\|_{L^\infty(]0,\, T[;
H^m_{3}(\RR^6))})
 \|W_{3}f(t)\|^2_{L^2(\RR^6)}.
\end{equation}
\end{lemm}

Notice that by using the above lemma with $\varepsilon=\kappa/4$ and
(\ref{2.4.1}), we get
$$
\frac{d}{d t} \|W_3 f(t)\|^2_{L^2(\RR^6)}\leq
\Big(C\|g_1\|^2_{L^\infty(]0,\, T[; H^m_{4}(\RR^6))} +K(\varepsilon,
\|g_2\|_{L^\infty(]0,\, T[; H^m_{4}(\RR^6))})\Big)
 \|W_{3}f(t)\|^2_{L^2(\RR^6)}.
$$
Then $\|W_{3}f(0)\|_{L^2(\RR^6)}=0$ implies
$\|W_{3}f(t)\|_{L^2(\RR^6)}=0$ for all $0\leq t\leq T$
which gives  Proposition \ref{prop-unique}.

\bigbreak \noindent{\bf Proof of Lemma \ref{lemm2.4.1}.} As for
\eqref{2.4.2}, we have
\begin{eqnarray*}
&&\Big(W_3\Gamma^t(g_1,\, f)\, ,W_3 f\Big)_{L^2(\RR^6)}\\
&=& \Big(W_3Q(\mu(t)g_1,\, f)\, ,W_3 f\Big)_{L^2(\RR^6)}+\iiiint
B\,\,(\mu(t)_*-\mu(t)'_*)
g'_{1 *} f' W^2_3 f d v_* d\sigma dv dx \\
&=&\Big(Q(\mu(t)g_1,\, W_3 f)\, ,W_3 f\Big)_{L^2(\RR^6)}+
\Big(W_3 Q(\mu(t)g_1,\, f)-Q(\mu(t)g_1,\, W_3 f)\,\, ,W_3 f\Big)_{L^2(\RR^6)}\\
&&+\iiiint B\,\,(\mu(t)_*-\mu(t)'_*)g'_{1 *} \big(W_3 f\big)' W_3 f d v_* d\sigma dv dx\\
&&+\iiiint B\,\,(\mu(t)_*-\mu(t)'_*)
g'_{1 *}\big(W_3-W'_3 \big) f' W_3 f d v_* d\sigma dv dx \\
&=&D_1+D_2+D_3+D_4.
\end{eqnarray*}
The term $D_1$ is similar to $B_1$ in the proof of Proposition
\ref{E-prop4.1}. By using $\mu(t)g_1\geq 0$, we have
$$
D_1\leq C  \|g_1(t)\|_{H^{3/2+\delta}(\RR^6_{x, \, v})}
\|f(t)\|_{L^2_{3}(\RR^6_{x, \,
v})}\|f(t)\|_{L^2_{3+\gamma^+}(\RR^6_{x, \, v})},
$$
for some small $\delta>0$. The term $D_2$ is similar to $B_2$ and
we can obtain
$$
|D_2|\leq C  \|g_1(t)\|_{H^{3/2+\delta}(\RR^6_{x, \, v})}
\|f(t)\|_{L^2_{3}(\RR^6_{x, \,
v})}\|f(t)\|_{L^2_{3+\gamma^+}(\RR^6_{x, \, v})}.
$$
The terms $D_3, D_4$ are similar to $I_1, I_2$ in the proof of
Proposition \ref{E-prop4.2}. Namely
$$
|D_3|+|D_4|\leq  C
\|g_1(t)\|_{H^{6/2+\delta}_{3+(\gamma+2s+\delta)^+} (\RR^6_{x, \,
v})} \|f(t)\|_{L^2_{3}(\RR^6_{x, \, v})}\|f(t)\|_{L^2_{3+
(\gamma+2s+\delta)^+}(\RR^6_{x, \, v})}.
$$
Thus, for any $0<t\leq T$ and $m>3$, we have
$$
\Big(W_3\Gamma^t(g_1,\, f)\, ,W_3 f\Big)_{L^2(\RR^6)}\leq C
\|g_1\|_{L^\infty(]0, T[; H^{m}_{4}(\RR^6_{x, \, v}))} \|W_3
f(t)\|_{L^2(\RR^6_{x, \, v})}\|W_4 f(t)\|_{L^2(\RR^6_{x, \, v})},
$$
which implies \eqref{2.4.2}. The left hand side of  \eqref{2.4.3} can be 
written as
\begin{eqnarray*}
&&\Big(W_3\,\Gamma^t(f,\, g_2)\, ,W_3\, f\Big)_{L^2(\RR^6)}\\
&=& \Big(W_3 Q(\mu(t)f,\, g_2)\, ,W_3 f\Big)_{L^2(\RR^6)}+\iiiint
B\,\,(\mu(t)_*-\mu(t)'_*)
f'_{*} g'_2 W^2_3 f d v_* d\sigma dv dx \\
&=&\Big(W_3 Q(\mu(t)f,\, g_2)\, ,W_3 f\Big)_{L^2(\RR^6)}+
\iiiint B\,\,(\mu(t)_*-\mu(t)'_*)f'_{ *} \big(W_3 g_2\big)' W_3 f d v_* d\sigma dv dx\\
&&+\iiiint B\,\,(\mu(t)_*-\mu(t)'_*)
f'_{*}\big(W_3-W'_3 \big) g'_2 W_3 f d v_* d\sigma dv dx \\
&=&E_1+E_2+E_3.
\end{eqnarray*}
Using Corollary \ref{lemm104} with $m=0, l=3$ gives
\begin{eqnarray*}
&&|E_1|\leq \int_{\RR^3_x} \|W_3 Q(\mu(t)f,\,
g_2)\|_{L^2(\RR^3_v)}\|W_3 f\|_{L^2(\RR^3_v)}
dx\\
&&\leq C\int_{\RR^3_x} \|\mu(t)f\|_{L^1_{3+(\gamma+2s)^+}(\RR^3_v)}
\|g_2\|_{H^{2s}_{3+(\gamma+2s)^+}(\RR^3_v)}\|W_3 f\|_{L^2(\RR^3_v)}
dx\\
&&\leq C \|g_2\|_{L^\infty(]0, T[\times\RR^3_x;\, H^{2s}_{3+
(\gamma+2s+\delta)^+}(\RR^3_v))} \|f(t)\|_{L^2(\RR^6)}
\|W_3 f(t)\|_{L^2(\RR^6)}\\
&&\leq C \|g_2\|_{L^\infty(]0, T[;\,
H^{3/2+2s+\delta}_{3+(\gamma+2s+\delta)^+}(\RR^3_v))} \|W_3
f(t)\|^2_{L^2(\RR^6)}.
\end{eqnarray*}
The term $E_2$ is similar to $D_3$, and we have
\begin{eqnarray*}
&&|E_2|\leq C \|f(t)\|_{L^2 (\RR^3_x;\, L^{1}_{(\gamma+2s+
\delta)^+}(\RR^3_v))} \|g_2\|_{L^\infty(]0, T[\times\RR^3_x;\,
L^2_{3+\gamma^+}(\RR^3_v))}
\|W_3 f(t)\|_{L^2(\RR^6)}\\
&&\leq C \|f(t)\|_{L^2_{3/2+\delta+ (\gamma+2s+\delta)^+}(\RR^6))}
\|g_2\|_{L^\infty(]0, T[;\, H^{3/2+\delta}_{3+\gamma^+}(\RR^6)) }
\|W_3 f(t)\|_{L^2(\RR^6)}\\
&&\leq C \|g_2\|_{L^\infty(]0, T[;\,
H^{3/2+\delta}_{4}(\RR^6)) }\|W_3f(t)\|_{L^2(\RR^6))}^2. 
\end{eqnarray*}
For the term $E_3$, we
can use \eqref{E-3.105+3+1} with $l=3$. Then
\begin{eqnarray*}
&&|E_3|\leq \iiiint b(\cos\theta)\,\,\left\langle v-v_*
\right\rangle^\gamma |\mu(t)_*-\mu(t)'_*|\,
|f'_{*}|\,\,\big|W_3 - W'_3 \big|\,\,| g'_2|\, |W_3 f|\, d v_* d\sigma dv dx \\
&\leq& C\iiiint \sin \Big(\frac{\theta}{2}\Big) b(\cos\theta)\,\,
|(W_{1+\gamma^+}f)'_{*}|\,| (W_{3+\gamma^+} g_2)'|\, |W_3 f|\, d v_* d\sigma dv dx \\
&+& C\iiiint \sin^{3}\Big(\frac{\theta}{2}\Big)
b(\cos\theta)\,\,\mu'_*(t)
|(W_{3+\gamma^+} f)'_{*}|\,| (W_{\gamma^+}g_2)'|\, |W_3 f|\, d v_* d\sigma dv dx \\
&+& C\iiiint \sin^{3}\Big(\frac{\theta}{2}\Big)
b(\cos\theta)\,\,\mu_*(t)
|(W_{3+\gamma^+} f)'_{*}|\,|(W_{\gamma^+}g_2)'|\, |W_3 f|\, d v_* d\sigma dv dx\\
&=&E_{3, 1}+E_{3, 2}+E_{3, 3}.
\end{eqnarray*}
Since  $0<2s<1$ is assumed, for any $\varepsilon>0$ there exists
$C_\varepsilon>0$ such that
\begin{eqnarray*}
&&|E_{3, 1}|\leq C\int_{\RR^3_x} \|f(t,\, x,\,
\cdot)\|_{L^1_{1+\gamma^+}(\RR^3_v)}\,\|g_2(t,\, x,\, \cdot)
\|_{L^2_4(\RR^3_v)}\, \| f(t,\, x,\, \cdot)\|_{L^2_{3}(\RR^3_v)}\, dx \\
&&\leq C\|f(t)\|_{L^2(\RR^3_x;\,
L^1_{1+\gamma^+}(\RR^3_v))}\,\|g_2\|_{L^\infty(]0, T[
\times\RR^3_x;\, L^2_4(\RR^3_v))}\, \| f(t)\|_{L^2_{3}(\RR^6_{x,
v})}
\\
&&\leq C\|f(t)\|_{L^2_{3/2+\delta+1+\gamma^+}(\RR^6_{x,
v})}\,\|g_2\|_{L^\infty(]0, T[;\,
H^{3/2+\delta}_4(\RR^6_{x, v}))}\, \| f(t)\|_{L^2_{3}(\RR^6_{x, v})}\\
&&\leq \left(\varepsilon\|W_4 f(t)\|^2_{L^2(\RR^6_{x, v})}+ C_\varepsilon
\|W_3 f(t)\|^2_{L^2(\RR^6_{x, v})}\right)\|g_2\|_{L^\infty(]0, T[;\,
H^{3/2+\delta}_4(\RR^6_{x, v}))}\,.
\end{eqnarray*}
Similarly
\begin{eqnarray*}
&&|E_{3, 2}|\leq C\int_{\RR^3_x} \|\mu(t) f(t,\, x,\,
\cdot)\|_{L^1_{1+\gamma^+}(\RR^3_v)}\,\|g_2(t,\, x,\, \cdot)
\|_{L^2_{\gamma^+}(\RR^3_v)}\, \| f(t,\, x,\, \cdot)\|_{L^2_{3}(\RR^3_v)}\, dx \\
&&\leq C \|g_2\|_{L^\infty(]0, T[;\, H^{3/2+\delta}_4(\RR^6_{x,
v}))}\| f(t) \|_{L^2(\RR^6_{x, v})}\|W_3 f(t)\|_{L^2(\RR^6_{x,
v})}\,.
\end{eqnarray*}
Since $3/2+(3+\gamma^+)>4$, we can not estimate $E_{3, 3}$ in the same way
as for
 $E_{3, 2}$. Instead, we have
\begin{eqnarray*}
|E_{3, 3}|&\leq& C\|W_{\gamma^+}g_2\|_{L^\infty(]0,
T[\times\RR^6_{x, v})} \iiiint \theta^{3} b(\cos\theta)\,\,\mu_*(t)
|(W_{3+\gamma^+}f)'_{*}|\,|\, |W_3 f|\, d v_* d\sigma dv dx\\
&&\leq C \|g_2\|_{L^\infty(]0, T[;\, H^{3+\delta}_3(\RR^6_{x,
v}))}\int_{\RR^3_x} \left(\iiint \theta^{1}
b(\cos\theta)\,\,\mu_*(t)
|W_3f|^2\,|\, d v_* d\sigma dv\right)^{\frac 12}\\
&& \times \left(\iiint \theta^{5} b(\cos\theta)\,\,\mu_*(t)
|(W_{3+\gamma^+}f)'_{*}|^2\,|\, d v_* d\sigma dv\right)^{\frac
12}dx\,.
\end{eqnarray*}
We now take the singular change of variables $v'_*\,\rightarrow\, v$.
The Jacobian is computed in \eqref{jacobian2} which is of the order of
$ \theta^{-2}$. Then this singular change of variables yields
\begin{eqnarray*}
&&\iiint \theta^{5} b(\cos\theta)\,\,\mu_*(t)
|(W_{3+\gamma^+}f)'_{*}|^2\,|\, d v_*  d \sigma dv\\
&\leq&C \iint D_1(v_*, v'_*) \,\,\mu_*(t) |(W_{3+\gamma^+}
f)'_{*}|^2\,|\, d v_* dv'_*,
\end{eqnarray*}
with $ D_1(v_*,v'_*)=\int_{S^2}\theta^{5-2}b(\cos\theta)d\sigma \le
C \int_{\pi/4}^{\pi/2}(\frac{\pi}{2}-\psi)^{-2-2s+5-2}d\psi \le C $.
Hence
\begin{eqnarray*}
&&\iiint \theta^{5} b(\cos\theta)\,\,\mu_*(t)
|(W_{3+\gamma^+}f)'_{*}|^2\,|\, d v_*  d \sigma dv\\
&\leq& C \| \mu(t)\|_{L^1(\RR^3_{v})}\|W_{3+\gamma^+}
f(t,\, x,\, \cdot )\|^2_{L^2(\RR^3_{ v})}.
\end{eqnarray*}
Therefore,
$$
|E_{3, 3}|\leq  C \|g_2\|_{L^\infty(]0, T[;\,
H^{3+\delta}_4(\RR^6_{x, v}))} \|W_{3}f\|_{L^2(\RR^6_{x,
v})}\|W_{3+\gamma^+}f\|_{L^2(\RR^6_{x, v})}.
$$
By combining the estimates on $E_1, E_2, E_3$, we have proved \eqref{2.4.3}.
Now the proof of Lemma \ref{lemm2.4.1} is complete.

\subsection{Proof of Theorem \ref{E-theo1}}
\setcounter{equation}{0}
\smallskip
Assume that $f_0\in
\cE^{k_0}_0(\RR^5)$. Then there exists $\rho_0>0$ such that $
e^{\rho_0\la v\ra^2}\, f_0\in H^{k_0}(\RR^6)$.
Choose $0<\rho<\rho_0$ and $\kappa>0$ small enough. By setting
$g_0=e^{\,\rho\la v\ra^2}f_0$, then $ g_0\in H^{k_0}_l(\RR^6)$ for all
$l\in\NN$. Theorem \ref{E-theo-0.2} asserts that the Cauchy problem
\eqref{E-Cauchy-B} with the initial datum $g_0$ admits a non-negative local solution
$$
g\in C^0 ([0, T_*]; H^{k_0}_{l}(\RR^6))\bigcap L^2 (]0, T_*[;
H^{k_0}_{l+1}(\RR^6)),\,\,\,\,\,\,\,\,\forall\,\, l\in\NN,
$$
with $T_* \in ]0,T_0]$  \,\, ($T_0=\frac{\rho}{2\kappa}$). Then
$$
f(t, x, v)=e^{-(\rho-\kappa t)\la v\ra^2}g(t, x, v)\in C^0 ([0,
T_*]; H^{k_0}_{l}(\RR^6))\bigcap L^2 (]0, T_*[;
H^{k_0}_{l}(\RR^6)),\,\,\,\,\,\,\,\,\forall\,\, l\in\NN,
$$
is a non-negative solution of the Cauchy problem \eqref{1.1b}. Since
for $0\leq t\leq T_*\leq T_0$,
\begin{equation}\label{exp-dec}
e^{\frac{\rho}{2}\la v\ra^2}\, f \in C^0([0, T_*];
H^{k_0}(\RR^6)),
\end{equation}
we can conclude $f\in \cE^{k_0}([0,T_*] \times \RR^6_{x, v})$, which leads to
 the local existence stated in Theorem \ref{E-theo1}.

Suppose now for some $f_0\in \cE^{4}_0(\RR^5)$, the Cauchy problem
\eqref{1.1b} admits two solutions $f_1 \in \cE^{4}([0,T_1] \times \RR^6_{x, v})$ and 
$f_2 \in \cE^{4}([0,T_2] \times \RR^6_{x, v})$. This implies that
 there exist
$\rho_0, \rho_1,  \rho_2 >0$  such that
$$
e^{\rho_0\la v\ra ^2} f_0\in H^{4}(\RR^6),
$$
and
$$
e^{\rho_1\la v\ra ^2} f_1\in C^0 ([0, T_1];\,
H^{4}(\RR^6)),\,\,\,\,\,\,\, e^{\rho_2\la v\ra ^2} f_2\in C^0
([0, T_2];\, H^{4}(\RR^6)).
$$
Take  $0<\rho<\min\{\rho_0,\, \rho_1,\, \rho_2\}$ and $\kappa>0$
sufficiently small such that $\frac{\rho}{2\kappa}>T_{**}=\min\{T_1, T_2\}$.
Then we have
$$
g_0=e^{\rho\la v\ra ^2} f_0\in H^{4}_l(\RR^6),
$$
for any $l\in\NN$, and
$$
g_1=e^{(\rho-\kappa t)\la v\ra ^2} f_1\in C^0 ([0, T_{**}];\,
H^{4}_l(\RR^6)),\,\,\,\,\,\,\, g_2= e^{(\rho-\kappa t)\la v\ra ^2}
f_2\in C^0 ([0, T_{**}];\, H^{4}_l(\RR^6)),
$$
are two solutions of the Cauchy problem \eqref{E-Cauchy-B} with the common initial
datum $g_0$. Then Proposition \ref{prop-unique} gives $g_1=g_2$, 
so that $f_1=f_2$ for $t\in [0, T_{**}]$.  Now the uniqueness of
solutions stated in Theorem \ref{E-theo1} is obvious since 
$T_1 =T_2 = T_{**}$.

On the other hand, in view of \eqref{exp-dec},
$\|f(t,\, x,\, \cdot\,)\|_{L^1}$ is continuous for $(t, x)\in
[0, T_*]\times\RR^3_x$. Therefore,  if for a compact $K\subset\RR^3$, we have
$$
\inf_{x\in K} \|f_0( x,\, \cdot\,)\|_{L^1}=c_0>0,
$$
then there exist $0<\tilde{T}_0\leq T_*$ and  a closed neighborhood of
$K$ denoted by $V_0$ in $\RR^3_x$ such that
$$
\inf_{(t, x)\in [0,\tilde{T}_0]\times V_0} \|f(t,\, x,\,
\cdot\,)\|_{L^1}\geq \frac{c_0}{2}.
$$
Now Theorem \ref{theo1} implies that
$$
f\in \bigcap_{l\in\NN}\cH^{+\infty}_l( ]0,\tilde{T}_0[\times
V_0\times\RR^3_v)\subset C^\infty( ]0,\tilde{T}_0[\times V_0;
\cS(\RR^3_v)).
$$

It remains to prove the uniqueness of solutions of Theorem
\ref{E-theo1} in the soft potential case $\gamma\leq 0$. In this
case, the uniqueness of solution can be proved
 in a  larger functional
space. We state it as follows.
\begin{prop}\label{Unique}
Assume that $0<s<1/2,\,\, \gamma\leq 0,\, 0<T\leq +\infty$ and $
m>2s+3/2, \, l>2s+3/2$. Let $f_0\geq 0, \, f_0\in
H^m_{l+2s}(\RR^6)$. Suppose that the Cauchy problem \eqref{1.1b}
admits two non-negative solutions
$$
f_1,\, f_2\in L^\infty(]0, T[;\,\, H^m_{l+2s}(\RR^6)).
$$
Then $f_1\,\equiv\, f_2$.
\end{prop}

\noindent{\bf Proof:} The proof is similar to the one for Proposition
\ref{prop-unique}. Set $F=f_1-f_2$, by using \eqref{1.1b}, we have
\begin{equation}\label{E-Cauchy}
\left\{\begin{array}{l}
F_t+v\cdot\nabla_x F =Q(f_1,\, F)+Q(F,\, f_2)\, ,\\
F|_{t=0}=0.
\end{array}\right.
\end{equation}
We can now take $W_l F$ as a test function to have
\begin{equation}\label{4.4.1}
\frac 12 \frac{d}{d t} \| F(t)\|^2_{L^2_l(\RR^6)}=\Big(W_l\,
Q(f_1,\, F)+W_l\, Q(F,\, f_2)\, ,W_l F\Big)_{L^2(\RR^6)}.
\end{equation}
Since $f_1\geq 0$ and $\gamma\leq 0$, similar to the
analysis on  $B_1$ in the
proof of Proposition \ref{E-prop4.1}, we have
$$
\Big(Q(f_1,\, W_l F)\, ,W_l F\Big)_{L^2(\RR^6)}\leq C
\|f_1(t)\|_{L^\infty(\RR^3_x; \,L^1(\RR^3_{v}))}
\|F(t)\|^2_{L^2_{l}(\RR^6_{x, \, v})}.
$$
Using \eqref{3.11-+} with $\gamma^+=0$ gives
$$
\left|\Big(W_l Q(f_1,\, F)-Q(f_1,\, W_l F)\,\, ,W_l
F\Big)_{L^2(\RR^6)}\right|\leq C  \|f_1(t)\|_{L^\infty(\RR^3_x;
\,L^2_l(\RR^3_{v}))}\|F(t)\|^2_{L^2_{l}(\RR^6_{x, \, v})},
$$
and
$$
\left|\Big(W_l Q(F,\, f_2)-Q(F,\, W_l f_2)\,\, ,W_l
F\Big)_{L^2(\RR^6)}\right|\leq C  \|F(t)\|_{L^2_l(\RR^6_{x, \, v})}
\|f_2(t)\|_{L^\infty(\RR^3_x;
\,L^2_l(\RR^3_{v}))}\|F(t)\|_{L^2_{l}(\RR^6_{x, \, v})}.
$$
Finally, for $l>3/2+2s$, we have
\begin{eqnarray*}
&&\left|\Big(Q(F,\, W_l f_2)\, ,W_l F\Big)_{L^2(\RR^6)}\right|\leq C
\|Q(F,\, W_l f_2)\|_{L^2(\RR^6)}\|F(t)\|_{L^2_{l}(\RR^6)}\\
&& \leq \|F(t)\|_{L^2_{l}(\RR^6)} \left(\int_{\RR^3_x}\|F(t, x,
\cdot\,)\|^2_{L^{1}_{2s}(\RR^3_{v})} \|f_2(t, x, \cdot\,
)\|^2_{H^{2s}_{l+2s}(\RR^3_{v})}\right)^{1/2}\\
&&\leq C \|F(t)\|^2_{L^2_{l}(\RR^6)}\|f_2(t)\|_{L^\infty(\RR^3_x;\;
H^{2s}_{l+2s}(\RR^3_{v}))}.
\end{eqnarray*}
Thus, we have, for any $0<t< T$ and $\delta>0$ small enough,
$$
\frac{d}{d t} \| F(t)\|^2_{L^2_l(\RR^6)}\leq C \left(
\|f_1\|_{L^\infty(]0, T[;\; H^{3/2+\delta}_{l}(\RR^6_{x, v}))}+
\|f_2\|_{L^\infty(]0, T[;\; H^{3/2+\delta+2s}_{l+2s}(\RR^6_{x,
v}))}\right)\|F(t)\|^2_{L^2_{l}(\RR^6)}.
$$
Therefore, $\| F(0)\|_{L^2_l(\RR^6)}=0$ implies $\|
F(t)\|_{L^2_l(\RR^6)}=0$ for all $t\in [0, T[$.

\bigskip
{\bf Acknowledgements:} The authors would like to express their sincere thanks to the referee
for his valuable comments. The research of the second author was
supported by  Grant-in-Aid for Scientific Research No.18540213,
Japan Society of the Promotion of Science. The last author's
research was supported by the General Research Fund of Hong Kong, CityU\#102606.


\end{document}